\documentclass[10pt]{article}
\usepackage{amsthm,amsfonts,amsmath,amscd,amssymb,latexsym,epsfig}
\usepackage[all]{xy}
\newtheorem{theorem}{Theorem}[section]
\newtheorem{corollary}[theorem]{Corollary}

\newtheorem{lemma}[theorem]{Lemma}
\newtheorem{proposition}[theorem]{Proposition}
\newtheorem{conjecture}[theorem]{Conjecture}

\newtheorem{definition}[theorem]{Definition}
\newtheorem{remark}[theorem]{Remark}

\newtheorem{example}[theorem]{Example}
\newcommand{\one}{{{\mathchoice {\rm 1\mskip-4mu l} {\rm 1\mskip-4mu l}
{\rm 1\mskip-4.5mu l} {\rm 1\mskip-5mu l}}}}
\newcommand{\dslash}{/\mskip-6mu/}
\newcommand{\A}{{\mathbb{A}}}
\newcommand{\B}{{\mathbb{B}}}
\newcommand{\C}{{\mathbb{C}}}
\newcommand{\E}{{\mathbb{E}}}

\newcommand{\HH}{{\mathbb{H}}}
\newcommand{\LL}{{\mathbb{L}}}

\newcommand{\Q}{{\mathbb{Q}}}
\newcommand{\R}{{\mathbb{R}}}
\renewcommand{\SS}{{\mathbb{S}}}
\newcommand{\T}{{\mathbb{T}}}
\newcommand{\Z}{{\mathbb{Z}}}
\newcommand{\Aa}{{\mathcal{A}}}   
\newcommand{\Bb}{{\mathcal{B}}}
\newcommand{\Dd}{{\mathcal{D}}}
\newcommand{\Ee}{{\mathcal{E}}}
\newcommand{\Ff}{{\mathcal{F}}}
\newcommand{\Gg}{{\mathcal{G}}}   
\newcommand{\Hh}{{\mathcal{H}}}

\newcommand{\Jj}{{\mathcal{J}}}

\newcommand{\Ll}{{\mathcal{L}}}   
\newcommand{\Mm}{{\mathcal{M}}}   

\newcommand{\Pp}{{\mathcal{P}}}

\newcommand{\Ss}{{\mathcal{S}}}
\newcommand{\Tt}{{\mathcal{T}}}
\newcommand{\Uu}{{\mathcal{U}}}

\newcommand{\Ww}{{\mathcal{W}}}
\newcommand{\Xx}{{\mathcal{X}}}
\newcommand{\Yy}{{\mathcal{Y}}}

\newcommand{\coker}{{\rm coker }}  
\newcommand{\im}{{\rm im }}        
\newcommand{\SPAN}{{\rm span }}    
\newcommand{\trace}{{\rm trace }}  
\newcommand{\id}{{\rm id}}         

\newcommand{\INDEX}{{\rm index}}   
\newcommand{\IND}{{\mathcal{IND}}} 
\newcommand{\Lie}{{\rm Lie}}          
\newcommand{\Vect}{{\rm Vect}}        
\newcommand{\Vol}{{\rm Vol}}          
\newcommand{\Sym}{{\rm Sym}}          
\newcommand{\Spin}{{\rm Spin}}        
\newcommand{\W}{{\rm W}}
\newcommand{\w}{{\rm w}}

\newcommand{\ch}{{\rm ch}}            
\newcommand{\td}{{\rm td}}            

\newcommand{\ad}{{\rm ad}}

\newcommand{\point}{{\rm pt}}
\newcommand{\ev}{{\rm ev}}

\newcommand{\MAX}{{\rm max}}

\newcommand{\dvol}{{\rm dvol}}
\newcommand{\loc}{{\rm loc}}

\newcommand{\eps}{{\varepsilon}}

\renewcommand{\i}{{\iota}}

\newcommand{\om}{{\omega}}

\newcommand{\Om}{{\Omega}}
\newcommand{\G}{{\rm G}}
\newcommand{\HG}{{\rm H}}
\newcommand{\EG}{{\rm EG}}
\newcommand{\BG}{{\rm BG}}

\newcommand{\SO}{{\rm SO}}
\newcommand{\U}{{\rm U}}

\renewcommand{\u}{{\mathfrak u}}

\newcommand{\g}{{\mathfrak g}}         
\newcommand{\h}{{\mathfrak h}}         
\renewcommand{\tt}{{\mathfrak t}}      
\newcommand{\Cinf}{C^{\infty}}

\newcommand{\SW}{{\rm SW}}
\newcommand{\Gr}{{\rm Gr}}

\newcommand{\Jreg}{{\mathcal{J}}_{\rm reg}}

\newcommand{\Hreg}{{\mathcal{H}}_{\rm reg}}

\newcommand{\PD}{{\rm PD}}

\newcommand{\QH}{{\rm QH}}
\newcommand{\GW}{{\rm GW}}
\newcommand{\inner}[2]{\langle #1, #2\rangle}

\def\NABLA#1{{\mathop{\nabla\kern-.5ex\lower1ex\hbox{$#1$}}}}
\def\Nabla#1{\nabla\kern-.5ex{}_{#1}}
\def\Tabla#1{\Tilde\nabla\kern-.5ex{}_{#1}}
\renewcommand{\Tilde}{\widetilde}

\newcommand{\p}{{\partial}}

\newcommand{\IMP}{\Longrightarrow}

\newcommand{\INTO}{\hookrightarrow}
\newcommand{\TO}{\longrightarrow}
\begin{document}

%


\title{The symplectic vortex equations \\
       and invariants of Hamiltonian group actions}

\author{Kai~Cieliebak,
        Stanford University\thanks{Partially supported 
        by National Science Foundation grant DMS-0072267} \\
        A.~Rita~Gaio,
        Fac. Ci\^encias - Univ. Porto\thanks{Partially supported 
        by the FCT grant PRAXIS XXI/BPD/22084/99} \\
        Ignasi~Mundet i Riera,
        Universidad Aut\'onoma de Madrid \\
        Dietmar~A.~Salamon,
        ETH-Z\"urich}

\date{3 October 2002}


\maketitle

\begin{abstract} 
In this paper we define invariants of Hamiltonian 
group actions for central regular values of the moment map.
The key hypotheses are that the moment map is proper
and that the ambient manifold is symplectically aspherical.
The invariants are based on the symplectic vortex equations.
Applications include an existence theorem for relative 
periodic orbits, a computation for circle actions on 
a complex vector space, and a theorem about the relation
between the invariants introduced here and the Seiberg--Witten
invariants of a product of a Riemann surface with a two-sphere.
\end{abstract}

\tableofcontents



\section{Introduction} 

In this paper we study the vortex equations with 
values in a symplectic manifold $(M,\om)$.
We assume that $(M,\om)$ is equipped with a Hamiltonian
action by a compact Lie group $\G$ that is generated by an 
equivariant moment map
$$
      \mu:M\to\g.
$$
The {\bf symplectic vortex equations} have the form 
\begin{equation}\label{eq:vortex}
      \bar\p_{J,A}(u)=0,\qquad *F_A+\mu(u)=\tau.
\end{equation}
Here $P\to\Sigma$ is a principal $\G$-bundle
over a compact Riemann surface, $u:P\to M$ is an equivariant 
smooth function, and $A$ is a connection on $P$. 
To define the terms in~(\ref{eq:vortex}) we must fix 
a $\G$-invariant almost complex structure on $M$,
a Riemannian metric on $\Sigma$, and an element 
$\tau\in Z(\g)$ in the center of the Lie algebra. 
The expression $\bar\p_{J,A}$ denotes the nonlinear Cauchy--Riemann
operator, twisted by the connection $A$, and $*$ denotes the 
Hodge $*$-operator on~$\Sigma$.  Equations~(\ref{eq:vortex})
were introduced in~\cite{CGS,GAIO,MUNDET}.
In the physics literature these equations are known as 
{\it gauged sigma models} in the case where the target space 
$M$ is a complex vector space. Special cases of the symplectic
vortex equations include pseudoholomorphic curves in symplectic 
manifolds ($\G=\{\one\}$),   
the usual vortex equations over $\Sigma$ 
($M=\C$ with the standard $S^1$-action~\cite{BRADLOW1,GP}),
Bradlow pairs ($M=\C^2$ with the standard $\U(2)$-action~\cite{BRADLOW2,T}),
anti-self-dual instantons over a product $\Sigma\times S$
($M$ is the infinite dimensional space of $\SO(3)$-connections
over $S$ and $\G$ is the gauge group),
and the Seiberg--Witten equations over $\Sigma\times S$ 
($M$ is the space of pairs, each consisting of a 
connection on a line bundle $L\to S$ and a holomorphic section, 
and $\G$ is the gauge group of~$L$). 
In the present paper the symplectic manifold $M$ is always 
finite dimensional. 

In the K\"ahler case the symplectic vortex equations admit 
an algebro geometric interpretation. For example, if $M$ is 
a complex vector space then the map $u$ defines a holomorphic 
section of a vector bundle over $\Sigma$ and the solutions 
of~(\ref{eq:vortex}) correspond to stable pairs. 

We impose the following conditions on the triple $(M,\om,\mu)$. 
\begin{description}
\item[(H1)] 
The moment map $\mu$ is proper. 
\item[(H2)] 
There exists a (strongly) convex structure $(f,J)$ on $M$.  
This means that $J$ is a $\G$-invariant and 
$\om$-compatible almost complex structure on $M$,
$f:M\to[0,\infty)$ is a proper $\G$-invariant function,
and there exists a continuous function 
$Z(\g)\to\R:\tau\mapsto c(\tau)$ such that
$$
      f(x)\ge c(\tau)\quad\IMP\quad
      \begin{array}{l}
      \inner{\Nabla{v}\nabla f(x)}{v}+\inner{\Nabla{Jv}\nabla f(x)}{Jv}\ge 0,\\
      df(x)J(x)L_x(\mu(x)-\tau) \ge 0
      \end{array}
$$
for every $x\in M$, $v\in T_xM$, and $\tau\in Z(\g)$.
Here $\nabla$ denotes the Levi-Civita connection of the 
metric $\inner{\cdot}{\cdot}=\om(\cdot,J\cdot)$
and $L_x:\g\to T_xM$ denotes the infinitesimal action.
\item[(H3)] 
The manifold $(M,\om)$ is symplectically aspherical, i.e. 
$$
     \int_{S^2}v^*\om=0
$$ 
for every smooth function $v:S^2\to M$.
\end{description}
Hypothesis~$(H2)$ is a natural generalization of the existence
of a plurisubharmonic function on noncompact symplectic 
manifolds~\cite{EG}. Both hypotheses~$(H1)$ and~$(H2)$ are natural 
in the context of this paper and are needed to obtain any kind of 
compactness theorem for the solutions of~(\ref{eq:vortex}).
Hypothesis~$(H3)$ constitutes a more severe restriction
and should in the future be removed or weakened. 
It implies that $M$ is noncompact whenever there exists a $\G$-orbit
of positive dimension. 
However, there are many interesting examples where all three
hypotheses are satisfied, e.g. linear actions on $\C^n$ with 
proper moment maps~\cite{CGS}.  The three hypotheses together 
guarantee that the moduli space of gauge equivalence classes 
of solutions of~(\ref{eq:vortex}) is compact.  As a result one can 
use these moduli spaces to define invariants which are analogous 
to the Gromov--Witten invariants in the nonequivariant case.
Let $B\in H_2^\G(M;\Z)$ denote the equivariant homology
class represented by the map $u$.  
Then the invariants take the form of a homomorphism
$$
     \Phi^{M,\mu-\tau}_{B,\Sigma}:H^*_\G(M;\Q)\to\Q,
$$
whenever $\tau$ is a central regular value of $\mu$. 
This homomorphism takes integer values on integral cohomology
classes whenever $\G$ acts freely on $\mu^{-1}(\tau)$. 
It depends only on the component of $\tau$ in the open set 
of central regular values.  We emphasize that the 
complex structure on $\Sigma$ is fixed in the definition of 
our invariant.  There should be natural extensions which involve 
varying complex structures on the domain $\Sigma$ and dispense 
with hypothesis~$(H3)$.  However, the definition of the invariants
in these cases will probably require a considerable amount of 
nontrivial analysis. For first steps in this direction
see~\cite{MUNDET,MUNDET2}. 

As a first application we establish the existence of relative 
periodic orbits for time dependent $\G$-invariant 
Hamiltonian systems.  This can be viewed as an equivariant 
version of a theorem of Gromov~\cite{GRO}. 

\medskip
\noindent{\bf Theorem~A.}
{\it Assume~$(H1-3)$ and let $\tau\in Z(\g)$ be a central
regular value of $\mu$ such that $\mu^{-1}(\tau)\ne\emptyset$.  
Then every time-dependent $1$-periodic $\G$-invariant Hamiltonian system
admits a contractible relative periodic orbit in $\mu^{-1}(\tau)$.}

\medskip
\noindent
If $\G$ is abelian then the hypothesis that $\tau$ is a regular
value can be dropped and we obtain a contractible relative periodic
orbit on every nonempty level set of $\mu$. It is natural to 
conjecture that this should continue to hold under 
hypothesis~$(H1)$ only.  Our proof of Theorem~A follows closely
Gromov's argument in~\cite{GRO} for the nonequivariant case. 
The pseudoholomorphic curves in Gromov's proof are 
replaced by the solutions of the perturbed symplectic vortex 
equations.  

In some cases the invariants can be computed explicitly.
We carry out such a computation for linear circle actions 
on $\C^n$.  Suppose $S^1$ acts on $\C^n$ with positive weights
$\ell_1,\dots,\ell_n$ and denote the correponding moment 
map by $\mu_\ell$.   Then there is only one nontrivial chamber 
for the regular values of $\mu_\ell$ and we denote by 
$\Phi^{\C^n,\mu_\ell}$ the invariant in this chamber.

\medskip
\noindent{\bf Theorem~B.}
{\it Let $\Sigma$ be a compact Riemann surface of genus $g$, 
$d\in\Z\cong H_2^{S^1}(\C^n;\Z)$ an integer, and $c\in
H^2_{S^1}(\C^n;\Z)\cong\Z$ the positive generator.  
Suppose that 
$$
      m := \sum_{\nu=1}^n(d\ell_\nu+1-g) + g-1 \ge 0.
$$
Then 
$$
      \Phi^{\C^n,\mu_\ell}_{d,g}(c^m)
      = \left(\sum_{\nu=1}^n\ell_\nu\right)^g
        \prod_{\nu=1}^n\ell_\nu^{-(d\ell_\nu+1-g)}.
$$
}

\medskip
\noindent
In the case $\ell_\nu=1$ and $d>2g-2$ this was proved by 
Bertram--Daskalopoulos--Wentworth~\cite{BDW}. 
The proof of Theorem~B involves the Atiyah--Singer index theorem 
for families of Cauchy--Riemann operators. 

Our invariants are related to the Seiberg--Witten invariants
of certain four-manifolds. The key observation is that 
the symmetric product of a Riemann surface $S$ 
can be interpreted as a symplectic quotient
of the infinite dimensional space whose elements are pairs,
each consisting of a connection and a holomorphic
section of a line bundle $L\to S$ of degree~$d$.  
In this situation the symplectic vortex equations, 
with $M$ replaced by the infinite dimensional space 
of which $\Sym^d(S)$ is a quotient, are the 
Seiberg--Witten equations on $\Sigma\times S$. 
When $d>2g_S-2$ one can write the symmetric product
as a quotient of a finite dimensional symplectic manifold 
$M_{d,S}$ (called the {\it vortex manifold of the pair} $(d,S)$)
by a Hamiltonian $S^1$ action with a moment map $\mu_{d,S}$
which satisfies~$(H1-3)$. The following theorem relates the invariants
of $(M_{d,S},\mu_{d,S})$ to the Seiberg--Witten invariants. It is a
special case of a result for general ruled surfaces in~\cite{OT}. 

\medskip
\noindent{\bf Theorem~C (\cite{OT}).}
{\it Let $S$ be the Riemann sphere and $\Sigma$ 
be a compact Riemann surface of genus $g$. 
Let $d$ and $k$ be nonnegative integers such that
$$
     m := d(1-g) + (d+1)k \ge 0.
$$
Then 
$$
     \Phi^{M_{d,S},\mu_{d,S}}_{k,\Sigma}(c^m) 
     = \SW_{\Sigma\times S}(\gamma_{k,d}),
$$
where $\gamma_{k,d}$ denotes the spin$^c$-structure 
determined by $k$ and $d$.  Moreover, if 
$$
     k > 2g-2,
$$
then 
$$
     \Phi^{M_{d,S},\mu_{d,S}}_{k,\Sigma}(c^m) 
     = \Phi^{M_{k,\Sigma},\mu_{k,\Sigma}}_{d,S}(c^m).
$$
}

\medskip
\noindent
Combining Theorems~B and~C one can recover the computation 
of the Seiberg--Witten invariants of product ruled surfaces 
by Li-Liu~\cite{LL} and Ohta--Ono~\cite{OO}. 

It is also interesting to examine the relation between 
our invariants and the Gromov--Witten invariants of the 
symplectic quotient 
$$
     \bar M := M\dslash\G(\tau) := \mu^{-1}(\tau)/\G
$$
whenever $\G$ acts freely on $\mu^{-1}(\tau)$.  
Such a relation was established in~\cite{GaSa} under the hypothesis
that the quotient is {\it monotone}.  Under this condition
(and hypotheses~$(H1-3)$) it is shown in~\cite{GaSa}
that there exists a surjective ring homomorphism
$$
    \phi:H^*_\G(M) \to \QH^*(\bar M)
$$
(with values in the quantum cohomology of the quotient)
such that 
$$
    \Phi^{M,\mu-\tau}_{B,\Sigma}(\alpha)
    = \GW^{\bar M}_{\bar B,\Sigma}(\phi(\alpha))
$$
for every $\alpha\in H^*_\G(M)$ and every $\bar B\in H_2(\bar M;\Z)$,
where $B$ denotes the image of $\bar B$ under the homomorphism
$H_2(\bar M;\Z)\to H_2^\G(M;\Z)$.  The proof is based on an adiabatic 
limit analysis which relates the solutions of the symplectic 
vortex equations in $M$ to pseudoholomorphic curves in the 
symplectic quotient $\bar M$.  This analysis is analogous to the 
proof of the Atiyah--Floer conjecture in~\cite{DS}. 

The present paper is organized as follows.  
In Section~\ref{sec:ham} we discuss the basic properties
of solutions to the symplectic vortex equations such as 
the energy identity, unique continuation, 
and apriori estimates under the convexity
hypothesis~$(H2)$.  Section~\ref{sec:CR} establishes the basic 
compactness and regularity theorems and Section~\ref{sec:trans}
discusses the Fredholm theory.  In Section~\ref{sec:integer}
we establish the integer invariants under the hypothesis
that $\G$ acts freely on $\mu^{-1}(\tau)$.  Section~\ref{sec:moduli}
is of preparatory nature. In it we recall some background from~\cite{CMS} 
about the equivariant Euler class of $\G$-moduli pro\-blems. 
Section~\ref{sec:rational} establishes the rational invariants 
in the presence of finite isotropy and discusses some relations 
between the invariants.  Theorems~A, B, and~C are proved 
in Sections~\ref{sec:periodic}, \ref{sec:proj}, and~\ref{sec:sw}.
Appendix~\ref{app:KW} establishes existence and uniqueness for a coupled
Kazdan-Warner equation that appears in the proof of Theorem~C.  
Appendix~\ref{app:slice} gives a proof of the local slice theorem for 
gauge group actions in a form needed for the compactness 
and regularity results of Section~\ref{sec:CR}.


\section{The geometry of symplectic vortices}\label{sec:ham}


\subsection{The vortex equations in a symplectic manifold}

Let $(M,\om)$ be a (not necessarily compact)
symplectic manifold and $\G$ be a compact 
Lie group which acts on $M$ by symplectomorphisms.
Let $\g=\Lie(\G)$ denote the Lie algebra and
$$
     \g\to\Vect(M,\om):\xi\mapsto X_\xi
$$
denote the infinitesimal action. We assume that
the action is Hamiltonian.  This means that the
action is generated by an equivariant map $\mu:M\to\g$
that satisfies
$$
     \i(X_\xi)\om = d\inner{\mu}{\xi}
$$
for every $\xi\in\g$. Here $\inner{\cdot}{\cdot}$
denotes an invariant inner product on $\g$. 
The function $\mu$ is called a moment map for the action. 

Let $P\to\Sigma$ be a principal $\G$-bundle 
over a compact connected oriented Riemann surface
$(\Sigma,J_\Sigma,\dvol_\Sigma)$.
We emphasise that the volume form and the 
complex structure on $\Sigma$ are fixed. 
Denote by $\Cinf_\G(P,M)$ the space of equivariant
functions $u:P\to M$ and by $\Aa(P)$ the space 
of connections on~$P$. We think of $A\in\Aa(P)$ 
as an equivariant Lie algebra valued $1$-form on $P$ which 
identifies the vertical tangent space with $\g$. 
Its curvature is a $2$-form $F_A$ on $\Sigma$
with values in the associated Lie algebra 
bundle $\g_P:=P\times_{\ad}\g$. 
In this paper we study the following system of nonlinear
first order partial differential equations,
for pairs $(u,A)\in\Cinf_\G(P,M)\times\Aa(P)$,
$$
     \bar\p_{J,A}(u) = 0,\qquad
     *F_A + \mu(u) = \tau.
$$
Here $\tau\in Z(\g)$ is an element in the centre
of the Lie algebra,
$\Jj_G(M,\om)$ denotes the space of $\G$-invariant
and $\om$-tame almost complex structures on~$M$,
and $J:\Sigma\to\Jj_\G(M,\om)$ is a smooth family 
of such almost complex structures.  The space of 
such families of almost complex structures
will be denoted by 
$$
     \Jj 
     := \Jj(\Sigma,M,\om,\mu)
     := \Cinf(\Sigma,\Jj_\G(M,\om)).
$$ 
The covariant derivative of $u$ with respect 
to the connection $A$ is the $1$-form $d_Au\in\Om^1(P,u^*TM)$
given by 
$$
     d_Au := du + X_A(u).
$$
This $1$-form is equivariant and horizontal
and hence descends to a $1$-form on $\Sigma$ 
with values in $u^*TM/\G$. The family of almost 
complex structures $J$ determines a $\G$-invariant
complex structure
$$
     J_u(p) := J(\pi(p),u(p))
$$
on the bundle $u^*TM\to P$ and hence a complex structure
on the bundle $u^*TM/\G\to\Sigma$, which will also be 
denoted by $J_u$. The term $\bar\p_{J,A}(u)$
denotes the $(0,1)$-part of this $1$-form
and so is a $(0,1)$-form on $\Sigma$ with values in
$u^*TM/\G$. Its lift to a $1$-form on $P$ with values
in $u^*TM$ will also be denoted by $\bar\p_{J,A}(u)$
and is given by 
$$
     \bar\p_{J,A}(u) 
     := \frac12\left(d_Au + J\circ d_Au\circ J_\Sigma\right).
$$
The right hand side is well defined since $d_Au$ 
is horizontal.  Namely, given a tangent vector $v\in T_pP$, 
lift the vector $J_\Sigma d\pi(p)v\in T_{\pi(p)}\Sigma$
to $T_pP$, and apply the linear map
$J(\pi(p),u(p))d_Au(p)$ to the lift.
The resulting vector in $T_{u(p)}M$ is independent
of the choice of the lift, because the $1$-form 
$d_Au$ vanishes on vertical tangent vectors. 
Equations~(\ref{eq:vortex}) were 
introduced in~\cite{CGS,GAIO,MUNDET}.
They are a generalized form of the vortex equations.
In the case of linear actions on $\C^n$ they are
known in the physics literature 
as {\it gauged sigma models}. 

\begin{remark}\label{rmk:moment}\rm
The space $\Cinf_\G(P,M)\times\Aa(P)$ is an infinite
dimensional Fr\'echet manifold and admits a natural 
symplectic structure.  The gauge group $\Gg(P)$ 
acts on this space by 
$$
      g^*(u,A)=(g^{-1}u,g^{-1}dg+g^{-1}Ag).
$$
This action is Hamiltonian and 
the function
\begin{equation}\label{eq:moment}
      \Cinf_\G(P,M)\times\Aa(P)
      \to\Cinf(\Sigma,\g_P):
      (u,A)\mapsto *F_A+\mu(u)
\end{equation}
is a moment map for this action (see~\cite{CGS}).
The space of solutions of~(\ref{eq:vortex})
is invariant under the action of $\Gg(P)$.
The quotient can be interpreted as a symplectic 
quotient whenever the space of pairs 
$(u,A)$ that satisfy $\bar\p_{J,A}(u)=0$
is a symplectic submanifold of 
$\Cinf_\G(P,M)\times\Aa(P)$.
\end{remark}


\subsection{Hamiltonian perturbations}\label{sec:H}

Le $\Cinf_\G(M)$ be the space 
of smooth $\G$-invariant functions on $M$.
A {\bf Hamiltonian perturbation} is a $1$-form 
$H\in\Om^1(\Sigma,\Cinf_\G(M))$.
One can think of $H$ as a $\G$-equivariant
section of the vector bundle 
$T^*\Sigma\times M\to\Sigma\times M$.
The space of Hamiltonian perturbations will be denoted
by
$$
     \Hh
     := \Hh(\Sigma,M,\om,\mu)
     := \Om^1(\Sigma,\Cinf_\G(M)).
$$
For $H\in\Hh$ and $\zeta\in T_z\Sigma$ we write 
$H_\zeta := H_z(\zeta)\in\Cinf_\G(M)$ 
and denote by $X_{H_\zeta}\in\Vect(M,\om)$
the $\G$-invariant Hamiltonian vector field 
of $H_\zeta$, i.e.
$$
     \i(X_{H_\zeta})\om = dH_\zeta.
$$
A Hamiltonian perturbation $H\in\Hh$ 
and a section $u\in\Cinf_\G(P,M)$ determine a $1$-form
$X_H(u)\in\Om^1(P,u^*TM)$ given by 
$$
     (X_H(u))_p(v) := X_{H_{d\pi(p)v}}(u(p)).
$$
This $1$-form is equivariant and horizontal and so is
$$
     d_{H,A}(u) := d_Au + X_H(u).
$$
Hence $X_H(u)$ and $d_{H,A}(u)$ descend to
$1$-forms on $\Sigma$ with values in $u^*TM/\G$.  
We denote
$$
     \bar\p_{J,H,A}(u) 
     :=  (d_{H,A}(u))^{0,1}
      =  \bar\p_{J,A}(u) + (X_H(u))^{0,1}
     \in \Om^{0,1}_{J_u}(\Sigma,u^*TM/\G)
$$
and replace~(\ref{eq:vortex})
by the perturbed equations
\begin{equation}\label{eq:vortex-ham}
     \bar\p_{J,H,A}(u)=0,\qquad
     *F_A + \mu(u) = \tau.
\end{equation}


\subsection{Energy}\label{sec:energy}

Fix a central element $\tau\in Z(\g)$, 
an almost complex structure $J\in\Jj(\Sigma,M,\om,\mu)$,
and a perturbation $H\in\Om^1(\Sigma,\Cinf_\G(M))$.
The {\bf energy} of a pair $(u,A)\in\Cinf_\G(P,M)\times\Aa(P)$
is defined by
$$
     E(u,A) 
     := \frac12\int_\Sigma
       \left(
         \left|d_{H,A}(u)\right|^2
         + \left|F_A\right|^2
         + \left|\mu(u)-\tau\right|^2
       \right)
       \dvol_\Sigma.
$$
This functional is invariant under the action
of the gauge group $\Gg(P)$.  Denote
by $[\om+\tau-\mu]\in H^2(M_\G;\R)$
the equivariant cohomology class determined by
the symplectic form $\om$ and the moment map
$\mu-\tau$ (see~\cite{CGS}), and denote by $[u]\in H_2(M_\G;\Z)$ 
the homology class determined by $u$.
More precisely, there is an equivariant classifying map 
$\theta:P\to\EG$ and hence the map $(u,\theta):P\to M\times\EG$
descends to a map 
$$
     u_\G:\Sigma\to M_\G:=M\times_\G\EG.
$$
The class $[u]\in H_2(M_\G;\Z)$ is defined
as the pushforward of the fundamental class 
$[\Sigma]$ under the map induced by $u_\G$. 
For every pair $(u,A)\in\Cinf_\rho(P,M)\times\Aa(P)$
the cohomology pairing between the classes
$[\om+\tau-\mu]$ and $[u]$ is given by 
$$
     \inner{[\om+\tau-\mu]}{[u]}
     = \int_\Sigma \left(u^*\om-d\inner{\mu(u)-\tau}{A}\right).
$$
This topological invariant appears in the following
energy identity. Another ingredient in this formula
is the curvature of the Hamiltonian connection $H$.
Since $\Sigma$ carries a volume form this curvature 
can be expressed as a function $\Om_H:\Sigma\times M\to\R$
which is invariant under the $\G$-action on $M$.
It is defined by the formula
$$
     \Om_H\,\dvol_\Sigma
     := d^\Sigma H + \frac12\{H\wedge H\}
     \;\;\in\;\;\Om^2(\Sigma,\Cinf_\G(M)),
$$
where $\{\cdot,\cdot\}$ denotes the Poisson bracket
for functions on $M$.  
The {\bf Hofer norm} of the curvature $\Om_H$ is defined by
$$
     \|\Om_H\|
     := \int_\Sigma\left(\sup_{x\in M}\Om_H(z,x)
        -\inf_{x\in M}\Om_H(z,x)\right)\dvol_\Sigma.
$$
This quantity is independent of the volume form of $\Sigma$. 
The next proposition states the basic energy identity.  
The first term on the right is the $L^2$-norm of the terms 
in equation~(\ref{eq:vortex-ham}), the second term is a topological 
invariant, and the last term is bounded by $\|\Om_H\|$.

\begin{proposition}\label{prop:energy}
For every $A\in\Aa(P)$ and every $u\in\Cinf_\G(P,M)$,
\begin{eqnarray}\label{eq:energy}
     E(u,A) 
&= &
    \int_\Sigma
       \left(\left|\bar\p_{J,H,A}(u)\right|^2
       + \frac12\left|*F_A+\mu(u)-\tau\right|^2\right)\dvol_\Sigma
    \nonumber\\
&& 
    +\,\inner{[\om+\tau-\mu]}{[u]}
    + \int_\Sigma\Om_H(u)\,\dvol_\Sigma. 
\end{eqnarray}
In particular,
$
     E(u,A) \le \inner{[\om+\tau-\mu]}{[u]} + \|\Om_H\|
$
for every solution of~(\ref{eq:vortex-ham}).
\end{proposition}

\begin{proof}  
Choose a holomorphic coordinate chart
$\phi:U\to\Sigma$, where $U\subset\C$ is an open set,
and let $\tilde\phi:U\to P$ be a lift of $\phi$,
that is $\pi\circ\tilde\phi=\phi$.   
Then $u$, $A$, and $H$ are in local coordinates given by 
$$
     u^\loc:=u\circ\tilde\phi,\qquad
     {\tilde\phi\,}^*A=\Phi\,ds+\Psi\,dt,\qquad
     \phi^*H = F\,ds + G\,dt
$$
where $\Phi,\Psi:U\to\g$ and $F,G:U\times M\to\R$. 
The pullback volume form on $U$ is
$
     \dvol_\Sigma^\loc = \lambda^2\,ds\wedge dt
$
for some function $\lambda:U\to(0,\infty)$
and the metric is $\lambda^2(ds^2+dt^2)$.
Hence
$$
\begin{array}{rcl}
     \phi^*\Om_H 
&= &
     \lambda^{-2}\left(\p_sG-\p_tF+\{F,G\}\right),  \\
     {\tilde\phi\,}^*F_A
&= &
     \left(\p_s\Psi-\p_t\Phi+[\Phi,\Psi]\right) 
           \,ds\wedge dt,  \\
     {\tilde\phi\,}^*d_Au
&= &
     \left(\p_su^\loc+X_\Phi(u^\loc)\right)\,ds 
     + \left(\p_tu^\loc+X_\Psi(u^\loc)\right)\,dt, \\
     {\tilde\phi\,}^*\bar\p_{J,A}(u)
&= &
     \frac12(\xi\,ds -J^\loc(s,t,u^\loc)\xi\,dt),
\end{array}
$$
where $\{F,G\}:=\om(X_F,X_G)$ denotes the 
Poisson bracket on $M$, 
$$
     \xi 
     := \p_su^\loc+X_\Phi(u^\loc)
       + J^\loc(s,t,u^\loc)
         \left(\p_tu^\loc+X_\Psi(u^\loc)\right),
$$
and $J^\loc(s,t,x):=J(\phi(s,t),x)$ for $(s,t,x)\in U\times M$.
In the following we shall drop the superscript ``$\loc$''.
Then~(\ref{eq:vortex-ham}) have the form
\begin{equation}\label{eq:new-loc}
\begin{array}{rcl}
    \p_su+X_\Phi(u) + X_F(u)   
    + J\left(\p_tu+X_\Psi(u)+X_G(u)\right) &= &0,  \\
    \p_s\Psi-\p_t\Phi+[\Phi,\Psi] + \lambda^2(\mu(u)-\tau) &= & 0,
\end{array}
\end{equation}
The pullback of the energy integrand under $\phi:U\to\Sigma$
is given by $e\,ds\wedge dt$ where $e:U\to\R$ is the function
\begin{eqnarray*}
    e 
&:= &
    \frac12\left|\p_su+X_\Phi(u)+X_F(u)\right|^2
        + \frac12\left|\p_tu+X_\Psi(u)+X_G(u)\right|^2  \\
&&
    + \,\frac{1}{2\lambda^2}
      \left|\p_s\Psi-\p_t\Phi+[\Phi,\Psi]\right|^2
    + \frac{\lambda^2}{2}\left|\mu(u)-\tau\right|^2  \\
&= &
    \frac12\left|\p_su+X_\Phi(u)+X_F(u)
     + J(\p_tu+X_\Psi(u)+X_G(u))\right|^2 \\
&&
    +\,\frac{\lambda^2}{2}
        \left|\lambda^{-2}\left(\p_s\Psi-\p_t\Phi+[\Phi,\Psi]\right) 
         + \mu(u)-\tau\right|^2  + R.
\end{eqnarray*}
The remainder term $R$ has the form
\begin{eqnarray*}
    R
&:=&
    \om\left(\p_su+X_\Phi(u)+X_F(u),
    \p_tu+X_\Psi(u)+X_G(u)\right) \\
&& 
    - \inner{\p_s\Psi-\p_t\Phi+[\Phi,\Psi]}{\mu(u)-\tau} \\
&= &
    \om(\p_su,\p_tu) 
    - \p_s\left(G(u)+\inner{\mu(u)-\tau}{\Psi}\right)
    + \p_t\left(F(u)+\inner{\mu(u)-\tau}{\Phi}\right)  \\
&&
    +\,\left(\p_sG-\p_tF+\{F,G\}\right)\circ u.
\end{eqnarray*}
This proves~(\ref{eq:energy}).
\end{proof}


\subsection{Unique continuation}\label{sec:unique}

A solution $(u,A)$ of~(\ref{eq:vortex-ham}) is called
{\bf horizontal} if $d_{H,A}(u)\equiv0$ and $\mu(u)\equiv\tau$.

\begin{lemma}\label{le:horizontal}
Let $(u,A)$ be a solution of~(\ref{eq:vortex-ham})
with $H=0$.  Then $(u,A)$ is horizontal if and only 
if the homology class $[u]\in H_2(M_\G;\Z)$ is torsion.
\end{lemma}

\begin{proof}
The ``if'' part follows from the energy identity. 
To prove the converse note that $A$ is flat for every 
horizontal solution $(u,A)$ of~(\ref{eq:vortex-ham}) 
and hence, in the case $H=0$, every equivariant cohomology 
class vanishes on $[u]$ (see~\cite{CGS}).
\end{proof}

\begin{lemma}\label{le:unique}
Suppose $u$ and $A$ satisfy~(\ref{eq:vortex-ham}) with $H=0$.
If $d_Au$ and $\mu(u)-\tau$ vanish to infinite order at 
some point $p_0\in P$ then $(u,A)$ is horizontal. 
\end{lemma}

\begin{proof}
Replacing $\mu$ by $\mu-\tau$ we may assume that $\tau=0$.
Consider the equations in their local coordinate
form~(\ref{eq:new-loc}). In the case $H=0$ they read
$$
     v_s+Jv_t=0,\qquad \kappa+\lambda^2\mu(u) = 0
$$
for $u:U\to M$ and $\Phi,\Psi:U\to\g$, where
$v_s,v_t:U\to u^*TM$ and $\kappa:U\to\g$ are 
defined by
$$
     v_s:=\p_su+L_u\Phi,\qquad
     v_t:=\p_tu+L_u\Psi,\qquad
     \kappa:=\p_s\Psi-\p_t\Phi+[\Phi,\Psi],
$$
where $L_x:\g\to T_xM$ denotes the infinitesimal
action given by $L_x\eta:=X_\eta(x)$. Let us denote 
$$
\begin{array}{rclrcl}
     \Nabla{A,s}\xi &:=& \Nabla{s}\xi + \Nabla{\xi}X_\Phi(u),&
     \Nabla{A,t}\xi &:=& \Nabla{t}\xi + \Nabla{\xi}X_\Psi(u),\\
     \Nabla{A,s}\eta &:=& \p_s\eta + [\Phi,\eta],&
     \Nabla{A,t}\eta &:=& \p_t\eta + [\Psi,\eta],
\end{array}
$$
for $\xi:U\to u^*TM$ and $\eta:U\to\g$.  Then
\begin{equation}\label{eq:vst}
     \Nabla{A,s}v_t-\Nabla{A,t}v_s = L_u\kappa = -\lambda^2L_u\mu(u)
\end{equation}
\begin{equation}\label{eq:nabla-eta}
     \Nabla{A,s}L_u\eta - L_u\Nabla{A,s}\eta =
     \Nabla{v_s}X_\eta(u),\qquad
     \Nabla{A,t}L_u\eta - L_u\Nabla{A,t}\eta = 
     \Nabla{v_t}X_\eta(u).
\end{equation}
Since $d\mu(u)=-L_u^*J$ we have
\begin{equation}\label{eq:mu}
     \Nabla{A,s}\mu(u)=d\mu(u)v_s=-L_u^*v_t,\qquad
     \Nabla{A,t}\mu(u)=d\mu(u)v_t=L_u^*v_s.
\end{equation}
It follows from~(\ref{eq:vst}), (\ref{eq:nabla-eta}), and~(\ref{eq:mu}) 
that
\begin{equation}\label{eq:vst-1}
     \Nabla{A,t}\left(\Nabla{A,t}v_s - \Nabla{A,s}v_t\right) 
     = (\p_t\lambda^2)L_u\mu(u) + \lambda^2d\mu(u)v_t
       + \lambda^2\Nabla{v_t}X_{\mu(u)}(u).
\end{equation}
Since
$
     0 = (\Ll_{X_\eta}J)\xi
       = (\Nabla{X_\eta}J)\xi + J\Nabla{\xi}X_\eta - \Nabla{J\xi}X_\eta,
$
we obtain 
$$
     \Nabla{A,s}J = \Nabla{v_s}J+\p_sJ,\qquad
     \Nabla{A,t}J = \Nabla{v_t}J+\p_tJ,
$$
and hence, using the identity $(\Nabla{v_s}J)v_t=(\Nabla{v_t}J)v_s$,
\begin{eqnarray}\label{eq:vst-2}
     \Nabla{A,s}v_s + \Nabla{A,t}v_t 
&= &
     \Nabla{A,t}(Jv_s) - \Nabla{A,s}(Jv_t) 
     \nonumber \\
&= &
     (\p_tJ)v_s - (\p_sJ)v_t + \lambda^2JL_u\mu(u).
\end{eqnarray}
This gives rise to an inequality of the form 
$$
    \left|\Nabla{A,s}\left(\Nabla{A,s}v_s + \Nabla{A,t}v_t
    \right)\right|
    \le c_1\left(|v_s| + |\p_sv_s| + |\p_tv_s| + |\mu(u)|\right).
$$
Moreover, by~(\ref{eq:vst-1}), 
$$
    \left|\Nabla{A,t}\left(\Nabla{A,t}v_s - \Nabla{A,s}v_t
    \right)\right|
    \le c_2\left(|v_s| + |\mu(u)|\right)
$$
and, by the curvature identity
$
     \Nabla{A,s}\Nabla{A,t}\xi-\Nabla{A,t}\Nabla{A,s}\xi
     = R(v_s,v_t)\xi + \Nabla{\xi}X_\kappa(u),
$
$$
    \left|\Nabla{A,s}\Nabla{A,t}v_t - \Nabla{A,t}\Nabla{A,s}v_t
    \right|
    \le c_3|v_s|.
$$
Putting these three inequalities together we obtain 
$$
    \left|\Nabla{A,s}\Nabla{A,s}v_s + \Nabla{A,t}\Nabla{A,t}v_s
    \right|
    \le c_4\left(|v_s| + |\p_sv_s| + |\p_tv_s| + |\mu(u)|\right).
$$
Moreover,
$$
     \Nabla{A,s}(L_u^*\xi)-L_u^*\Nabla{A,s}\xi = \rho(v_s,\xi),\qquad
     \Nabla{A,t}(L_u^*\xi)-L_u^*\Nabla{A,t}\xi = \rho(v_t,\xi),
$$
where $\rho\in\Om^2(M,\g)$ is defined by 
$
     \inner{\eta}{\rho(\xi_1,\xi_2)}:=\inner{\Nabla{\xi_1}X_\eta}{\xi_2}
$
(see~\cite{GaSa}). Hence, by~(\ref{eq:mu}) and~(\ref{eq:vst}), 
$$
     \Nabla{A,s}\Nabla{A,s}\mu(u) + \Nabla{A,t}\Nabla{A,t}\mu(u)  
     = \lambda^2L_u^*L_u\mu(u) - 2\rho(v_s,v_t). 
$$
Hence there exists a constant $c$ such that
$$
     |\Delta v_s|\le c\left(
     |v_s| + |\p_sv_s| + |\p_tv_s| + |\mu(u)|
     \right),\qquad
     |\Delta \mu(u)|\le c\left(|\mu(u)| + |v_s|\right).
$$
Hence it follows from Aronszajn's theorem that, if $v_s$ and $\mu(u)$
vanish to infinite order at a point in $U$ and $U$
is connected, then $v_s$ and $\mu(u)$ vanish identically
on $U$. 
\end{proof}

\begin{lemma}\label{le:horizontal1}
Suppose $u$ and $A$ satisfy~(\ref{eq:vortex-ham}) with $H=0$.
If there exists an open set $U\subset P$ such that
$d_Au(p)=0$ and $L_{u(p)}:\g\to T_{u(p)}M$ 
is injective for every $p\in U$ then $(u,A)$ is horizontal.
\end{lemma}

\begin{proof} 
By~(\ref{eq:vst}), 
$
     L_u(\mu(u)-\tau) = 0
$
in $U$.  Hence $\mu(u)=\tau$ in $U$ and hence, by Lemma~\ref{le:unique},
we have $d_Au\equiv0$ and $\mu(u)\equiv\tau$.
\end{proof}


\subsection{Convexity}\label{sec:convex}

\begin{definition}\label{def:convex}
A {\bf convex structure} on $(M,\om,\mu)$ is a pair 
$(f,J)$ where $J\in\Jj(M,\om)$ is a $\G$-invariant 
$\om$-compatible almost complex structure on $M$
and $f:M\to[0,\infty)$ is a smooth function 
satisfying the following conditions. 
\begin{description}
\item[(C1)]  $f$ is $\G$-invariant and proper.
\item[(C2)]  There exists a constant $c_0>0$ such that
$$
      f(x)\ge c_0\qquad\IMP\qquad
      \inner{\Nabla{\xi}\nabla f(x)}{\xi} 
      + \inner{\Nabla{J\xi}\nabla f(x)}{J\xi} \ge 0
$$
for every $x\in M$ and every $\xi\in T_xM$.
Here $\nabla$ denotes the Levi-Civita connection of the 
metric $\inner{\cdot}{\cdot}=\om(\cdot,J\cdot)$.
\item[(C3)]  There exists a constant $c_0>0$ such that
$$
      f(x) \ge c_0\qquad\IMP\qquad
      df(x)J(x)L_x\mu(x) \ge 0
$$
for every $x\in M$.
\end{description}
The second hypothesis says that the upward gradient flow
of $f$ expands the metric outside of a sufficiently
large compact set.  It is sometimes useful to assume 
condition~$(C3)$ for all moment maps $\mu-\tau$.
\begin{description}
\item[(C3')]  There exists a continuous function
$Z(\g)\to\R:\tau\mapsto c_0(\tau)$ such that
$$
      f(x) \ge c_0(\tau)\qquad\IMP\qquad
      df(x)J(x)L_x(\mu(x)-\tau)\ge 0
$$
for every $\tau\in Z(\g)$ and every $x\in M$.
\end{description}
A convex structure $(f,J)$ that satisfies~$(C3')$ is called 
{\bf strongly convex}. 
\end{definition}

\begin{lemma}\label{le:apriori}
Fix a homology class $B\in H_2(M_\G;\Z)$ and 
let $(f,J_0)$ be a convex structure for $(M,\om,\mu-\tau)$.  Denote
$$
     M_0 := \left\{x\in M\,|\,f(x)\le c_0\right\},
$$
where $c_0=c_0(\tau)$ is chosen such that $(C2)$ and~$(C3')$ 
are satisfied and
\begin{equation}\label{eq:apriori}
     f(x) > c_0 \qquad\IMP\qquad
     \left|\mu(x)-\tau\right|^2 > 
     \frac{\inner{[\om+\tau-\mu]}{B}}{\Vol(\Sigma)}
\end{equation}
Let $P$ be a principal $\G$-bundle over a compact Riemann surface 
$\Sigma$, suppose that $J\in\Jj$ agrees 
with $J_0$ outside $M_0$, and let $H\in\Hh$ be a Hamiltonian 
perturbation with support in $M_0$.  Then every solution 
$(u,A)$ of~(\ref{eq:vortex-ham}) representing 
the class $B$ satisfies $u(P)\subset M_0$.
\end{lemma}

\begin{proof}
Assume without loss of generality that $\tau=0$
and continue the notation of the proof of Lemma~\ref{le:unique}.
Since $\nabla f$ is a $\G$-invariant vector field
we have $[\nabla f,X_\eta]=0$ for every $\eta\in\g$ and hence 
$$
     \Nabla{A,s}\nabla f(u) = \Nabla{v_s}\nabla f(u),\qquad
     \Nabla{A,t}\nabla f(u) = \Nabla{v_t}\nabla f(u).
$$
Let
$
     \Delta:=\p_s^2+\p_t^2
$
denote the standard Laplacian. If $f(u)>c_0$ then
\begin{eqnarray*}
    \Delta f(u)
&= &
    \p_s\inner{\nabla f(u)}{v_s} + \p_t\inner{\nabla f(u)}{v_t}  \\
&= &
    \inner{\Nabla{A,s}\nabla f(u)}{v_s} 
    + \inner{\Nabla{A,t}\nabla f(u)}{v_t}  
    + \inner{\nabla f(u)}{\Nabla{A,s}v_s + \Nabla{A,t}v_t}  \\
&= &
    \inner{\Nabla{v_s}\nabla f(u)}{v_s} 
    + \inner{\Nabla{v_t}\nabla f(u)}{v_t}  
    + \lambda^2\inner{\nabla f(u)}{JL_u\mu(u)}  \\
&\ge &
     \inner{\Nabla{v_s}\nabla f(u)}{v_s} 
    + \inner{\Nabla{v_t}\nabla f(u)}{v_t}.
\end{eqnarray*}
Here the third equality follows from~(\ref{eq:vst-2})
and the last inequality from~$(C3)$. 
Now suppose, by contradiction, that 
$
    m:=\max_P f\circ u>c_0
$
and choose a local coordinate chart as above 
near a point where $f\circ u$ attains its maximum.
Since $f\circ u$ is subharmonic it follows from the 
mean value inequality that $f\circ u=m$ in a neighbourhood
of this maximum.  Hence the subset of $P$ where $f\circ u=m$
is open and closed, and hence $f\circ u\equiv m$.  
Hence, by~(\ref{eq:apriori}),
$$
      E(u,A)
      \ge \int_\Sigma\left|\mu(u)\right|^2\dvol_\Sigma
      \ge \Vol(\Sigma)\inf_P\left|\mu\circ u\right|^2
      > \inner{[\om-\mu]}{B}.
$$
Since $\Om_H$ vanishes on the image of $u$,
this contradicts the energy identity. 
\end{proof}

\begin{example}[\cite{CGS}]\label{ex:convex-lin}\rm
Consider the linear action of a compact Lie group
$\G$ on $\C^n$ by a homomorphism $\rho:\G\to\U(n)$
with proper moment map $\mu_\rho:\C^n\to\g$ given by 
$$
     \mu_\rho(x) := \pi_\rho\left(-\frac{i}{2}xx^*\right)
$$
where $\pi_\rho:=\dot\rho^*:\u(n)\to\g$ is the dual 
operator of the Lie algebra homomorphism 
$\dot\rho:\g\to\u(n)$ with respect to the 
inner product $\inner{A}{B}:=\trace(A^*B)$
on $\u(n)$.  A strongly convex structure for $\rho$
is the pair $(f,J)$ with $J=i$ and 
$$
     f(x) = \frac12|x|^2.
$$
To see this note that $\nabla f(x)=x$ and 
$$
     df(x)JL_x(\mu_\rho(x)-\tau) 
     = \inner{\mu_\rho(x)}{\mu_\rho(x)-\tau}
$$
Note that the pair $(f,i)$ with $f(x)=|\mu(x)|^2/2$ need not 
be a convex structure.  An example is the action of $\T^2$
on $\C^2$ by $(t_1,t_2)\cdot(x_1,x_2)=(t_1x_1,t_2x_2)$.
\end{example}

\begin{example}[Contact boundaries]\label{ex:liouville}\rm
Suppose $(M,\om)$ is a compact symplectic manifold with boundary
$\p M$, equipped with a Hamiltonian $\G$-action
generated by a moment map $\mu:M\to\R$.  Suppose that 
$X\in\Vect(M)$ is a $\G$-invariant vector field which points
out on the boundary and satisfies
$$
     \Ll_X\om = \om,\qquad \om(X_\mu,X)\ge 0
$$
near $\p M$.  Such a vector field gives rise to a convex structure
as follows.  Let $\phi_t$ denote the flow of $X$ and choose a
$\G$-invariant $\om$-compatible almost complex structure $J$ 
on $M$ such that
$$
     d\phi_t(x)J(x)=J(\phi_t(x))d\phi_t(x),\qquad
     \om(X(x),J(x)X(x))=1,
$$
and
$
     \om(X(x),J(x)v)=0
$
for $x\in\p M$, $v\in T_x\p M$, and $-\eps<t\le0$. 
Then the function $f:M\to\R$, defined by 
$$
     f(\phi_t(x)) := t
$$
for $x\in\p M$ and $-\eps<t\le 0$ defines a convex structure 
near $\p M$.  Its gradient is the vector field $X$ and its
covariant Hessian is half the identity.  Moreover, the manifold
can be extended by attaching a cylindrical end of the form 
$\p M\times[0,\infty)$ with the obvious extensions of the 
symplectic and almost complex structures to obtain a noncompact
manifold as above. 
\end{example}

\begin{example}[Convex fibrations]\label{ex:convfibr}\rm
Let $\G$ and $\HG$ be compact connected Lie groups
with Lie angebras $\g=\Lie(\G)$ and $\h=\Lie(\HG)$. 
Let $(M,\om)$ be a (not necessarily compact)
symplectic manifold equipped with Hamiltonian
action by both Lie groups $\G$ and $\HG$, generated
by moment maps $\mu_\G:N\to\g$ and $\mu_\HG:N\to\h$.
We assume throughout that the action of $\G$ 
commutes with the action of $\HG$. This is equivalent to
the condition
$$
     \mu_\G(hx)=\mu_\G(x),\qquad
     \mu_\HG(gx)=\mu_\HG(x)
$$
for $g\in\G$, $h\in\HG$, and $x\in M$.
Let $(S,\sigma)$ be a compact symplectic manifold 
and $\pi_S:Q\to S$ be a principal $\HG$-bundle.
We assume that $Q$ is equipped with a connection
$B\in\Aa(Q)\subset\Om^1(Q,\h)$ with {\it nonpositive curvature}. 
This means that there exists a $\sigma$-compatible
almost complex structure $J_S\in\Jj(S,\sigma)$ such that
\begin{equation}\label{eq:negative}
     d\pi(q)w'=J_Sd\pi(q)w\qquad\IMP\qquad
     \inner{F_B(w,w')}{\mu_{\HG}(x)}\le 0.
\end{equation}
Then the manifold
$$
     \tilde M := Q\times_\HG M
$$
carries a symplectic form $\tilde\om\in\Om^2(\tilde M)$ 
whose pullback under the projection
$\pi:Q\times M\to\tilde M$ is given by 
$$
     \pi^*\tilde\om 
     := \pi_2^*\om - d\inner{B}{\mu_{\HG}} + \pi_1^*\sigma,
$$
where $\pi_1:Q\times M\to S$ and $\pi_2:Q\times M\to M$
denote the obvious projections. A moment map 
$\tilde\mu:\tilde M\to\g$ for the obvious $\G$-action on
$\tilde M$ is given by 
$$
     \tilde\mu([q,x]) := \mu_\G(x)
$$
for $q\in Q$ and $x\in M$. Here $[q,x]=[qh,h^{-1}x]$
denotes the equivalence class of the pair $(q,x)$
in $Q\times_\HG M$. Note that if $\mu_\G$ is proper
then so is $\tilde\mu$. Note also that every
$\HG$-invariant and $\om$-compatible 
almost complex structure $J\in\Jj_{\HG}(M,\om)$
induces an almost complex $\tilde J\in\Jj(\tilde M,\tilde\om)$
which acts by $J_S$ on the horizontal
subbundle and by $J$ on the vertical subbundle
of $T\tilde M$.  If $J$ is invariant under both 
$\G$ and $\HG$, then $\tilde J$ in invariant under the 
$\G$-action on $\tilde M$. 

Now suppose that $(J,f)$ is a convex structure 
for the $\G$-action on $M$ (as in Definition~\ref{def:convex}).
Suppose also that $J$ and $f$ are $\HG$-invariant. 
Then the above almost complex structure $\tilde J$
and the function $\tilde f:\tilde M\to[0,\infty)$ given by 
$$
     \tilde f([q,x]) := f(x)
$$
define a convex structure for the $\G$-action on $\tilde M$.  
To see this, note that the gradient $\nabla{\tilde f}$ is given by 
$\nabla{\tilde f}([q,x]) = [0,\nabla f(x)]$. 
Let $\phi_t:M\to M$ denote the (upward)
gradient flow of $f$. Then the gradient flow of 
$\tilde f$ is given by 
$$
     \tilde\phi_t([q,x])=[q,\phi_t(x)].
$$
for $q\in Q$ and $x\in M$. Hence 
$
     d\tilde\phi_t([q,x])[w,\xi]
     = [w,d\phi_t(x)\xi].
$
In particular, the image of a horizontal tangent vector 
$\tilde\xi_0:=[w,-Y_{B_q(w)}(x)]$ under $d\tilde\phi_t([q,x])$
is the horizontal vector
$$
     \tilde\xi_t := [w,-Y_{B_q(w)}(\phi_t(x))]
     \in T_{[q,\phi_t(x)]}\tilde M,
$$
and so
$$
     |\tilde\xi_t|^2
     = |d\pi_S(q)w|^2 
       + |d\phi_t(x)Y_{B_q(w)}|^2
       - \inner{F_B(w,w')}{\mu_{\HG}(x)},
$$
where $w'\in T_qQ$ satisfies $d\pi_S(q)w'=J_Sd\pi_S(q)w$. 
Here we have used the fact that $\phi_t$ commutes with the 
action of $\HG$ and so $\mu_\HG\circ\phi_t=\mu_\HG$.
It follows from the hypotheses on $(J,f)$ that the function
$t\mapsto|\tilde\xi_t|$ is nondecreasing
whenever $\tilde f([q,\phi_t(x)])$ is sufficiently large.
\end{example}


\section{Compactness and Regularity}\label{sec:CR}


\subsection{Regularity}\label{sec:regular}

The next theorem asserts that every weak solution
of equations~(\ref{eq:vortex-ham}) 
is gauge equivalent to a strong solution. 
For an integer $\ell\ge 1$ we denote by 
$$
     \Jj^\ell=\Jj^\ell(\Sigma,M,\om,\mu)
$$
the space of almost complex structures of class $C^\ell$
and by 
$$
     \Hh^\ell=\Hh^\ell(\Sigma,M,\om,\mu)
$$ 
the space of Hamiltonian
perturbations of class $C^\ell$ (see Section~\ref{sec:H}).
Thus $\Hh^\ell$ is the vector space 
of $\G$-equivariant $\C^\ell$-sections of the 
vector bundle $T^*\Sigma\times M\to\Sigma\times M$. 
For $\ell=\infty$ we write 
$
     \Jj^\infty=:\Jj
$
and
$
     \Hh^\infty=:\Hh.
$ 
Consider the symplectic fibre bundle 
$$
     \tilde M:=P\times_\G M\to\Sigma
$$
with fibres diffeomorphic to $M$. 
There is a one-to-one correspondence between
sections $\tilde u:\Sigma\to\tilde M$ and 
$\G$-equivariant functions $u:P\to M$ via
$$
     \tilde u\circ\pi(p)=[p,u(p)]
$$
for $p\in P$. For a positive integer $k$ and a constant 
$p>2$ we denote by $W^{k,p}_\G(P,M)$ the Banach manifold
of all continuous $\G$-equivariant functions
$u:P\to M$ such that the corresponding section
$\tilde u:\Sigma\to\tilde M$ is of class $W^{k,p}$.

\begin{theorem}\label{thm:regular}
Fix a constant $p>2$ and let $\ell$ be either
a positive integer or be equal to $\infty$. 
Let $J\in\Jj^\ell$ and $H\in\Hh^{\ell+1}$. 
If $u\in W^{1,p}_\G(P,M)$ and 
$A\in\Aa^{1,p}(P)$ satisfy~(\ref{eq:vortex-ham})
then there exists a gauge transformation
$g\in\Gg^{2,p}(P)$ such that
$g^{-1}u$ and $g^*A$ are of class $W^{\ell+1,p}$.
For $\ell=\infty$ this means that $g^{-1}u$ and $g^*A$
are smooth. 
\end{theorem}

\begin{proof}
Let 
$
     (u,A)\in 
     W^{1,p}(\Sigma,P\times_\G M)\times\Aa^{1,p}(P)
$
be a solution of~(\ref{eq:vortex-ham}).
Assume first that there exists a smooth 
connection $A_0\in\Aa(P)$ such that 
\begin{equation}\label{eq:slice}
     d_{A_0}^*(A-A_0) = 0.
\end{equation}
Under this assumption we shall prove that the pair 
$(u,A)$ is of class $W^{\ell+1,p}$.  Denote
$
     \alpha:=A-A_0 \in W^{1,p}(\Sigma,T^*\Sigma\otimes\g_P).
$
Then the first equation 
in~(\ref{eq:vortex-ham}) has the form
\begin{equation}\label{eq:reg1}
     \bar\p_{J,A_0}(u)
     =  - \bigl(X_\alpha(u) + X_H(u)\bigr)^{0,1},
\end{equation}
where the $(0,1)$-part of the $1$-form 
$
     X_\alpha(u) - X_H(u)
$
on $\Sigma$ with values in $u^*TM/\G$
is understood with respect to $J_u$.
The second equation in~(\ref{eq:vortex-ham})
and~(\ref{eq:slice}) together have the form
\begin{equation}\label{eq:reg2}
    d_{A_0}\alpha 
    = - F_{A_0} - \frac{1}{2}[\alpha\wedge\alpha]
    + (\tau-\mu(u))\,\dvol_\Sigma,\qquad
    d_{A_0}^*\alpha = 0.
\end{equation}
We prove by induction that $u$ and $A$ are of class
$W^{k,p}$ for every integer $k\le\ell+1$. 
For $k=1$ this holds by assumption.
If $u$ and $A$ are of class $W^{k,p}$ for some $k\le\ell$
then, by~(\ref{eq:reg2}), $d_{A_0}\alpha$ and $d_{A_0}^*\alpha$ 
are of class $W^{k,p}$ and hence $\alpha$ is of class
$W^{k+1,p}$.  Moreover, by~(\ref{eq:reg1}),
$\bar\p_{J,A_0}(u)$ is of class $W^{k,p}$
and the complex structure $J_u$ on the bundle
$u^*TM/\G$ is also of class $W^{k,p}$.
Hence $u$ is of class $W^{k+1,p}$
(see~\cite[Proposition~B.4.7]{MS1}).  
This completes the induction.  
Hence the pair $(u,A)$ is of class $W^{\ell+1,p}$,
and is smooth in the case $\ell=\infty$.

Thus we have proved the theorem under the assumption that $A$
satisfies~(\ref{eq:slice}) for some smooth connection $A_0$.  
In general, it follows from the local slice theorem
(see Theorem~\ref{thm:slice}) that there exists
a smooth connection $A_0$ and a gauge transformation
$g\in\Gg^{2,p}(P)$ such that 
$$
     d_A^*(g_*A_0-A)=0.
$$
Then $g^*A$ satisfies~(\ref{eq:slice}) and hence 
the pair $(g^{-1}u,g^*A)$ is of class $W^{\ell+1,p}$,
and is smooth in the case $\ell=\infty$.
\end{proof}


\subsection{Compactness with bounded derivatives}
\label{sec:compact}

In this section we prove a compactness result 
for solutions of~(\ref{eq:vortex-ham}) with values
in a fixed compact subset of $M$ under the hypothesis
that the first derivatives satisfy a uniform 
$L^\infty$-bound.  We assume that $\om_\nu\in\Om^2(M)$
is a sequence of symplectic forms on $M$ converging 
uniformly with all derivatives to a symplectic form $\om$ 
and that
$
     \mu_\nu:M\to\g
$
is a sequence of moment maps
(corresponding to a sequence of $\om_\nu$-Hamiltonian 
$\G$-actions on $M$) that converges uniformly with 
all derivatives to the moment map $\mu$. 
We assume that $\om_\nu$ agrees with $\om$ and 
$\mu_\nu$ agrees with $\mu$ outside of a compact set.
We assume further that $\dvol_{\Sigma,\nu}$ is 
a sequence of volume forms on $\Sigma$ converging
in the $\Cinf$-topology to $\dvol_\Sigma$
and $J_{\Sigma,\nu}$ is a sequence of complex structures
on $\Sigma$ converging in the $\Cinf$-topology
to $J_\Sigma$. 

\begin{theorem}\label{thm:compact}
Let $\ell$ be either a positive integer 
or be equal to $\infty$. Suppose that
$$
    (J_\nu,H_\nu)
    \in\Jj^\ell(\Sigma,M,\om_\nu,\mu_\nu)
    \times\Hh^{\ell+1}(\Sigma,M,\om_\nu,\mu_\nu)
$$
is a sequence such that $J_\nu$ 
converges to $J\in\Jj^\ell(\Sigma,M,\om,\mu)$ 
in the $C^\ell$-norm on every compact set, 
and that $H_\nu$ converges to $H\in\Hh^{\ell+1}(\Sigma,M,\om,\mu)$ 
in the $C^{\ell+1}$-norm on every compact set.
Suppose further that $\tau_\nu\in Z(\g)$ converges to $\tau$. 
For every $\nu$ let $(u_\nu,A_\nu)\in W^{1,p}_\G(P,M)\times\Aa^{1,p}(P)$
be a solution of~(\ref{eq:vortex-ham})
with $(\mu,J_\Sigma,\dvol_\Sigma,J,H,\tau)$ replaced by 
$(\mu_\nu,J_{\Sigma,\nu},\dvol_{\Sigma,\nu},J_\nu,H_\nu,\tau_\nu)$.
Suppose that there exist a constant $c>0$ 
and a compact set $K\subset M$ such that
$$
     u_\nu(P)\subset K,\qquad 
     \left\|d_{A_\nu}u_\nu\right\|_{L^\infty}\le c
$$
for every $\nu$. Then there exists a sequence of 
gauge transformations $g_\nu\in\Gg^{2,p}(P)$
such that the sequence 
$(g_\nu^{-1}u_\nu,g_\nu^*A_\nu)$
has a $C^\ell$-convergent subsequence. 
\end{theorem}

\begin{lemma}\label{le:estimate}
Fix positive integers $k$ and $n$, a real number $p>2$,
an open set $U\subset\C$, and a compact subset $K\subset U$.
Let $\Jj_n\subset\R^{2n\times 2n}$ denote the set 
of complex structures on $\R^{2n}$.
Then for every constant $c_0>0$ there exists 
a constant $c=c(c_0,K,U,n,k,p)>0$ 
such that the following holds.  
If $J\in W^{k,p}(U,\Jj_n)$ satisfies
$$
     \left\|J\right\|_{W^{k,p}(U)}\le c_0
$$
then every function $u\in W^{k+1,p}(U,\R^{2n})$
satisfies the inequality
$$
     \left\|u\right\|_{W^{k+1,p}(K)}
     \le c\left(
         \left\|\p_su+J\p_tu\right\|_{W^{k,p}(U)}
         + \left\|u\right\|_{W^{k,p}(U)}
         + \left\|u\right\|_{W^{1,\infty}(U)}
         \right).
$$
\end{lemma}

\begin{proof}
We argue by contradiction. Suppose that there exist sequences
$
     J_\nu\in W^{k,p}(U,\Jj_n)
$
and
$
     u_\nu\in W^{k+1,p}(U,\R^{2n}) 
$
such that
$$
     \left\|J_\nu\right\|_{W^{k,p}(U)}\le c_0,\qquad
     \left\|u_\nu\right\|_{W^{k+1,p}(K)}\to\infty
$$
and
$$
     \left\|\p_su_\nu+J_\nu\p_tu_\nu\right\|_{W^{k,p}(U)}
     + \left\|u_\nu\right\|_{W^{k,p}(U)}
     + \left\|u_\nu\right\|_{W^{1,\infty}(U)}
     \le 1.
$$
Passing to a subsequence, if necessary, 
we may assume that there is a complex structure 
$J_0\in W^{k,p}(U,\Jj_n)$ such that $J_\nu$ 
converges to $J_0$ in the weak $W^{k,p}$-topology 
and in the strong $C^0$-topology. 
Choose a smooth cutoff function 
$
     \beta:U\to[0,1]
$
with compact support such that $\beta|_K\equiv 1$
and define
$$
     v_\nu:=\beta u_\nu.
$$
Then $v_\nu$ is bounded in $W^{k,p}$ and $W^{1,\infty}$
and it satisfies the identity
$$
    \p_sv_\nu+J_0\p_tv_\nu 
    = \beta(\p_su_\nu+J_\nu\p_tu_\nu)
      + (\p_s\beta)u_\nu + (\p_t\beta)J_\nu u_\nu
      + (J_0-J_\nu)\p_tv_\nu.
$$
The elliptic estimate for the operator 
$\p_s+J_0\p_t$ has the form
$$
    \left\|v\right\|_{W^{k+1,p}(U)}
    \le c\left(
    \left\|\p_sv+J_0\p_tv\right\|_{W^{k,p}(U)}
    + \left\|v\right\|_{W^{k,p}(U)}
    \right)
$$
for some constant $c>0$ and every function
$v:U\to\R^{2n}$ with compact support 
(see for example~\cite[Proposition~B.4.7]{MS1}).
Hence
\begin{eqnarray*}
    \left\|v_\nu\right\|_{W^{k+1,p}(U)} 
&\le &
    c\biggl(
    \left\|\beta(\p_su_\nu+J_\nu\p_tu_\nu)\right\|_{W^{k,p}(U)} \\
&&    
    +\,
    \left\|(\p_s\beta)u_\nu + (\p_t\beta)J_\nu u_\nu\right\|_{W^{k,p}(U)} \\
&&
    +\,\left\|(J_0-J_\nu)\p_tv_\nu\right\|_{W^{k,p}(U)}
    +\,\left\|v_\nu\right\|_{W^{k,p}(U)}
    \biggr)  \\
&\le &
    c'\biggl(
    \left\|\p_su_\nu+J_\nu\p_tu_\nu\right\|_{W^{k,p}(U)}
    + \left\|u_\nu\right\|_{W^{k,p}(U)} \\
&&
    +\,
      \left\|J_0-J_\nu\right\|_{W^{k,p}(U)}
      \left\|\p_tv_\nu\right\|_{L^\infty(U)} \\
&&
    +\,
      \left\|J_0-J_\nu\right\|_{L^\infty(U)}
      \left\|v_\nu\right\|_{W^{k+1,p}(U)}
    \biggr)
\end{eqnarray*}
for every $\nu$.
If $\nu$ is sufficiently large then 
$c'\left\|J_0-J_\nu\right\|_{L^\infty(U)}\le 1/2$
and hence 
\begin{eqnarray*}
    \left\|v_\nu\right\|_{W^{k+1,p}(U)}    
&\le &
    2c'
    \left(
     \left\|\p_su_\nu+J_\nu\p_tu_\nu\right\|_{W^{k,p}(U)}
     + \left\|u_\nu\right\|_{W^{k,p}(U)}
    \right) \\
&&
    +\, 4c'c_0\left\|\p_tv_\nu\right\|_{L^\infty(U)}.
\end{eqnarray*}
This contradicts the fact that
$
    \left\|u_\nu\right\|_{W^{k+1,p}(K)}    
$
diverges to $\infty$.
\end{proof}

\begin{proof}[Proof of Theorem~\ref{thm:compact}]
The proof consists of three steps.

\medskip
\noindent{\bf Step~1}
{\it 
Fix a constant $p>2$
and a smooth reference connection $A_0\in\Aa(P)$.
We may assume without loss of generality that
the sequence $A_\nu-A_0$ is bounded in $W^{1,p}$.
}

\medskip
\noindent
Since 
$
      F_{A_\nu} = (\tau_\nu-\mu_\nu(u_\nu))\dvol_{\Sigma,\nu}
$
it follows from the assumptions that 
$$
      \sup_\nu\left\|F_{A_\nu}\right\|_{L^\infty} < \infty.
$$
Hence, by Uhlenbeck compactness~\cite{UHL,WEHRHEIM}, 
there exists a sequence of gauge 
transformations $g_\nu\in\Gg^{2,p}(P)$ such that
the connections $g_\nu^*A_\nu\in\Aa^{1,p}(P)$ 
satisfy a uniform $W^{1,p}$-bound.
Replace the sequence $(u_\nu,A_\nu)$
by $(g_\nu^{-1}u_\nu,g_\nu^*A_\nu)$.

\medskip
\noindent{\bf Step~2}
{\it
We may assume without loss of generality
that $A_\nu$ converges weakly in $W^{1,p}$
and strongly in $C^0$ to a connection
$A\in\Aa^{\ell+1,p}(P)$, that $u_\nu$ converges weakly
in $W^{1,p}$ and strongly in $C^0$ to a
section $u\in W^{\ell+1,p}_\G(P,M)$, and that}
$$
     d_A^*(A_\nu-A) = 0.
$$

\medskip
\noindent
By Step~1, the sequence $A_\nu-A_0$ is bounded
in $W^{1,p}$ and, by assumption, the sequence
$u_\nu$ is bounded in $W^{1,p}$.
Hence, by the theorems of Alaoglu and Rellich, we may assume, 
after passing to a subsequence if necessary, 
that $A_\nu$ converges weakly in $W^{1,p}$ 
and strongly in $C^0$ to a connection $A_\infty\in\Aa^{1,p}(P)$, 
and $u_\nu$ converges weakly in $W^{1,p}$ 
and strongly in $C^0$ to a section $u_\infty\in W^{1,p}(\Sigma,E)$.
Since $\bar\p_{J_\nu,H_\nu,A_\nu}(u_\nu)$
converges weakly in $L^p$ to $\bar\p_{J,H,A_\infty}(u_\infty)$
and $F_{A_\nu}$ converges weakly in $L^p$ to $F_{A_\infty}$
it follows that the limit $(u_\infty,A_\infty)$ 
satisfies~(\ref{eq:vortex-ham}).  
By Theorem~\ref{thm:regular}, there exists 
a gauge transformation $g\in\Gg^{2,p}(P)$
such that the pair
$$
     (u,A) := (g^{-1}u_\infty,g^*A_\infty)
$$ 
is of class $W^{\ell+1,p}$. 
Moreover, $g^*A_\nu$ converges weakly in $W^{1,p}$ 
and strongly in $C^0$ to $A$ 
and $g^{-1}u_\nu$ converges weakly in $W^{1,p}$ 
and strongly in $C^0$ to $u$.
By the local slice theorem
(Theorem~\ref{thm:slice}), 
there exists a sequence of gauge 
transformations $h_\nu$ such that 
$$
    d_A^*(h_\nu^*g^*A_\nu-A)=0
$$
and
$$
    \lim_{\nu\to\infty}\left\|
    h_\nu^*g^*A_\nu-A
    \right\|_{L^p}=0,\qquad
    \sup_\nu\left\|
    h_\nu^*g^*A_\nu-A
    \right\|_{W^{1,p}}
    < \infty.
$$
Passing to a subsequence, if necessary, we may assume
that $h_\nu^*g^*A_\nu$ converges weakly in the
$W^{1,p}$-norm and strongly in the $C^0$-norm. 
The limit is necessarily equal to $A$. 
Moreover, the sequence $h_\nu$ is uniformly 
bounded in the $W^{2,p}$-norm.
Passing to a further subsequence we may assume that
$h_\nu$ converges weakly 
in $W^{2,p}$ and strongly in $W^{1,p}$
to a gauge tranformation $h\in\Gg^{2,p}(P)$.
This gauge transformation satisfies
$$
    h^*A=A
$$
and hence is of class $W^{\ell+2,p}$. It follows 
that $h_\nu^{-1}g^{-1}u_\nu$ converges 
weakly in $W^{1,p}$ and strongly in $C^0$ 
to $h^{-1}u$. Now replace $A_\nu$ by $h_\nu^*g^*A_\nu$,
$u_\nu$ by $h_\nu^{-1}g^{-1}u$, and $u$ by $h^{-1}u$. 

\medskip
\noindent{\bf Step~3}
{\it The sequence $(u_\nu,A_\nu-A)$ 
is bounded in $W^{\ell+1,p}$.}

\medskip
\noindent
We prove, by induction, that $(u_\nu,A_\nu-A)$ 
is bounded in $W^{k,p}$ for $1\le k\le \ell+1$. 
For $k=1$ this was proved in Step~1.
Let 
$
     \alpha_\nu:=A_\nu-A
$
and assume, by induction, 
that the sequence $(u_\nu,\alpha_\nu)$ is bounded
in $W^{k,p}$ for some $k\in\{1,\dots,\ell\}$. 
In local $J_{\Sigma,\nu}$-holomorphic coordinates
on $\Sigma$ and local coordinates on $M$ 
the equation $\bar\p_{J_\nu,H_\nu,A_\nu}(u_\nu)=0$
has the form
$$
    \p_su_\nu+J_\nu(s,t,u_\nu)\p_tu_\nu
    = - X_{F_\nu+\inner{\mu_\nu}{\Phi_\nu}}(u_\nu)
      - J_\nu(s,t,u_\nu) X_{G_\nu+\inner{\mu_\nu}{\Psi_\nu}}(u_\nu)
$$
where $A_\nu=\Phi_\nu\,ds+\Psi_\nu\,dt$ and
$H_\nu=F_\nu\,ds+G_\nu\,dt$.
This local equation holds in an open set $U\subset\C$,
the function $u_\nu:U\to\R^{2n}$ takes values in an 
open set $V\subset\R^{2n}$, the function  
$J_\nu:U\times V\to\Jj_n$ is of class $C^\ell$ 
(with a uniform $C^\ell$-bound),
and the functions $F_\nu,G_\nu:U\times V\to\R$
are of class $C^{\ell+1}$ 
(with uniform $C^{\ell+1}$-bounds).
Since $u_\nu$ is uniformly bounded in $W^{k,p}$
so is the complex structure 
$U\to\Jj_n:s+it\mapsto J_\nu(s,t,u_\nu(s,t))$.
Moreover, the sequences $\Phi_\nu,\Psi_\nu:U\to\g$ 
are bounded in $W^{k,p}$.  Hence the sequence
$\p_su_\nu+J_\nu(s,t,u_\nu)\p_tu_\nu:U\to\R^{2n}$
is bounded in $W^{k,p}$.  By assumption, 
the sequence $u_\nu:U\to V$ is bounded in $W^{1,\infty}$.
Hence, by Lemma~\ref{le:estimate}, the sequence 
$u_\nu$ is bounded in $W^{k+1,p}$. 
Now, by Step~2, we have
$$
    d_A\alpha_\nu 
    = - F_A 
      - \frac{1}{2}[\alpha_\nu\wedge\alpha_\nu]
      + (\tau_\nu-\mu_\nu(u_\nu))\,\dvol_{\Sigma,\nu},\qquad
    d_A^*\alpha_\nu = 0.
$$
Since the sequences $\alpha_\nu$ and $u_\nu$ are 
bounded in $W^{k,p}$ it follows that 
$d_A\alpha_\nu$ is bounded in $W^{k,p}$ and hence
the sequence $\alpha_\nu$ is bounded in $W^{k+1,p}$.
This completes the induction.
Hence, by Rellich's theorem, the sequence
$(u_\nu,A_\nu)$ has a $C^\ell$-convergent 
subsequence.
\end{proof}


\subsection{Bubbling}\label{sec:bubble}

The following theorem removes the hypothesis
of a uniform $L^\infty$ bound on the first 
derivatives in Theorem~\ref{thm:compact}.
The manifold $(M,\om)$ is called
{\bf symplectically aspherical} if 
$$
      \int_{S^2}v^*\om = 0
$$
for every smooth map $v:S^2\to M$.
This implies that there is no nonconstant 
$J$-holomorphic sphere (for any almost
complex structure on $M$ that is tamed by $\om$).

\begin{theorem}\label{thm:bubble}
Suppose that $M$ is symplectically aspherical
and fix a compact subset $K\subset M$. 
Let $\ell$ be either a positive integer 
or be equal to $\infty$. Let $\om_\nu$, $\mu_\nu$,
$\dvol_{\Sigma,\nu}$, $J_{\Sigma,\nu}$, 
$J_\nu$, $H_\nu$, and $\tau_\nu$ be as in 
Theorem~\ref{thm:compact}. For every $\nu$ let
$
     (u_\nu,A_\nu)\in W^{1,p}_\G(P,M)\times\Aa^{1,p}(P)
$
be a solution of~(\ref{eq:vortex-ham})
with $(\mu,J_\Sigma,\dvol_\Sigma,J,H,\tau)$ replaced by 
$(\mu_\nu,J_{\Sigma,\nu},\dvol_{\Sigma,\nu},J_\nu,H_\nu,\tau_\nu)$
such that
$$
     u_\nu(P)\subset K
$$
for every $\nu$ and
\begin{equation}\label{eq:e-bound}
     \sup_\nu\inner{[\om_\nu-\mu_\nu]}{[u_\nu]}
     < \infty. 
\end{equation}
Then there exists sequence of 
gauge transformations $g_\nu\in\Gg^{2,p}(P)$
such that the sequence 
$(g_\nu^{-1}u_\nu,g_\nu^*A_\nu)$
has a $C^\ell$-convergent subsequence. 
\end{theorem}

\begin{proof}
By Theorem~\ref{thm:compact}, it suffices to prove
that
\begin{equation}\label{eq:du-Linfty}
     \sup_\nu\left\|d_{A_\nu}u_\nu\right\|_{L^\infty} < \infty.
\end{equation}
Fix a constant $p>2$.  By Step~1 in the proof of 
Theorem~\ref{thm:compact}, we may assume that the sequence 
$A_\nu-A_0$ satisfies a uniform $W^{1,p}$-bound 
for some (and hence every) smooth connection $A_0$.
Now suppose, by contradiction, that the sequence
$
     \left\|d_{A_\nu}u_\nu\right\|_{L^\infty}
$
is unbounded.  Passing to a subsequence,
if necessary, we may assume that this sequence 
diverges to $\infty$. Choose a sequence $p_\nu\in P$ such that
$$
     c_\nu := \left|d_{A_\nu}u_\nu(p_\nu)\right|
     = \left\|d_{A_\nu}u_\nu\right\|_{L^\infty}\to\infty.
$$
Passing to a subsequence, we may assume without loss
of generality that $p_\nu$ converges. 
Let
$
     p_\infty:=\lim_{\nu\to\infty}p_\nu
$
and
$
     z_\infty:=\pi(p_\infty).
$
Choose a convergent sequence of local 
$J_{\Sigma,\nu}$-holomorphic coordinates 
$s+it$ on $\Sigma$ near $z_\infty$ 
and lift these to a convergent sequence of local sections
of $P$ that pass at the origin through $p_\infty$.
In this local frame equations~(\ref{eq:vortex-ham})
have the form 
$$
\begin{array}{rcl}
     \p_su_\nu+X_{F_\nu+\inner{\mu_\nu}{\Phi_\nu}}(u_\nu)
     + J_\nu(s,t,u_\nu)
     \left(
     \p_tu_\nu+X_{G_\nu+\inner{\mu_\nu}{\Psi_\nu}}(u_\nu)
     \right) 
     &= &0, \\
     \p_s\Psi_\nu - \p_t\Phi_\nu + [\Phi_\nu,\Psi_\nu]
     + \lambda^2(\mu(u_\nu)-\tau_\nu)
     &= & 0.
\end{array}
$$
Now consider the rescaled sequence 
$$
     v_\nu(s,t) := u_\nu(\eps_\nu s,\eps_\nu t),\qquad
     \eps_\nu:=\frac{1}{c_\nu}.
$$
This sequence satisfies
$$
     \p_sv_\nu + J_\nu(\eps_\nu s,\eps_\nu t,v_\nu)\p_tv_\nu 
     = -\eps_\nu w_\nu,
$$
where
\begin{eqnarray*}
     w_\nu
&:=&
     X_{\Phi_\nu(\eps_\nu s,\eps_\nu t)}(v_\nu)
     + X_{F_\nu}(\eps_\nu s,\eps_\nu t,v_\nu)  \\
&&
     +\,J_\nu(\eps_\nu s,\eps_\nu t,v_\nu)
     \biggl(
     X_{\Psi_\nu(\eps_\nu s,\eps_\nu t)}(v_\nu)
     + X_{G_\nu}(\eps_\nu s,\eps_\nu t,v_\nu)
     \biggr).
\end{eqnarray*}
The sequences $\Phi_\nu$ and $\Psi_\nu$
satisfy uniform $W^{1,p}$-bounds and, 
by construction, the sequence $v_\nu$ satisfies 
a uniform $W^{1,\infty}$-bound on every compact set. 
Hence the sequence $w_\nu$ satisfies a 
uniform $W^{1,p}$-bound on every compact set.
By Lemma~\ref{le:estimate},
the sequence $v_\nu$ satisfies a uniform 
$W^{2,p}$-bound on every compact set. 
It follows that $v_\nu$ has a subsequence
which converges strongly in $C^1$ on every compact set. 
The limit is a nonconstant pseudo-holomorphic 
curve $v:\C\to K$ with respect to the 
almost complex structure 
$
     J_\infty:=J(z_\infty,\cdot).
$
We prove that it has finite energy. 
To see this note that, for every $R>0$, 
\begin{eqnarray*}
&&
     E(v;B_R) \\
&&=
     \int_{B_R}\left|\p_sv\right|_{J_\infty}^2 \\
&&=
     \lim_{\nu\to\infty}
     \int_{B_R}\left|
      \p_sv^\nu
      + \eps_\nu X_{\Phi_\nu(\eps_\nu s,\eps_\nu t)}(v_\nu)
      + \eps_\nu X_{F_\nu}(\eps_\nu s,\eps_\nu t,v_\nu)
     \right|_{J_\nu(\eps_\nu s,\eps_\nu t,v_\nu)}^2 \\
&&=
     \lim_{\nu\to\infty}
     \int_{B_{\eps_\nu R}}\left|
      \p_su_\nu 
      + X_{\Phi_\nu}(u_\nu)
      + X_{F_\nu}(s,t,u_\nu)
     \right|_{J_\nu(s,t,u_\nu)}^2 \\
&&\le 
     \limsup_{\nu\to\infty} E_\nu(u_\nu,A_\nu).
\end{eqnarray*}
By~(\ref{eq:e-bound}) and the energy 
identity in Proposition~\ref{prop:energy}, 
the sequence $E(u_\nu,A_\nu)$ is bounded.  
Hence $v$ has finite energy, and hence, 
by the removable singularity theorem
(see for example~\cite[Theorem~4.2.1]{MS1}),
it extends to a nonconstant $J_\infty$-holomorphic
sphere in $K$.  Since $M$ is symplectically 
aspherical such a $J_\infty$-holomorphic sphere does not exist. 
This contradiction proves~(\ref{eq:du-Linfty}).
\end{proof}

Combining Theorem~\ref{thm:bubble} 
with the apriori estimate of Lemma~\ref{le:apriori}
we obtain the following compactness result for 
the moduli space of solutions of~(\ref{eq:vortex-ham}).

\begin{corollary}\label{cor:compact}
Suppose $(M,\om,\mu)$ is symplectically aspherical and admits
a strongly convex structure $(f,J_0)$ as in Section~\ref{sec:convex}.
Let $\ell$ be either a positive integer or be equal to $\infty$. 
Let $\om_\nu$, $\mu_\nu$, $\dvol_{\Sigma,\nu}$, $J_{\Sigma,\nu}$, 
$J_\nu$, $H_\nu$, and $\tau_\nu$ be as in Theorem~\ref{thm:compact}.
Let $c_0:Z(\g)\to(0,\infty)$ be as in~$(C3')$, suppose 
that~(\ref{eq:apriori}) holds with $c_0=c_0(\tau)$, that 
each $J_\nu$ agrees with $J_0$ outside of the compact set
$$
     M_0 := \left\{x\in M\,|\,f(x)\le c_0(\tau)+1\right\},
$$
and that each Hamiltonian $H_\nu$ is supported in $M_0$.
Then, for every sequence $(u_\nu,A_\nu)$
of $W^{1,p}$-solutions of~(\ref{eq:vortex-ham}),
with the tuple $(\mu,J_\Sigma,\dvol_\Sigma,J,H,\tau)$ replaced by
$(\mu_\nu,J_{\Sigma,\nu},\dvol_{\Sigma,\nu},J_\nu,H_\nu,\tau_\nu)$,
such that $u_\nu$ represents a fixed equivariant
homology class there exists a sequence of gauge 
transformations $g_\nu\in\Gg^{2,p}(P)$ such that 
$(g_\nu^{-1}u_\nu,g_\nu^*A_\nu)$ has a 
$C^\ell$-convergent subsequence. 
\end{corollary}

\begin{proof}
By Lemma~\ref{le:apriori}, $u_\nu(P)\subset M_0$ for every $\nu$. 
Hence the result follows from Theorem~\ref{thm:bubble} with $K=M_0$.
\end{proof}


\section{Fredholm theory}\label{sec:trans}

Fix a symplectic $2n$-manifold $(M,\om)$ with a 
Hamiltonian $\G$-action and moment map 
$\mu:M\to\g$, a compact Riemann surface 
$(\Sigma,J_\Sigma,\dvol_\Sigma)$,
an almost complex structure $J\in\Jj$,
a Hamiltonian perturbation $H\in\Hh$,
an equivariant homology class
$
     B \in H_2(M_\G;\Z),
$
and a principal $\G$-bundle $P\to\Sigma$
whose characteristic class $[P]\in H_2(\BG;\Z)$
is the image of $B$ under the 
homomorphism $H_2(M_\G;\Z)\to H_2(\BG;\Z)$. 
In this section we examine the moduli space
\begin{eqnarray*}
     \Tilde{\Mm}_{B,\Sigma}
&:= &
     \Tilde{\Mm}_{B,\Sigma}(\tau;J,H)  \\
&:= &
     \left\{(u,A)\in W^{1,p}_\G(P,M)\times\Aa^{1,p}(P)\,|\,
     (u,A)\mbox{ satisfy }(\ref{eq:vortex-ham}),\,
     [u]=B
     \right\}.
\end{eqnarray*}
The quotient by the action of the gauge group will
be denoted by 
$$
     \Mm_{B,\Sigma}
     := \Mm_{B,\Sigma}(\tau;J,H)
     := \Tilde{\Mm}_{B,\Sigma}(\tau;J,H)/\Gg^{2,p}(P).
$$
In this section we prove that,
for a generic Hamiltonian perturbation, the
subspace $\Mm_{B,\Sigma}^*$ of irreducible solutions 
is a finite dimensional manifold. 


\subsection{Regular and irreducible solutions}\label{sec:irr}

Let $J\in\Jj^\ell$ and $H\in\Hh^\ell$ 
(see Sections~\ref{sec:H} and~\ref{sec:regular} for the notation).
We denote by 
$
        \G_x := \left\{g\in\G\,|\,gx=x\right\}
$
the isotopy subgroup of $x\in M$.

\begin{definition}\label{def:irreducible}
A solution $(u,A)\in\W^{1,p}_\G(P,M)\times\Aa^{1,p}(P)$
of~(\ref{eq:vortex-ham}) is called {\bf regular} if 
\begin{equation}\label{eq:regular}
     d_A\eta=0,\quad L_u\eta=0
     \qquad\IMP\qquad
     \eta=0
\end{equation}
for every $\eta\in W^{2,p}(\Sigma,\g_P)$. 
It is called {\bf irreducible} if there exists a point
$p\in P$ such that 
$$
     \G_{u(p)}=\{\one\},\qquad
     \im\,L_{u(p)}\cap\im\,JL_{u(p)} = \{0\}.
$$
\end{definition}

Note that every irreducible solution is regular. 
Note also that an element $\tau\in Z(\g)$ is a regular value
of the restriction of the moment map 
$$
     \W^{1,p}_\G(P,M)\times\Aa^{1,p}(P)\to L^p(\Sigma,\g_P):
     (u,A)\mapsto *F_A+\mu(u)
$$
to the space of pairs $(u,A)$ such that $\bar\p_{J,H,A}(u)=0$
and $[u]=B$ if and only if every pair 
$(u,A)\in\Tilde{\Mm}_{B,\Sigma}(\tau;J,H)$ is regular. 
The next lemma shows that, if $\tau$ is a regular value 
of $\mu$ then regularity can be achieved by choosing a Riemann
surface with large volume. 

\begin{lemma}\label{le:irreducible}
Let $\tau\in Z(\g)$.

\smallskip
\noindent{\bf (i)}
If $\tau$ is a regular value of $\mu$
then there exists a constant $\delta>0$
with the following significance.  
If $B$, $\Sigma$, and $H$ satisfy
\begin{equation}\label{eq:regular1}
     \frac{\inner{[\om+\tau-\mu]}{B}+\|\Om_H\|}{\Vol(\Sigma)} 
     \le\delta
\end{equation}
then every solution $(u,A)$ of~(\ref{eq:vortex-ham})
(for every $J\in\Jj^\ell$) with $[u]=B$ is regular.

\smallskip
\noindent{\bf (ii)}
If $\G$ acts freely on $\mu^{-1}(\tau)$
then there exists a constant $\delta>0$
with the following significance.  
If $B$, $\Sigma$, and $H$ satisfy~(\ref{eq:regular1})
then every solution $(u,A)$ of~(\ref{eq:vortex-ham})
(for every $J\in\Jj^\ell$) with $[u]=B$ is irreducible.
\end{lemma}

\begin{proof}
Choose $\delta>0$ such that 
$$
     |\mu(x)-\tau|^2 \le \delta
     \qquad\IMP\qquad \ker\,L_x=\{0\}.
$$
Let $(u,A)\in\W^{1,p}_\G(P,M)\times\Aa^{1,p}(P)$, 
be a solution of~(\ref{eq:vortex-ham}) such that $[u]=B$.
Then, by Proposition~\ref{prop:energy}, we have
\begin{eqnarray*}
      \inf_{p\in P}\left|\mu(u(p))-\tau\right|^2
&\le &
      \frac{1}{\Vol(\Sigma)}
      \int_\Sigma\left|\mu(u)-\tau\right|^2\,\dvol_\Sigma  \\
&\le &
      \frac{E(u,A)}{\Vol(\Sigma)}  \\
&\le &
      \frac{\inner{[\om+\tau-\mu]}{B}+\|\Om_H\|}{\Vol(\Sigma)}  \\
&\le &
      \delta.
\end{eqnarray*}
The last inequality follows from~(\ref{eq:regular1}). 
Hence there exists a point $p_0\in P$ such that 
$|\mu(u(p_0))-\tau|^2\le\delta$ and so, 
by definiton of $\delta$, the linear map
$L_{u(p_0)}:\g\to T_{u(p_0)}M$ is injective. 
Now suppose that $\eta\in W^{1,p}(\Sigma,\g_P)$ satisfies
$$
      d_A\eta=0,\qquad L_u\eta=0.  
$$
Then $\eta(p_0)=0$ and hence $\eta\equiv 0$.  This proves~(i).
To prove~(ii) choose $\delta>0$ such that 
$$
     |\mu(x)-\tau|^2 \le \delta
     \qquad\IMP\qquad
     \G_x=\{\one\},\;\;\im\,L_x\cap\im\,J(z,x)L_x = \{0\}
$$
for all $(z,x)\in\Sigma\times M$
and argue as in the proof of~(i). 
\end{proof}

Given $\tau\in Z(\g)$ and $(J,H)\in\Jj\times\Hh$, 
we denote the set of irreducible solutions 
of~(\ref{eq:vortex-ham}) by 
$$
     \Tilde{\Mm}_{B,\Sigma}^*
:= 
     \Tilde{\Mm}_{B,\Sigma}^*(M,\om,\mu,\tau;J,H) 
:= 
     \left\{(u,A)\in\Tilde{\Mm}_{B,\Sigma}\,|\,
     (u,A)\mbox{ is irreducible}\right\}
$$
and the quotient space by 
$$
     \Mm_{B,\Sigma}^*
     := \Mm_{B,\Sigma}^*(M,\om,\mu,\tau;J,H)
     := \Tilde{\Mm}_{B,\Sigma}^*(M,\om,\mu,\tau;J,H)/\Gg^{2,p}(P).
$$

\begin{remark}\label{rmk:regular}\rm
The regularity criterion of Lemma~\ref{le:irreducible}
is useful in certain situations (e.g. for the adiabatic
limit argument in~\cite{GaSa}).  However, the condition
is rather restrictive and in many cases the solutions 
are regular under much more general hypotheses.
For example, in the case of linear torus actions, 
one can consider the element $\tau_0\in\g$ defined by
\begin{equation}\label{eq:tau0}
    \tau_0 := \tau-\frac{1}{\Vol(\Sigma)}\int_\Sigma F_A
\end{equation}
for $A\in\Aa(P)$.  This element is independent 
of the connection $A$. Suppose $M=\C^n$ and $\G$ is abelian 
and acts linearly on $M$. If $\G$ and acts freely
(respectively with finite isotropy) on $\mu^{-1}(\tau_0)$ 
then the gauge group acts freely (respectively with finite 
isotropy) on the space of solutions of~(\ref{eq:vortex-ham})
for every Hamiltonian perturbation.  To see this note
that, for every subgroup $\HG\subset\G$, the set
$$
     M^\HG:=\{x\in M\,|\,\HG\subset\G_x\}
$$
is a linear subspace of $M=\C^n$ and so $\mu(M^H)$ 
is a closed convex cone. Applying this to the subgroup 
$$
     \HG:=\bigcap_{p\in P}\G_{u(p)},
$$
where $*F_A+\mu(u)=\tau$, we find that 
$$
     \tau_0 = \frac{1}{\Vol(\Sigma)}
     \int_\Sigma\mu(u)\,\dvol_\Sigma\in\mu(M^H).
$$
Hence, if $\G$ acts freely on $\mu^{-1}(\tau_0)$,
it follows that $\HG=\{\one\}$ and, if $\tau_0$ is a regular
value of $\mu$, it follows that $\HG$ is finite.  
Since $\HG$ is isomorphic to the isotropy subgroup of 
the pair $(u,A)$, this proves the claim. 
\end{remark}


\subsection{The linearized operator}\label{sec:fredholm}

\subsubsection*{Cauchy-Riemann operators}

Fix an almost complex structure $J\in\Jj$ and a Hamiltonian
perturbation $H\in\Hh$.  We begin with a discussion of 
the Cauchy-Riemann operator on the vector bundle 
$
      E_u := u^*TM/\G\to\Sigma
$
associated to a pair $(u,A)\in\Cinf_\G(P,M)\times\Aa(P)$.
This operator will be denoted by 
$
      D_{u,A}:\Cinf(\Sigma,E_u)\to\Om^{0,1}(\Sigma,E_u).
$
It is obtained by differentiating the first equation 
in~(\ref{eq:vortex-ham}) with respect to $u$ and is given by 
\begin{equation}\label{eq:CR}
     D_{u,A}\xi
     = (\Nabla{H,A}\xi)^{0,1} 
       - \frac12J(\Nabla{\xi}J)\p_{J,H,A}(u)
\end{equation}
for $\xi\in\Cinf(\Sigma,E_u)$.
Here $\nabla$ denotes the Levi-Civita connection
of the metric $\om(\cdot,J\cdot)$ on $M$.
Since $J$ depends on the basepoint $z\in\Sigma$
so does the connection $\nabla$. 
The connection $\Nabla{H,A}$ on $E_u$ is given by 
\begin{equation}\label{eq:nablaAsigma}
     \Nabla{H,A}\xi 
     = \nabla\xi + \Nabla{\xi}X_{H,A}(u),
\end{equation}
where the $1$-form $X_{H,A}:TP\to\Vect(M,\om)$ 
is given by 
$$
     (X_{H,A})_p(v):=X_{A_p(v)}+X_{H_{d\pi(p)v}}
$$
for $v\in T_pP$. 

\begin{remark}\label{rmk:local}\rm
In conformal coordinates $s+it$ on $\Sigma$ the connection
$\Nabla{H,A}$ has the form 
$$
     \Nabla{H,A,s}\xi 
     = \Nabla{s}\xi + \Nabla{\xi}X_\Phi + \Nabla{\xi}X_F,\qquad
     \Nabla{H,A,t}\xi 
     = \Nabla{t}\xi + \Nabla{\xi}X_\Psi + \Nabla{\xi}X_G.
$$
Here $\Om\subset\C$ is an open set, $u:\Om\to M$ is a smooth
function, $\nabla$ denotes the Levi-Civita connection of the metric
$\inner{\cdot}{\cdot}_{s,t}=\om(\cdot,J_{s,t}\cdot)$, and
$$
     A=\Phi\,ds+\Psi\,dt,\qquad H=F\,ds+G\,dt
$$
where $\Phi,\Psi:\Om\to\g$ and $F,G:\Om\times M\to\R$. 
Thus the Cauchy--Riemann operator has the form 
$$
     D_{u,A}\xi=\frac12(\xi'\,ds+J\xi'\,dt),
$$
where 
$$
     \xi' = \Nabla{H,A,s}\xi + J\Nabla{H,A,t}\xi 
            - \frac12J(\Nabla{\xi}J)(v_s-Jv_t)
$$
and 
$$
     v_s:=\p_su+L_u\Phi+X_F(u),\qquad v_t:=\p_tu+L_u\Psi+X_G(u).
$$
The covariant derivatives of $J=J(s,t,u(s,t))$ with respect to the 
connection $\Nabla{H,A}$ are given by 
$$
     \Nabla{H,A,s}J = \Nabla{v_s}J + \p_sJ - \Ll_{X_F}J,\qquad
     \Nabla{H,A,t}J = \Nabla{v_t}J + \p_tJ - \Ll_{X_G}J.
$$
\end{remark}

\begin{remark}\label{rmk:Tabla}\rm
We obtain a Hermitian connection $\Tabla{H,A}$ on $u^*TM$ by the 
formula
$$
     \Tabla{H,A}\xi := \Nabla{H,A}\xi - \frac12J(\Nabla{H,A}J)\xi.
$$
The complex linear part of $D_{u,A}$ is given by 
$\xi\mapsto(\Tabla{H,A}\xi)^{0,1}$ and, moreover, 
$$
     D_{u,A}\xi = (\Tabla{H,A}\xi)^{0,1} 
     + \frac14N(\xi,\p_{J,H,A}(u))
     + \frac12(J(\Ll_{X_H}J-\dot J)\xi)^{0,1}
$$
(see~\cite{GaSa}). 
Here 
$$
     N(\xi_1,\xi_2)=2J(\Nabla{\xi_2}J)\xi_1-J(\Nabla{\xi_1}J)\xi_2
$$
denotes the Nijenhuis tensor of $J=J_z$. 
\end{remark}

\subsubsection*{An abstract setting}

Consider the infinite dimensional Banach manifold 
$$
      \Bb := \Bb^{k,p} := W^{k,p}_\G(P,M)\times\Aa^{k,p}(P).
$$
The tangent space of $\Bb$ at $(u,A)\in\Bb$ is given by 
$$
      T_{(u,A)}\Bb 
      = W^{k,p}(\Sigma,E_u)\times 
        W^{k,p}(\Sigma,T^*\Sigma\otimes\g_P),\qquad
      E_u:=u^*TM/\G.
$$
The almost complex structure $J\in\Jj$ determines
a complex structure $J_u$ on $E_u$ and
hence a vector bundle $\Ee=\Ee^{k-1,p}\to\Bb$ with fibres
$$
      \Ee_{(u,A)} 
      := \Ee_u 
      := W^{k-1,p}(\Sigma,\Lambda^{0,1}_{J_u}T^*\Sigma\otimes E_u)
         \times W^{k-1,p}(\Sigma,\g_P).
$$
The action of the gauge group $\Gg^{2,p}(P)$ 
on $\Bb$ lifts to an action on $\Ee$.  
For every Hamiltonian perturbation $H\in\Hh$
there is a $\Gg^{2,p}(P)$-equivariant section 
$$
      \Ff=\Ff_{H,J}:\Bb\to\Ee
$$
given by 
$$
      \Ff_{H,J}(u,A) := (\bar\p_{J,H,A}(u),*F_A+\mu(u)-\tau).
$$
The space $\Tilde{\Mm}_{B,\Sigma}(M,\om,\mu,\tau;J,H)$ is the zero 
set of this section. 

\subsubsection*{The linearized operator}

The vertical differential of $\Ff$
at a zero $(u,A)$ gives rise to an operator
$$
     \Dd_{u,A}:T_{(u,A)}\Bb^{k,p}\to
     \Ee^{k-1,p}_u\oplus W^{k-1,p}(\Sigma,\g_P)
$$
given by 
\begin{equation}\label{eq:DuA}
     \Dd_{u,A}
     \left(\begin{array}{c}\xi \\ \alpha\end{array}\right)
     = \left(
       \begin{array}{c}
         D_{u,A}\xi + (L_u\alpha)^{0,1} \\
         L_u^*\xi - d_A^*\alpha  \\
         d\mu(u)\xi + *d_A\alpha
       \end{array}
       \right).
\end{equation}
Here the linear map $L_x:\g\to T_xM$ is given by
the infinitesmal action, i.e.
$$
     L_x\eta := X_\eta(x)
$$
for $x\in M$ and $\eta\in\g$, and $L_x^*:T_xM\to\g$
denotes its dual operator with respect to the 
given invariant inner product on $\g$ and the 
inner product $\om(\cdot,J(z,x)\cdot)$ on $T_xM$.
Note that this inner product, in general, depends both 
on $x\in M$ and on the point $z\in\Sigma$.  

\begin{proposition}\label{prop:fredholm}
Assume $J\in\Jj^\ell$ and $H\in\Hh^{\ell+1}$
and let $k\in\{1,\dots,\ell+1\}$ and $p>2$. 
Then the operator 
$
     \Dd_{u,A}:T_{(u,A)}\Bb^{k,p}
     \to\Ee_u^{k-1,p}\oplus W^{k-1,p}(\Sigma,\g_P),
$ 
defined by~(\ref{eq:DuA}), is a Fredholm operator
for every pair $(u,A)\in\Bb^{k,p}$. 
It has real index
$$
     \INDEX\,\Dd_{u,A} 
     = (n-\dim\G)\chi(\Sigma) + 2\inner{c_1^\G(TM)}{[u]},
$$
where $c_1^\G(TM)\in H^2(M_\G;\Z)$ denotes the first Chern 
class of the vertical tangent bundle 
$TM\times_\G\EG\to M\times_\G\EG=M_\G$.
\end{proposition}

\begin{proof}
The operator 
$$
     \Om^1(\Sigma,\g_P)
     \to\Om^0(\Sigma,\g_P)\oplus\Om^0(\Sigma,\g_P):
     \alpha\mapsto (-d_A^*\alpha,*d_A\alpha)
$$
has index $-\chi(\Sigma)\dim\G$ and,
by the Riemann-Roch theorem, 
the Cauchy-Rie\-mann operator 
$
     D_{u,A}:\Cinf(\Sigma,E_u)\to\Om^{0,1}(\Sigma,E_u)
$  
has index 
$
     n\chi(\Sigma) + 2c_1(E_u),
$
where 
$
     c_1(E_u):=\inner{c_1^\G(TM)}{[u]}
$
denotes the first Chern number 
of the complex vector bundle $E_u\to\Sigma$. 
The operator $\Dd_{u,A}$ is a compact perturbation 
of the direct sum of these operators.
\end{proof}

\subsubsection*{The adjoint operator}

The formal $L^2$-adjoint operator
$$
     \Dd_{u,A}^*:\Ee^{k+1,p}_u\oplus W^{k+1,p}(\Sigma,\g_P)
     \to T_{(u,A)}\Bb^{k,p}
$$
is given by 
\begin{equation}\label{eq:DuA*}
     \Dd_{u,A}^*
     \left(\begin{array}{c}\eta \\ \phi \\ \psi\end{array}\right)
     = \left(
       \begin{array}{c}
         D_{u,A}^*\eta + L_u\phi + JL_u\psi  \\
         L_u^*\eta - d_A\phi-*d_A\psi  \\
       \end{array}
       \right)
\end{equation}
for $\eta\in\Om^{0,1}_{J_u}(\Sigma,E)$ and 
$\phi,\psi\in\Om^0(\Sigma,\g_P)$.

\begin{proposition}\label{prop:DD*}
If $u$ and $A$ satisfy~(\ref{eq:vortex-ham})
then
$$
     \Dd_{u,A}\Dd_{u,A}^*
     \left(\begin{array}{c}\eta \\ \phi \\ \psi\end{array}\right)
     = \left(
       \begin{array}{c}
         D_{u,A}D_{u,A}^*\eta + (L_uL_u^*\eta)^{0,1}
         + (D_{u,A}J-JD_{u,A})L_u\psi  \\
         d_A^*d_A\phi + L_u^*L_u\phi \\
         d_A^*d_A\psi + L_u^*L_u\psi 
         + L_u^*(D_{u,A}J-JD_{u,A})^*\eta
       \end{array}
       \right)
$$
for $\eta\in\Om^{0,1}_{J_u}(\Sigma,E_u)$ and 
$\phi,\psi\in\Om^0(\Sigma,\g_P)$.
\end{proposition}

\begin{proof}
We shall abbreviate $D:=D_{u,A}$
and use the identities 
$$
     d_A^*\alpha=-*d_A*\alpha,\qquad 
     *d_Ad_A\phi=[*F_A,\phi],
$$
$$
     d\mu(u)J = L_u^*,\qquad
     d\mu(u)L_u\phi = [\phi,\mu(u)]
$$
for $\alpha\in\Om^1(\Sigma,\g_P)$ and 
$\phi\in\Om^0(\Sigma,\g_P)$.
With this understood we obtain
$$
     \Dd_{u,A}\Dd_{u,A}^*(\eta,\phi,\psi)
     = (\tilde\eta,\tilde\phi,\tilde\psi),
$$
where 
\begin{eqnarray*}
     \tilde\eta 
&= &
     D(D^*\eta+L_u\phi+JL_u\psi)
     + (L_u(L_u^*\eta-d_A\phi-*d_A\psi))^{0,1} \\
&= &
     DD^*\eta + (L_uL_u^*\eta)^{0,1} \\
&&
       +\, DL_u\phi-(L_ud_A\phi)^{0,1}
       + DJL_u\psi-(L_u*d_A\psi)^{0,1}, \\
     \tilde\phi
&= &
     L_u^*(D^*\eta+L_u\phi+JL_u\psi)
     -\,d_A^*(L_u^*\eta-d_A\phi-*d_A\psi) \\
&= &
     d_A^*d_A\phi + L_u^*L_u\phi
       + L_u^*D^*\eta - d_A^*L_u^*\eta
       + [*F_A+\mu(u),\psi],  \\
     \tilde\psi
&= &
     d\mu(u)(D^*\eta+L_u\phi+JL_u\psi)
     + *d_A(L_u^*\eta-d_A\phi-*d_A\psi) \\
&= &
     d_A^*d_A\psi + L_u^*L_u\psi 
       + L_u^*J^*D^*\eta+d_A^**L_u^*\eta
       - [*F_A+\mu(u),\phi].
\end{eqnarray*}
The assertion then follows from the fact that
\begin{equation}\label{eq:JLu}
     J(L_u\alpha)^{0,1} = (L_u*\alpha)^{0,1},\qquad
     L_u^*J^*\eta = -*L_u^*\eta,
\end{equation}
for $\alpha\in\Om^1(\Sigma,\g_P)$ and
$\eta\in\Om^{0,1}_{J_u}(\Sigma,E_u)$,
and
\begin{equation}\label{eq:DLu}
     \bar\p_{J,H,A}(u)=0\qquad\IMP\qquad
     DL_u\phi = (L_ud_A\phi)^{0,1}
\end{equation}
for $\phi\in\Om^0(\Sigma,\g_P)$. 
The first equation in~(\ref{eq:JLu}) follows from 
the fact that 
$$
     *\alpha = -\alpha\circ J_\Sigma
$$
for every $1$-form $\alpha$ on $\Sigma$
(with values in any vector bundle)
and hence 
\begin{eqnarray*}
     (L_u\alpha)^{0,1}
&= & 
     \frac12(L_u\alpha+JL_u(\alpha\circ J_\Sigma))  \\
&= &
     \frac12((L_u*\alpha)\circ J_\Sigma - JL_u*\alpha)  \\
&= &
     -J(L_u*\alpha)^{0,1}.
\end{eqnarray*}
The second equation in~(\ref{eq:JLu}) follows
from the first by duality. 
Next we observe that the operator 
$(u,A)\mapsto\bar\p_{J,H,A}(u)$ is a section of the 
bundle over $\Bb$ with fibres $\Om^{0,1}_{J_u}(\Sigma,E_u)$. 
Its vertical differential at a zero $(u,A)$ is the operator
$
     (\xi,\alpha)\mapsto D_{u,A}\xi+(L_u\alpha)^{0,1}.
$
Since the section $(u,A)\mapsto\bar\p_{J,H,A}(u)$
is equivariant under the action of $\Gg(P)$ it follows
that the pair $(\xi,\alpha)=(L_u\phi,-d_A\phi)$ 
is contained in the kernel of the 
vertical differential at any zero $(u,A)$.   
This proves~(\ref{eq:DLu}).
\end{proof}

\begin{proposition}\label{prop:onto}
Let $\tau\in Z(\g)$. Then the following holds.

\smallskip
\noindent{\bf (i)}
If $(u,A)$ is a regular solution of~(\ref{eq:vortex-ham}) 
then 
$
     \Dd_{u,A}^*(\eta,\phi,\psi)=0
$
implies
$
     \phi=0.
$

\smallskip
\noindent{\bf (ii)}
If $H=0$, $\dot J=0$, $J$ is integrable, 
and $(u,A)$ is a regular solution of~(\ref{eq:vortex-ham}) 
then 
$
     \Dd_{u,A}^*(\eta,\phi,\psi)=0
$
implies
$
     \phi=\psi=0,
$
$
     D_{u,A}^*\eta=0,
$
and
$
     L_u^*\eta=0.
$

\smallskip
\noindent{\bf (iii)}
If $H=0$, $\dot J=0$, $J$ is integrable, 
$(u,A)$ is a regular solution of~(\ref{eq:vortex-ham}), 
and $L_{u(p)}:\g\to T_{u(p)}M$
is onto for some $p\in P$ then $\Dd_{u,A}$ is onto.

\smallskip
\noindent{\bf (iv)}
If $H=0$, $\dot J=0$, $d_Au\equiv0$, $\mu(u)\equiv\tau$, 
$\Sigma=S^2$, and $(u,A)$ is irreducible then $\Dd_{u,A}$ is onto.
\end{proposition}

\begin{proof}
If $u$ and $A$ satisfy~(\ref{eq:vortex-ham}) and
$
     \Dd_{u,A}^*(\eta,\phi,\psi)=0
$
then, by Proposition~\ref{prop:DD*}, $d_A\phi=0$ and $L_u\phi=0$.
Since $(u,A)$ is regular it follows that $\phi=0$. 
This proves~(i). If $u$ and $A$ satisfy~(\ref{eq:vortex}) 
and $J$ is integrable then 
$$
     D_{u,A}J=JD_{u,A}
$$
and, by Proposition~\ref{prop:DD*}, we have
$$
     \Dd_{u,A}\Dd_{u,A}^*
     \left(\begin{array}{c}\eta \\ \phi \\ \psi\end{array}\right)
     = \left(
       \begin{array}{c}
         D_{u,A}D_{u,A}^*\eta + (L_uL_u^*\eta)^{0,1} \\
         d_A^*d_A\phi + L_u^*L_u\phi \\
         d_A^*d_A\psi + L_u^*L_u\psi
       \end{array}
       \right).
$$
This proves~(ii). To prove~(iii) suppose that
$(u,A)$ is a regular solution of~(\ref{eq:vortex}),
$J$ is integrable, and $\Dd_{u,A}^*(\eta,\phi,\psi)=0$. 
Then, by~(ii), $\phi=\psi=0$ and 
$
     D_{u,A}^*\eta=0,
$
$
     L_u^*\eta=0.
$
Since $L_{u(p)}^*$ is injective for some $p\in P$
it follows that $\eta$ vanishes on some open set.
Hence, by unique continuation, $\eta\equiv0$.
Thus we have proved that $\Dd_{u,A}^*$ is injective. 
Hence $\Dd_{u,A}$ has a dense range and hence,
by elliptic regularity, $\Dd_{u,A}$ is onto. 
To prove~(iv) note that, by Remark~\ref{rmk:Tabla},
The operator $D_{u,A}$ is complex linear whenever
$H=0$, $\dot J=0$, and $d_Au=0$. 
Hence $\Dd_{u,A}^*(\eta,\phi,\psi)=0$ implies 
$\phi=\psi=0$ and 
$
     D_{u,A}^*\eta=0,
$
$
     L_u^*\eta=0.
$
Since $\eta$ is a $(0,1)$-form we have
$$
     d\mu(u)\eta=-L_u^*J\eta=L_u^*(\eta\circ J_\Sigma)=0.
$$ 
Since $d_Au=0$ the image of $u$ is an orbit of some point 
$x_0\in M$ under the $\G$-action.  Since $(u,A)$ is irreducible,
$\G_{x_0}=\{1\}$. Hence $\eta$ defines an element of the 
cokernel of the Cauchy--Riemann operator along the constant
function $\bar u\equiv[x_0]:\Sigma\to\bar M=M\dslash\G(\tau)$. 
Since $\Sigma=S^2$, it follows from the Riemann--Roch theorem, 
that $\eta=0$.  
\end{proof}


\subsection{Transversality}\label{sec:trans1}

In this section we establish transversality for the 
irreducible solutions of~(\ref{eq:vortex-ham})
for generic Hamiltonian perturbations.

\begin{definition}\label{def:regular}
A pair 
$
     (J,H)\in\Jj\times\Hh
$
is called {\bf regular 
(for the sixtuple $(B,\Sigma,M,\om,\mu,\tau)$)} 
if the operator $\Dd_{u,A}$ is surjective for every 
$(u,A)\in\Tilde{\Mm}_{B,\Sigma}^*(M,\om,\mu,\tau;J,H)$,
i.e. for every irreducible solution 
of~(\ref{eq:vortex-ham}) representing the class $B$.
\end{definition}

Given an almost complex structure $J\in\Jj$ 
we denote by
$
     \Hreg(\tau,J)\subset\Hh
$
the set of Hamiltonian perturbations $H$ such that
the pair $(J,H)$ is regular.  Given $\tau\in Z(\g)$ and 
$B\in H_2(M_\G;\Z)$ we denote by
$
     \Jreg(\tau,B)\subset\Jj
$
the set of almost complex structures $J$ such that
the pair $(J,H=0)$ is regular for $B$.  If $(J,H)$ is regular
then it follows from Proposition~\ref{prop:fredholm}
and the infinite dimensional implicit function theorem 
that $\Mm_{B,\Sigma}^*(M,\om,\mu,\tau;J,H)$ 
is a finite dimensional smooth manifold of dimension
\begin{equation}\label{eq:dim}
     \dim\Mm_{B,\Sigma}^*(M,\om,\mu,\tau;J,H)
     = (n-\dim\G)\chi(\Sigma) + 2 \inner{c_1^\G(TM)}{B}.
\end{equation}

\begin{theorem}\label{thm:trans}
{\bf (i)}
For every $J\in\Jj$ and every $\tau\in Z(\g)$
the set $\Hreg(\tau,J)$ is a countable intersection
of open and dense subsets of $\Hh$.

\smallskip
\noindent{\bf (ii)}
Assume $B$ is not a torsion class.  Then, for every 
$\tau\in Z(\g)$, the set $\Jreg(\tau,B)$ is a countable 
intersection of open and dense subsets of $\Jj$.
\end{theorem}

\begin{proof}
Fix a sufficiently large integer $\ell$,
a constant $c>0$, a compact set $K\subset M$, 
and a real number $p>2$.
Consider the space of pairs
$(u,A)\in W^{1,p}_\G(P,M)\times\Aa^{1,p}(P)$
that satisfy
\begin{equation}\label{eq:uK1}
     u(P)\subset K,\qquad
     \left\|d_Au\right\|_{L^\infty}\le c
\end{equation}
and, for some $p_0\in P$ and all $\eta_1,\eta_2\in\g$,
\begin{equation}\label{eq:uK2}
     \inf_{\G_x\ne\{\one\}}|u(p_0)-x|\ge \frac{1}{c},\qquad
     \left|\eta_1\right|+\left|\eta_2\right|
     \le c\left|L_{u(p_0)}\eta_1+JL_{u(p_0)}\eta_2\right|.
\end{equation}
Denote
$$
     \Tilde{\Mm}^{c,K}(H) 
     := \left\{(u,A)\in\Tilde{\Mm}_{B,\Sigma}(M,\om,\mu,\tau;J,H)
        \,|\,u\mbox{ and }A\mbox{ satisfy }(\ref{eq:uK1}-\ref{eq:uK2})
        \right\}.
$$
By Theorem~\ref{thm:compact}, the moduli space 
$$
     \Mm^{c,K}(H)
     := \Tilde{\Mm}^{c,K}(H)/\Gg^{2,p}(P)
$$
is compact, and it consists entirely of irreducible
solutions of~(\ref{eq:vortex-ham}). 
We shall examine the universal moduli space
$$
      \Uu^{\ell,c,K}
      := \Tilde{\Uu}^{\ell,c,K}/\Gg^{2,p}(P),
$$
$$
      \Tilde{\Uu}^{\ell,c,K}
      := \left\{(u,A,H)\,|\,
         H\in\Hh^\ell,\,
         (u,A)\in\Tilde{\Mm}^{c,K}(H),\,
         (\ref{eq:uK1}-\ref{eq:uK2})
         \right\}.
$$
We prove that $\Uu^{\ell,c,K}$ is a Banach manifold.  
To see this we must show that the linearized operator
$$
     \Dd_{u,A,H}:\Bb_{(u,A)}^{1,p}\oplus\Hh^\ell
     \to \Ee_u^p\oplus L^p(\Sigma,\g_P)
$$
given by 
\begin{equation}\label{eq:D-H}
     \Dd_{u,A,H}
     \left(\begin{array}{c}
     \xi \\ 
     \alpha  \\
     \hat H
     \end{array}\right)
     =
     \Dd_{u,A}\left(
       \begin{array}{c}
         \xi \\ \alpha
       \end{array}
       \right)
     + \left(
       \begin{array}{c}
         (X_{\hat H}(u))^{0,1} \\
         0  \\
         0
       \end{array}
       \right),
\end{equation}
is surjective for every triple 
$(u,A,H)\in\Tilde{\Uu}^{\ell,c,K}$.
The proof of Theorem~\ref{thm:regular} shows that
we may assume, without loss of generality, that 
$u$ and $A$ are of class $W^{\ell,p}$. 
By Proposition~\ref{prop:fredholm},
$\Dd_{u,A}$ is a Fredholm operator
and hence it suffices to prove that 
$\Dd_{u,A,H}$ has a dense range. 
Let $1/p+1/q = 1$ and assume that the triple 
$$
    (\eta,\phi,\psi)
    \in L^q(\Sigma,\Lambda^{0,1}_{J_u}T^*\Sigma\otimes E_u) 
      \times
      L^q(\Sigma,\g_P) 
      \times  
      L^q(\Sigma,\g_P)
$$
is $L^2$ orthogonal to the image of $\Dd_{u,A,H}$.
Then, in particular, the triple $(\eta,\phi,\psi)$
is $L^2$-orthogonal to the image of $\Dd_{u,A}$.
Since $u$ and $A$ are of class $W^{\ell,p}$ 
and $H$ is of class $C^\ell$ it follows from elliptic
regularity that $\eta$, $\phi$, and $\psi$ are of class
$W^{\ell,p}$ (and hence of class $C^{\ell-1}$) and
\begin{equation}\label{eq:trans1}
     \Dd_{u,A}^*(\eta,\phi,\psi)=0.
\end{equation}
Moreover,
\begin{equation}\label{eq:trans2}
     \int_\Sigma \inner{\eta}{X_{\hat H}(u)}\,\dvol_\Sigma=0
\end{equation}
for every $\hat H\in\Hh^\ell$.
Since $u$ is irreducible, there exists a $p_0\in P$ such that 
\begin{equation}\label{eq:p0}
     \G_{u(p_0)}=\{\one\},\qquad
     \im\,L_{u(p_0)}\cap JL_{u(p_0)}=\{0\}.
\end{equation}
We prove that the linear map 
$\eta_{p_0}:T_{p_0}P\to T_{u(p_0)}M$
is equal to zero.  Suppose otherwise
that $\eta_{p_0}\ne 0$ and choose 
$v_0\in T_{p_0}P$ such that 
$
     \eta_{p_0}(v_0)\ne 0.
$
Since the linear map $\eta_{p_0}:T_{p_0}P\to T_{u(p_0)}M$
vanishes on $p_0\cdot\g$ we have 
$v_0\notin p_0\cdot\g$.  Denote 
$$
     x_0:=u(p_0),\qquad
     z_0:=\pi(p_0),\qquad
     \zeta_0:=d\pi(p_0)v_0\ne 0,
$$
and choose $v_1\in T_{p_0}P$ such that
$$
     \zeta_1 := d\pi(p_0)v_1 = -J_\Sigma\zeta_0,\qquad
     \eta_{p_0}(v_1) = J(z_0,x_0)\eta_{p_0}(v_0).
$$
The last identity follows from the fact that 
$\eta$ is a $(0,1)$-form.  By~(\ref{eq:p0}), 
$\eta_{p_0}(v_0)$ and $\eta_{p_0}(v_1)$
cannot both lie in the image of the map
$L_{x_0}:\g\to T_{x_0}M$. Let us assume, without 
loss of generality, that 
$$
     \eta_{p_0}(v_1)\notin\im\,L_{x_0}.
$$
Since $\G_{x_0}=\{\one\}$ there 
exists a $\G$-invariant neighbourhood 
$U_0\subset M$ of $x_0$ such that 
$\G_x=\{\one\}$ for every $x\in U_0$. 
Since $\G_x=\{\one\}$ for every $x\in U_0$ 
there exists a smooth $\G$-invariant function
$\hat H_0:M\to\R$, supported in $U_0$, 
such that 
$$
     d\hat H_0(x_0)\eta_{p_0}(v_1) > 0.
$$
Hence 
\begin{eqnarray*}
     \inner{X_{\hat H_0}(u(p_0))}{\eta_{p_0}(v_0)}
&= &
     \om(X_{\hat H_0}(x_0),J(z_0,x_0)\eta_{p_0}(v_0))  \\
&= &
     \om(X_{\hat H_0}(x_0),\eta_{p_0}(v_1))  \\
&= &
     d\hat H_0(x_0)\eta_{p_0}(v_1)  \\
&> &
     0.
\end{eqnarray*}
Now let $\zeta:\Sigma\to T\Sigma$ be a vector field
such that 
$
     \zeta(z_0)=\zeta_0.
$ 
Choose a neighbourhood $V_0$ of $z_0$ such that,
for every $p\in P$ and every $v\in T_pP$,
$$
     \pi(p)\in V_0,\quad d\pi(p)v=\zeta(\pi(p))
     \qquad\IMP\qquad
     \inner{X_{\hat H_0}(u(p))}{\eta_p(v)} > 0
$$
Now choose a cutoff function $\beta:\Sigma\to[0,1]$
which is supported in $V_0$ and satisfies
$\beta(z_0)=1$. Define $\hat H\in\Hh$ by the conditions
$$
     \hat H_{\zeta(z)}=\beta(z)\hat H_0,\qquad
     \hat H_{J_\Sigma\zeta(z)}=0.
$$
Then $\inner{\eta}{X_{\hat H}(u)}>0$ at the point $z_0$ and
$\inner{\eta}{X_{\hat H}(u)}\ge 0$ everywhere. 
Hence 
$$
     \int_\Sigma 
     \inner{\eta}{X_{\hat H}(u)}\,\dvol_\Sigma
     >0
$$
in contradiction to~(\ref{eq:trans2}).  
Thus we have proved that, for every $p\in P$,
$$
     \G_{u(p)}=\{\one\}\qquad\IMP\qquad \eta(p)=0.
$$
Moreover, by~(\ref{eq:trans1}) 
and Proposition~\ref{prop:onto},
we have that
$
     \phi=0.
$
Since $\ker L_{u(p)}=0$ whenever $\G_{u(p)}=\{\one\}$, 
it follows from~(\ref{eq:DuA*}) and~(\ref{eq:trans2})
that 
$$
     \G_{u(p)}=\{\one\}\qquad\IMP\qquad \psi(p)=0.
$$
Hence $\eta$ and $\psi$ vanish simultaneously
on some open subset of $P$.  Since 
$$
     \Dd_{u,A}^*(\eta,0,\psi)=0
$$
it follows by unique continuation for first order 
elliptic operators that $\eta=0$ and $\psi=0$.  

Thus we have proved that the operator 
$\Dd_{u,A,H}$ has a dense range for every
$(u,A,H)\in\Tilde{\Uu}^{\ell,c,K}$, as claimed. 
Since $\Dd_{u,A}$ is Fredholm,  
the operator $\Dd_{u,A,H}$ has a right inverse. 
Hence, by the implicit function theorem, 
$\Uu^{\ell,c,K}$ is a separable Banach manifold
of class $C^{\ell-1}$. The projection
$$
     \pi^{\ell,c,K}:\Uu^{\ell,c,K}\to\Hh^\ell
$$
is a Fredholm map of index 
$
     (n-\dim\G)\chi(\Sigma)+2\inner{c_1^\G(TM)}{B}.
$ 
Hence, for $\ell$ sufficiently large, it follows from
the Sard-Smale theorem, that the set 
$\Hreg^{\ell,c,K}(\tau,J)\subset\Hh^\ell$
of regular values of $\pi^{\ell,c,K}$ is dense in $\Hh^\ell$. 
Moreover, the moduli space 
$$
     \Mm^{\ell,c,K}(H) 
     := \left\{[u,A]\,|\,([u,A],H)\in\Uu^{\ell,c,K}\right\}
$$
is compact for every $H$. Hence the set 
$\Hreg^{\ell,c,K}(\tau,J)$ is open and dense in $\Hh^\ell$.
Hence the set
$
     \Hreg^{c,K}(\tau,J) = \Hreg^{\ell,c,K}(\tau,J)\cap\Hh
$
is dense in $\Hh^\ell$ and hence is also dense in $\Hh$.
That it is open follows again from compactness.
Hence the set
$$
     \Hreg(\tau,J) = \bigcap_{c,K}\Hreg^{c,K}(\tau,J)
$$
is a countable intersection of open and dense 
subsets of $\Hh$. This proves~(i).

We sketch the proof of~(ii). 
Assume $B$ is not a torsion class.  
Then, by Lemma~\ref{le:horizontal}, every 
solution $(u,A)$ of~(\ref{eq:vortex-ham}) with $H=0$ 
representing the class $B$ satisfies $d_Au\not\equiv0$
and so, by Lemma~\ref{le:horizontal1}, $d_Au\ne 0$ 
almost everywhere.  Hence,
for every irreducible solution $(u,A)$ 
of~(\ref{eq:vortex-ham}) with $H=0$ that represents
the class $B$, there exists a point $p_0\in P$ such that
$$
     \G_{u(p_0)}=\{\one\},\qquad
     \im\,L_{u(p_0)}\cap\im\,JL_{u(p_0)}=\{0\},\qquad
     d_{A}(u(p_0))\ne 0.
$$
With this understood the proof of assertion~(ii)
is almost word by word the same as that of~(i)
and will be omitted. 
\end{proof}


\subsection{Cobordisms}\label{sec:cobordism}

Let $(\tau_\lambda,J_\lambda)\in Z(\g)\times\Jj$
and $H_\lambda\in\Hreg(\tau_\lambda,J_\lambda)$
for $\lambda=0,1$ .  
For every smooth homotopy 
$$
     \{\tau_\lambda,J_\lambda,H_\lambda\}_{0\le\lambda\le 1}
     \in Z(\g)\times\Jj\times\Hh
$$
from $(\tau_0,J_0,H_0)$ to $(\tau_1,J_1,H_1)$ we consider the 
space
$$
     \Ww^*_{B,\Sigma}
     := \Ww^*_{B,\Sigma}(\{\tau_\lambda,J_\lambda,H_\lambda\}_\lambda)
     := \bigcup_{0\le\lambda\le 1}
        \{\lambda\}\cup
        \Mm^*_{B,\Sigma}(\tau_\lambda;J_\lambda,H_\lambda).
$$

\begin{definition}\label{def:homotopy}
A homotopy $\{\tau_\lambda,J_\lambda,H_\lambda\}_\lambda$
is called {\bf regular (for the tuple $(B,\Sigma,M,\om,\mu)$)}
if $H_\lambda\in\Hreg(\tau_\lambda,J_\lambda)$ for $\lambda=0,1$
and, for every triple $(\lambda,[u,A])\in\Ww^*_{B,\Sigma}$,
we have
$$
     \Ee_u
     = \im\,\Dd_{u,A}
       + \SPAN^\R \zeta_{\lambda,u,A},
$$
where
\begin{equation}\label{eq:cobordism}
     \zeta_{\lambda,u,A} 
     := \left(
       \begin{array}{c}
         \Nabla{\lambda}(\bar\p_{J_\lambda,H_\lambda,A}(u)) \\
         0 \\
         -\p_\lambda\tau_\lambda
       \end{array}
       \right).
\end{equation}
\end{definition}

Here the expression 
$\Nabla{\lambda}(\bar\p_{J_\lambda,H_\lambda,A}(u))$
is independent of the (Hermitian) connection used to define it.
If $\{\tau_\lambda,J_\lambda,H_\lambda\}_\lambda$
is a regular homotopy then the moduli space
$\Ww_{B,\Sigma}^*(\{\tau_\lambda,J_\lambda,H_\lambda\}_\lambda)$
is a smooth finite dimensional manifold with boundary
$$
     \p\Ww_{B,\Sigma}^*(\{\tau_\lambda,J_\lambda,H_\lambda\}_\lambda)
     = \{0\}\times\Mm^*_{B,\Sigma}(\tau_0;J_0,H_0)
       \cup \{1\}\times\Mm^*_{B,\Sigma}(\tau_1;J_1,H_1).
$$
Let $\Hh^{[0,1]}(H_0,H_1)$ denote the space of smooth paths
$[0,1]\to\Hh:\lambda\mapsto H_\lambda$ with fixed 
endpoints $H_0$ and $H_1$. Given a homotopy
$\{\tau_\lambda,J_\lambda\}_\lambda$ from
$(\tau_0,J_0)$ to $(\tau_1,J_1)$ denote 
by 
$
     \Hreg^{[0,1]}(H_0,H_1;\{\tau_\lambda,J_\lambda\}_\lambda)
     \subset \Hh^{[0,1]}(H_0,H_1)
$
the set of all smooth homotopies $\{H_\lambda\}_\lambda$
from $H_0$ to $H_1$ such that the triple 
$\{\tau_\lambda,J_\lambda,H_\lambda\}_\lambda$
is regular. The next theorem asserts the set
of regular homotopies is of the second category
in the sense of Baire. 

\begin{theorem}\label{thm:cobordism}
Assume $\G=T$ is a torus.
Let $\{\tau_\lambda,J_\lambda\}_{0\le\lambda\le1}$
be a smooth homotopy in $\Cinf(\Sigma,Z(\g))\times\Jj$
and let $H_\lambda\in\Hreg(\tau_\lambda,J_\lambda)$ 
for $\lambda=0,1$. Then the set 
$\Hreg^{[0,1]}(H_0,H_1;\{\tau_\lambda,J_\lambda\}_\lambda)$
is a countable intersection of open and dense subsets
of $\Hh^{[0,1]}(H_0,H_1)$. 
\end{theorem}

\begin{proof}
The proof is similar to that of Theorem~\ref{thm:trans}
and we only sketch the main points.  
Denote by $\hat\Hh^\ell$ the set of all $C^\ell$
homotopies from $H_0$ to $H_1$.
Fix a constant $c>0$ and a compact set $K\subset M$
and consider the universal moduli space
of all gauge equivalence classes of quadruples 
$$
     (\lambda,u,A,\{H_\lambda\}_\lambda)
     \in [0,1]\times W^{1,p}_\G(P,M)\times\Aa^{1,p}(P)\times\hat\Hh^\ell
$$
such that 
$
     (u,A)\in\Tilde{\Mm}^*_{B,\Sigma}(\tau_\lambda;J_\lambda,H_\lambda)
$
and $u$ and $A$ satisfy~(\ref{eq:uK1}-\ref{eq:uK2}).  
The proof of Theorem~\ref{thm:trans} shows that
this space is a separable Banach manifold.
The projection 
$$
     (\lambda,[u,A],\{H_\lambda\}_\lambda)
     \mapsto \{H_\lambda\}_\lambda
$$
is then a Fredholm map and a smooth homotopy 
$\{H_\lambda\}_\lambda$ is a regular value of 
this projection for every triple $(K,c,\ell)$
if and only if 
$$
     \{H_\lambda\}_\lambda\in
     \Hreg^{[0,1]}(H_0,H_1;\{\tau_\lambda,J_\lambda\}_\lambda).
$$
Hence the result follows from the Sard-Smale theorem.
\end{proof}


\subsection{Orientation}\label{sec:orient}

In this subsection we shall prove that the moduli spaces
$\Mm_{B,\Sigma}^*(M,\om,\mu;J,H)$ carry natural 
orientations. Consider the determinant line bundle
$$
     \det(\Dd)\to\Bb=\Cinf_\G(P,M)\times\Aa(P)
$$ 
whose fibre over $(u,A)$ is the $1$-dimensional
real vector space
$$
     \det(\Dd_{u,A})
     := \Lambda^\MAX(\ker \Dd_{u,A})
        \otimes\Lambda^\MAX(\ker\Dd_{u,A}^*)
$$
(See, for example,~\cite[Appendix~A]{SAL2} for a detailed
exposition.) Here 
$$
     \Dd_{u,A}=\Dd_{u,A}:T_{(u,A)}\Bb\to\Ee_u
$$
denotes the Fredholm operator 
given by equation~(\ref{eq:DuA}).
The orientation of the moduli spaces 
is an immediate consequence of the
following proposition. The proof is reminiscent
of the arguments in~\cite{MS1,SAL2}.

\begin{proposition}\label{prop:orient}
The determinant line bundle $\det(\Dd)\to\Bb$ admits
a natural $\Gg$-invariant orientation. 
\end{proposition}

\begin{proof}
The tangent space 
$$
     T_{(u,A)}\Bb = \Cinf(\Sigma,E_u)\oplus\Om^1(\Sigma,\g_P)
$$
admits a natural complex structure given by 
the complex structure $J_u$ on $E_u$ and 
the Hodge $*$-operator on $\Om^1(\Sigma)$. 
The fibre 
$$
    \Ee_u 
    = \Om^{0,1}_{J_u}(\Sigma,E_u)
      \oplus\Om^0(\Sigma,\g_P)\oplus\Om^0(\Sigma,\g_P)
$$
also admits a complex structure given by the complex structure
$J_u$ on $E_u$ and by the map $(\phi,\psi)\mapsto(-\psi,\phi)$ on
$\Om^0(\Sigma,\g_P)\oplus\Om^0(\Sigma,\g_P)$. 
The only term in the formula~(\ref{eq:DuA}) for 
the operator $\Dd_{u,A}$ that is not necessarily 
complex linear, is the operator 
$$
     D_{u,A}:\Cinf(\Sigma,E_u)\to\Om^{0,1}_{J_u}(\Sigma,E_u).
$$
However, by Remark~\ref{rmk:Tabla}, the complex anti-linear
part of $D_{u,A}$ is of zeroth order and is therefore compact.
Hence $\Dd_{u,A}$ is a compact perturbation 
of a complex linear operator and hence admits
a natural orientation. 
It follows that the real line bundle
$\det(\Dd)\to\Bb$ admits a natural orientation
(see~\cite{MS1}). 

Now let $g\in\Gg(P)$ and choose a pair $(u,A)\in\Bb$.
Linearizing the action of the gauge group 
gives rise to isomorphisms
\begin{equation}\label{eq:iso}
\begin{array}{rcl}
     \ker\Dd_{u,A}\to\ker\Dd_{g^{-1}u,g^*A}
     &:&
     (\xi,\alpha)\mapsto(g^{-1}\xi,g^{-1}\alpha g), \\
     \ker\Dd^*_{u,A}\to\ker\Dd^*_{g^{-1}u,g^*A}
     &:&
     (\eta,\phi,\psi)\mapsto(g^{-1}\eta,g^{-1}\phi g,g^{-1}\psi g).
\end{array}
\end{equation}
We prove that the resulting isomorphism of determinant
lines is orientation preserving with respect
to the natural orientations introduced above.
To see this, we assume first that $H=0$ 
and $J(z,\cdot)$ is integrable near $u(p)$ 
for every $z\in\Sigma$ and every $p\in P$ 
with $\pi(p)=z$. Then the operators $\Dd_{u,A}$
and $\Dd_{g^{-1}u,g^*A}$ are both complex linear,
and hence the orientations of $\det(\Dd_{u,A})$
and $\det(\Dd_{g^{-1}u,g^*A})$ both agree
with the orientations induced by the complex
structures.  Hence, in this case, 
the result follows from the fact that the 
maps~(\ref{eq:iso}) are complex linear.  
In general, the result
follows from the fact that the spaces $\Hh$ and $\Jj$
are both connected and hence the isomorphism
$\det(\Dd_{u,A})\to\det(\Dd_{g^{-1}u,g^*A})$
is orientation preserving for some pair 
$(H,J)$ if and only if it is orientation
preserving for every pair $(H,J)$. 
\end{proof}

\begin{proposition}\label{prop:orient-cob}
Let 
$
     \{H_\lambda\}_\lambda\in
     \Hreg^{[0,1]}(H_0,H_1;\{\tau_\lambda,J_\lambda\}_\lambda).
$
Then the moduli space 
$
     \Ww_{B,\Sigma}^*(\{\tau_\lambda,J_\lambda,H_\lambda\}_\lambda
$
is an oriented cobordism from 
$\Mm_{B,\Sigma}^*(\tau_0;J_0,H_0)$ to
$\Mm_{B,\Sigma}^*(\tau_1;J_1,H_1)$. 
\end{proposition}

\begin{proof}
The tangent space of $\Ww_{B,\Sigma}^*$ at a triple
$(\lambda,u,A)$ is the kernel of the operator
$$
     \hat\Dd_{\lambda,u,A}:\R\times T_{(u,A)}\Bb\to\Ee_u
$$
given by 
$$
     \hat\Dd_{\lambda,u,A}(\hat\lambda,\xi,\alpha)
     := \Dd_{u,A}(\xi,\alpha) + \hat\lambda\zeta_{\lambda,u,A}
$$
for $\hat\lambda\in\R$ and $(\xi,\alpha)\in T_{(u,A)}\Bb$,
where $\zeta_{\lambda,u,A}\in\Ee_u$ is given 
by~(\ref{eq:cobordism}).  Since $\Dd_{\lambda,u,A}$ 
is surjective, the orientation of the kernel is determined
by the orientation of the determinant line.
Thus we must examine the determinant line bundle
$$
     \det(\hat\Dd)\to [0,1]\times\Bb
$$
whose fibre over a triple $(\lambda,u,A)\in[0,1]\times\Bb$
is $\det(\Dd_{\lambda,u,A})$. The homotopy
$t\mapsto t\zeta_{\lambda,u,A}$ yields a natural
isomorphism
$$
     \det(\hat\Dd)
     \cong{\rm pr}_1^*T[0,1]\otimes{\rm pr}_2^*\det(\Dd)
     \cong{\rm pr}_2^*\det(\Dd),
$$
where ${\rm pr_1}:[0,1]\times\Bb\to[0,1]$
and ${\rm pr_2}:[0,1]\times\Bb\to\Bb$
denote the obvious projections.  
This is because the tangent space $T_\lambda[0,1]$
is canonically isomorphic to $\R$ and, for $t=0$, we have
$$
     \ker\hat\Dd\cong \R\times\ker\,\Dd,\qquad
     \coker\hat\Dd\cong\coker\Dd.
$$
Hence $\det(\hat\Dd)$ inherits the orientation of
$\det(\Dd)$ and, since the orientation of $\det(\Dd)$ 
is invariant under the action of $\Gg$ 
so is the orientation of $\det(\hat\Dd)$.
It follows that the manifold $\Ww^*_{B,\Sigma}$ 
admits a natural orientation. 

Now choose a triple 
$
     (\lambda,[u,A])\in\Ww^*_{B,\Sigma}
$ 
such that $\Dd_{u,A}$ is onto. 
Then a positively oriented basis of
$T_{(\lambda,[u,A])}\Ww^*_{B,\Sigma}$ has the 
form 
$$
    (1,\xi_0,\alpha_0),(0,\xi_1,\alpha_1),\dots,
    (0,\xi_k,\alpha_k),
$$
where the vectors $(\xi_1,\alpha_1),\dots,(\xi_k,\alpha_k)$
form a positively oriented basis of the kernel of
$\Dd_{u,A}$. With the standard convention for orienting the boundary
(the outward unit normal vector comes first) the result 
follows.
\end{proof}


\section{Integer invariants}\label{sec:integer}

Let $(M,\om,\mu)$ be a symplectic manifold of (real) dimension $2n$
equipped with a Hamiltonian action by a compact Lie group $\G$ which 
is generated by a moment map $\mu:M\to\g$.
Suppose that the triple $(M,\om,\mu)$ satisfies~$(H1-3)$.
We shall define rational invariants of the triple $(M,\om,\mu)$ 
for central regular values of the moment map
under these hypotheses.  Conditions~$(H1-2)$ are needed 
to prove that the moduli spaces are compact.  
It should be possible to remove condition~$(H3)$,
however, the construction of the invariants without this
condition will probably require considerably more analysis than
has been carried out in the present paper.
This would include a full version of compactness for the solutions
of~(\ref{eq:vortex-ham}) without loss of energy and with 
preservation of the homotopy class in the limit, as well
as the construction of virtual moduli cycles analogous
to the definition of the Gromov--Witten invariants for 
general symplectic manifolds as in~\cite{FO,LRR,LT,RUAN}.
On the other hand, there are many interesting examples
that satisfy~$(H1-3)$, such as linear actions on 
complex vector spaces. (In this case~$(H1)$ implies~$(H2)$,
and~$(H3)$ is obvious.) In the present section we define 
integer invariants under the additional assumption
that $\G$ acts freely  on $\mu^{-1}(\tau)$ for some
central element $\tau\in Z(\g)$. 
This hypothesis will in Section~\ref{sec:rational}
be replaced by the assumption that $\mu$ has a 
central regular value. 


\subsection{Smooth moduli spaces}

Fix an equivariant homology class~$B\in H_2(M_\G;\Z)$ and 
a compact Riemann surface $(\Sigma,J_\Sigma,\dvol_\Sigma)$. 
Recall from Section~\ref{sec:trans} the definition of the 
moduli space 
$$
      \Mm(\tau;J,H)
      = \Mm_{B,\Sigma}(\tau;J,H)
      = \Tilde{\Mm}_{B,\Sigma}(\tau;J,H)/\Gg(P)
$$
of gauge equivalence classes of solutions of~(\ref{eq:vortex-ham})
for a (family of) almost complex structures
$J\in\Jj$ and a Hamiltonian perturbation $H\in\Hh$. 
Consider the set 
$$
     Z_0(\g):=\left\{\tau\in Z(\g)\,|\,\G\mbox{ acts freely on }
     \mu^{-1}(\tau)\right\}.
$$
By~$(H1)$, this set is open and we assume here that it is nonempty.
Choose a smooth function $\delta:Z_0(\g)\to(0,\infty)$ 
such that 
\begin{equation}\label{eq:delta}
     |\mu(x)-\tau|^2\le\delta(\tau)\quad\IMP\quad
     \G_x=\{1\},\;\;\im\,L_x\cap\im\,J(z,x)L_x=\{0\}
\end{equation}
for all $(z,x)\in \Sigma\times M$ and $\tau\in Z_0(\g)$. If 
\begin{equation}\label{eq:irr}
     \frac{\inner{[\om+\tau-\mu]}{B} + \|\Om_H\|}{\Vol(\Sigma)} 
     \le\delta(\tau)
\end{equation}
then the moduli space $\Mm_{B,\Sigma}(\tau;J,H)$ consists 
entirely of irreducible solutions (Lemma~\ref{le:irreducible})
and hence is a smooth manifold for a generic 
Hamiltonian function $H$ that satisfies~(\ref{eq:irr})
(Theorem~\ref{thm:trans}).  Moreover, $\Mm_{B,\Sigma}(\tau;J,H)$ 
is compact whenever $H$ has compact support and 
$J$ agrees with the almost complex
structure $J_0$ of hypothesis~$(H2)$ ouside of some compact 
subset of $M$ (Corollary~\ref{cor:compact}).  
Let us denote by $\Jj_0$ the space of 
almost complex structures that agree with $J_0$ outside 
of a compact set, by $\Hh(\tau;\delta)$ the space of compactly 
supported Hamiltonian perturbations that satisfy~(\ref{eq:irr}),
and, for $\tau\in Z_0(\g)$ and $J\in\Jj_0$, denote by 
$$
     \Hreg(\tau,J;\delta) := \Hreg(\tau,J)\cap\Hh(\tau;\delta)
$$
the subset of regular perturbations in the sense of 
Definition~\ref{def:regular}.  Let
$$
     \Bb^*\subset\Bb = \Cinf_\G(P,M)\times\Aa(P)
$$
denote the subset of irreducible pairs $(u,A)\in\Bb$
(see Definition~\ref{def:irreducible}).  By the local slice theorem 
(see Theorem~\ref{thm:slice}) the quotient space 
$\Bb^*/\Gg$ is an infinite dimensional 
Fr\^echet manifold (determined by $\mu$, $P$, and $\Sigma$).  
The next theorem summarizes our results about the moduli
spaces $\Mm_{B,\Sigma}(\tau;J,H)$.  

\begin{theorem}\label{thm:moduli}
Assume~$(H1-3)$ and let $\tau\in Z_0(\g)$. 
Then the following holds.
\begin{description}
\item[(i)] 
For every pair $(\tau,J)\in Z_0(\g)\times\Jj_0$
the set $\Hreg(\tau,J;\delta)$ is open and dense 
in $\Hh(\tau;\delta)$. 
\item[(ii)]
For every pair $(\tau,J)\in Z_0(\g)\times\Jj_0$ 
and every $H\in\Hreg(\tau,J;\delta)$ the moduli space 
$
     \Mm_{B,\Sigma}(\tau;J,H)
$ 
is a compact smooth naturally oriented submanifold of 
$\Bb^*/\Gg$ of dimension 
$$
     \dim\Mm_{B,\Sigma}(\tau;J,H)
     = 2m: = (n-\dim\G)\chi(\Sigma) + 2\inner{c_1^\G(TM)}{B}.
$$
\item[(iii)] 
For $\lambda=0,1$ let $(\tau_\lambda,J_\lambda)\in Z_0(\g)\times\Jj_0$
and $H_\lambda\in\Hreg(\tau_\lambda,J_\lambda;\delta)$.
Suppose that $\tau_0$ and $\tau_1$ belong  to the same 
component of $Z_0(\g)$. Then the moduli spaces 
$\Mm(\tau_0;J_0,H_0)$ and $\Mm(\tau_1;J_1,H_1)$ are 
oriented cobordant in $\Bb^*/\Gg$, i.e. there exists a compact 
oriented submanifold
$$
     \Ww\subset[0,1]\times\Bb^*/\Gg
$$
of dimension $2m+1$ such that 
$$
     \p\Ww 
     = \{1\}\times\Mm(\tau_1;J_1,H_1)
     - \{0\}\times\Mm(\tau_0;J_0,H_0).
$$
\end{description}
\end{theorem}

\begin{proof}
By Corollary~\ref{cor:compact}, the moduli space 
$\Mm(\tau;J,H)$ is a compact subset 
of $\Bb/\Gg$ and, by Lemma~\ref{le:irreducible}, 
it consists entirely of irreducible solutions 
of~(\ref{eq:vortex-ham}) for every $H\in\Hh(\tau;\delta)$.
Hence the set $\Hreg(\tau,J;\delta)$ is open in $\Hh(\tau;\delta)$.
By Theorem~\ref{thm:trans}, it is dense.  This proves~(i).

That $\Mm^*(\tau;J,H)$ is a smooth submanifold
of $\Bb^*/\Gg$ of dimension $2m$
for $H\in\Hreg(\tau,J;\delta)$ follows from the definitions, 
from Proposition~\ref{prop:fredholm},
and from the implicit function theorem.
That $\Mm(\tau;J,H)=\Mm^*(\tau;J,H)$ follows
from Lemma~\ref{le:irreducible},
and that $\Mm(\tau;J,H)$ is orientable
follows from Proposition~\ref{prop:orient}.
This proves~(ii).  Assertion~(iii) follows
from Theorem~\ref{thm:cobordism},
Corollary~\ref{cor:compact}, and
Proposition~\ref{prop:orient-cob}.
\end{proof}


\subsection{Definition of the invariants}\label{sec:invariant}

\subsubsection*{The evaluation map}

The group $\Gg\times\G$ acts freely 
on the product $\Bb^*\times P$ by 
$$
     (g,h)^*(u,A,p) := (g^{-1}u,g^*A,pg(p)^{-1}h)
$$
for $g\in\Gg$, $h\in\G$, $(u,A)\in\Bb^*$, and $p\in P$. 
Hence there is a principal $\G$-bundle 
$$
     \Pp:= (\Bb^*\times P)/\Gg \TO \Bb^*/\Gg\times\Sigma.
$$
The classifying map $\Bb^*/\Gg\times\Sigma\to\BG$ of this bundle 
lifts to a map $\theta:\Bb^*\times P\to\EG$ that 
is $\Gg$-invariant and $\G$-equivariant:
$$
     \theta(g^{-1}u,g^*A,pg(p)^{-1}h)=h^{-1}\theta(u,A,p).
$$
Likewise, the evaluation map 
$
     \Bb^*\times P\to M:(u,A,p)\mapsto u(p)
$
is $\Gg$-invariant and $\G$-equivariant.
These two maps together give rise to a map
$$
     \ev_\G:\Bb^*/\Gg\times\Sigma \TO M_\G := M\times_\G\EG
$$
given by
$$
     \ev_\G([u,A,p]) := [u(p),\theta(u,A,p)].
$$
The composition of $\ev_\G$ 
with the projection $\rho_M:M_\G\to\BG$
is the classifying map of $\Pp$:
$$
\xymatrix{
     \Bb^*/\Gg\times\Sigma\ar[r]^{\;\;\;\;\;\ev_\G} 
     \ar[rd] & M_\G\ar[d]^{\rho_M} \\
     & \BG  
}.
$$
If $M$ is contractible, 
then the projection $\rho_M:M_\G\to\BG$ is a homotopy
equivalence. 

\subsubsection*{The projection}

Now fix a point $p_0\in P$, denote by 
$
    \Gg_0:=\{g\in\Gg\;|\;g(p_0)=\id\}
$
the {\em based gauge group}, and consider the space
$$
     \Aa_\Gg := \Aa\times_\Gg\EG,
$$
where $\Gg$ acts by 
$
     g^*(A,e) := (g^*A,g(p_0)^{-1}e).
$
This space can be identified with $\Aa/\Gg_0\times_\G\EG$.
Since $\Gg$ acts freely on $\Bb^*$ there 
is a principal $\G$-bundle 
$$
     \Pp_0 := \Bb^*/\Gg_0 \TO \Bb^*/\Gg.
$$
The classifying map of this bundle 
lifts to a $\Gg$-equivariant map $\theta_0:\Bb^*\to\EG$,
which is equal to the restriction of $\theta$ to 
$\Bb^*\times\{p_0\}$. It satisfies
$$
     \theta_0(g^{-1}u,g^*A)=g(p_0)^{-1}\theta_0(u,A)
$$
and gives rise to a projection
$$
     \pi_\Aa:\Bb^*/\Gg\to\Aa_\Gg
$$
given by
$$
     \pi_\Aa([u,A]):=[A,\theta_0(u,A)].
$$ 
The composition of $\pi_\Aa$ with the projection
$\rho_\Aa:\Aa_\Gg\to\BG$ is the classifying map of $\Pp_0$:
$$
\xymatrix{
     \Bb^*/\Gg\ar[r]^{\pi_\Aa} 
     \ar[rd] & \Aa_\Gg\ar[d]^{\rho_\Aa}  \\
     & \BG 
}.
$$

\subsubsection*{The invariants}

We define invariants of the sixtuple
$(M,\om,\mu,\tau,B,\Sigma)$ with $\tau\in Z_0(\g)$
by integrating suitable cohomology classes over the moduli 
space $\Mm_{B,\Sigma}(\tau;J,H)$.
Such cohomology classes can be obtained by 
pulling back equivariant cohomology classes
on $M$ under the evaluation map 
$\ev_\G$ and equivariant cohomology classes 
on $\Aa/\Gg_0$ under the projection $\pi_\Aa$.
Let $\alpha\in H^*(\Aa_\Gg;\Z)$, 
$\beta_1,\dots,\beta_k\in H^*(M_\G;\Z)$, 
and $\gamma_1,\dots,\gamma_k\in H_*(\Sigma;\Z)$ 
such that 
$$
      \deg(\alpha) + \sum_{i=1}^k\deg(\beta_i) 
      - \sum_{i=1}^k\deg(\gamma_i)
      = 2m,
$$
and define 
\begin{eqnarray*}
&&
     \Phi_{B,\Sigma}^{M,\mu-\tau}
     (\alpha;\beta_1,\dots,\beta_k;\gamma_1,\dots,\gamma_k) \\
&&
     := \int_{\Mm_{B,\Sigma}(\tau;J,H)}
        \pi_\Aa^*\alpha\smile\ev_\G^*\beta_1/\gamma_1\smile\cdots
        \smile\ev_\G^*\beta_k/\gamma_k.
\end{eqnarray*}
Here the map 
$
     H^q(\Bb^*/\Gg\times\Sigma;\Z)\times H_i(\Sigma;\Z)
     \to H^{q-i}(\Bb^*/\Gg;\Z):
     (\beta,\gamma)\mapsto\beta/\gamma
$
denotes the slant product, $J\in\Jj_0$, and 
$H\in\Hreg(\tau,J;\delta)$, where the function
$\delta:Z_0(\g)\to(0,\infty)$ satisfies~(\ref{eq:delta}).

\begin{theorem}\label{thm:invariant}
Assume~$(H1-3)$ and let $\tau\in Z_0(\g)$ (i.e. $\G$ acts
freely on $\mu^{-1}(\tau)$). The invariant 
$
     \Phi_{B,\Sigma}^{M,\mu-\tau}
     (\alpha;\beta_i;\gamma_i) 
$
is independent of the almost complex structure
$J$ and the Hamiltonian perturbation $H$ used to
define it. It depends only on the triple 
$(M,\om,\mu)$, on the (co)homology classes
$B,\alpha,\beta_i,\gamma_i$, and on the 
component of $\tau$ in $Z_0(\g)$.
\end{theorem}

\begin{proof}
The space $\Bb^*/\Gg$ depends on $M$, the $\G$-action, 
$\Sigma$, and $P$. The invariant is defined by pairing an
integral cohomology class on $\Bb^*/\Gg$, 
determined by $\alpha,\beta_i,\gamma_i$, with the homology class 
$[\Mm_{B,\Sigma}(\tau;J,H)]\in H_*(\Bb^*/\Gg;\Z)$.
By Theorem~\ref{thm:moduli}, the latter is independent 
of $J$ and $H$ and depends only on the component 
of $\tau$ in $Z_0(\g)$. That it is also independent of the 
metric on $\Sigma$ follows by a similar cobordism argument. 
\end{proof}

\begin{remark}\label{rmk:Hsmall}\rm
The hypothesis that the Hamiltonian is small
(compared to the volume of $\Sigma$) is quite 
restrictive. If we allow for more general (abstract) 
perturbations of the symplectic vortex equations, then 
the condition that all solutions of~(\ref{eq:vortex-ham})
are irreducible can be replaced by the weaker condition 
that the gauge group acts freely on the space 
of solutions of~(\ref{eq:vortex-ham}).  
In the case of linear torus actions 
this condition is satisfied
for every Hamiltonian perturbation whenever
$\G$ acts freely on $\mu^{-1}(\tau_0)$, where
$\tau_0\in\g$ is defined by~(\ref{eq:tau0})
(see Remark~\ref{rmk:regular}).
\end{remark}


\section{$\G$-moduli problems and the Euler class}\label{sec:moduli}

In this section we review the results of~\cite{CMS}
about the Euler class of $\G$-moduli problems.   
They play a crucial role in the definition of the 
rational invariants in the presence of finite isotropy.

\begin{definition}\label{def:moduli}
Let $\G$ be a compact oriented Lie group.
A {\bf $\G$-moduli problem} is a triple 
$(\Bb,\Ee,\Ss)$ with the following properties.
\begin{description}
\item[$\bullet$]
$\Bb$ is a Hilbert manifold 
equipped with a smooth $\G$-action.
\item[$\bullet$]
$\Ee$ is a Hilbert space bundle over $\Bb$, also 
equipped with a smooth $\G$-action, such that 
$\G$ acts by isometries on the fibres of $\Ee$
and the projection $\Ee\to\Bb$ is $\G$-equivariant.
\item[$\bullet$]
$\Ss:\Bb\to\Ee$ is a smooth $\G$-equivariant Fredholm 
section of constant Fredholm index 
such that the determinant bundle $\det(\Ss)\to\Bb$ 
is oriented, $\G$ acts by orientation preserving isomorphisms
on the determinant bundle, and the zero set
$$
    \Mm := \left\{x\in\Bb\,|\,\Ss(x) = 0\right\}
$$
is compact.
\end{description}
A finite dimensional $\G$-moduli problem $(B,E,S)$ 
is called {\bf oriented} if $B$ and $E$ are oriented 
and $\G$ acts on $B$ and $E$ by orientation preserving 
diffeomorphisms.
A $\G$-moduli problem $(\Bb,\Ee,\Ss)$ is called {\bf regular}
if the isotropy subgroup 
$
    \G_x := \left\{g\in\G\,|\,g^*x=x\right\}
$
is finite for every $x\in\Mm$.
\end{definition}

$\G$-moduli problems form a category as follows.

\begin{definition}\label{def:morphism}
Let $(\Bb,\Ee,\Ss)$, $(\Bb',\Ee',\Ss')$ 
be $\G$-moduli problems.  
A {\bf morphism} from $(\Bb,\Ee,\Ss)$ to $(\Bb',\Ee',\Ss')$ 
is a pair $(\psi,\Psi)$ with the following properties.
$
     \psi:\Bb_0\to\Bb'
$
is a smooth $\G$-equivariant embedding
of a neighbourhood $\Bb_0\subset\Bb$ of $\Mm$ into $\Bb'$,
$
     \Psi:\Ee_0:=\Ee|_{\Bb_0}\to\Ee'
$
is a smooth injective bundle homomorphism and a lift of~$\psi$, 
and the sections $\Ss$ and $\Ss'$ satisfy
$$
     \Ss'\circ\psi=\Psi\circ\Ss,\qquad 
     \Mm'=\psi(\Mm).
$$
Moreover, the linear operators $d_x\psi:T_x\Bb\to T_{\psi(x)}\Bb'$ 
and $\Psi_x:\Ee_x\to\Ee'_{\psi(x)}$ induce isomorphisms
\begin{equation}\label{eq:ker-coker}
        d_x\psi:\ker\Dd_x\to\ker\Dd'_{\psi(x)},\qquad 
        \Psi_x:\coker\Dd_x\to\coker\Dd'_{\psi(x)},
\end{equation}
for $x\in\Mm$, and the resulting isomorphism
from $\det(\Dd)$ to $\det(\Dd')$ 
is orientation preserving. 
\end{definition}

Let $(\Bb,\Ee,\Ss)$ and $(\Bb',\Ee',\Ss')$ 
be $\G$-moduli problems and suppose that there exists a 
morphism from $(\Bb,\Ee,\Ss)$ to $(\Bb',\Ee',\Ss')$.
Then the indices of $\Ss$ and $\Ss'$ agree.
Moreover, $(\Bb,\Ee,\Ss)$ is regular if and only if 
$(\Bb',\Ee',\Ss')$ is regular. 

\begin{definition}\label{def:G-homotopy}
Two regular $\G$-moduli problems $(\Bb,\Ee_i,\Ss_i)$, $i=0,1$,
(over the same base) are called {\bf homotopic} if there exists a
$\G$-equivariant Hilbert space bundle 
$
     \Ee\to[0,1]\times\Bb
$
and a $\G$-equivariant smooth section $\Ss:[0,1]\times\Bb\to\Ee$
such that $\Ee_i=\Ee|_{\{i\}\times\Bb}$ and 
$\Ss_i=\Ss|_{\{i\}\times\Bb}$ for $i=0,1$,
the triple $(\Bb,\Ee_t,\Ss_t)$, defined by 
$\Ee_t:=\Ee|_{\{t\}\times\Bb}$ and 
$\Ss_t=\Ss|_{\{t\}\times\Bb}$, 
is a regular $\G$-moduli problem
for every $t\in[0,1]$, and the set 
$
     \Mm:=\left\{(t,x)\in[0,1]\times\Bb\,|\,\Ss_t(x)=0\right\}
$
is compact.
\end{definition}

The following theorem in proved in~\cite{CMS}. 
It states the properties of the Euler class. 
We denote by $H^*_\G(\Bb)$ the equivariant cohomology
with real coefficients.

\begin{theorem}\label{thm:euler}
There exists a functor, called the 
{\bf Euler class}, which assigns to each
compact oriented Lie group $\G$ and each regular
$\G$-moduli problem $(\Bb,\Ee,\Ss)$ a homomorphism
$
      \chi^{\Bb,\Ee,\Ss}:H^*_\G(\Bb)\to\R
$
and satisfies the following.
\begin{description}
\item[(Functoriality)]
If $(\psi,\Psi)$ is a morphism from 
$(\Bb,\Ee,\Ss)$ to $(\Bb',\Ee',\Ss')$ then
$
      \chi^{\Bb,\Ee,\Ss}(\psi^*\alpha)
      = \chi^{\Bb',\Ee',\Ss'}(\alpha)
$
for every $\alpha\in H^*_\G(\Bb')$.
\item[(Thom class)]
If $(B,E,S)$ is a finite dimensional oriented regular $\G$-moduli 
problem and $\tau\in\Om^*_\G(E)$ is an equivariant Thom form 
supported in an open neighbourhood $U\subset E$ of the zero 
section such that $U\cap E_x$ is convex for every $x\in B$,
$U\cap\pi^{-1}(K)$ has compact closure for every
compact set $K\subset B$, and $S^{-1}(U)$ has compact closure, 
then
$$
      \chi^{B,E,S}(\alpha) = \int_{B/\G}\alpha\wedge S^*\tau
$$
for every closed form $\alpha\in\Om^*_\G(B)$.
\item[(Transversality)]
If $\Ss$ is transverse to the zero section then
$$
      \chi^{\Bb,\Ee,\Ss}(\alpha) = \int_{\Mm/\G}\alpha
$$
for every $\alpha\in H^*_\G(B)$, where $\Mm:=\Ss^{-1}(0)$.
\item[(Homotopy)]
If $(\Bb,\Ee_0,\Ss_0)$ and $(\Bb,\Ee_1,\Ss_1)$ are 
regular homotopic $\G$-moduli problems then
$
      \chi^{\Bb,\Ee_0,\Ss_0}(\alpha)
      = \chi^{\Bb,\Ee_1,\Ss_1}(\alpha)
$
for every $\alpha\in H^*_\G(\Bb)$.
\item[(Subgroup)]
If $(\Bb,\Ee,\Ss)$ is a regular $\G$-moduli problem and $\HG\subset\G$ 
is a normal subgroup acting freely on $\Bb$ then
$$
      \chi^{\Bb/\HG,\Ee/\HG,\Ss/\HG}(\alpha)
      = \chi^{\Bb,\Ee,\Ss}(\pi^*\alpha)
$$
for every $\alpha\in H^*_{\G/\HG}(\Bb/\HG)$, 
where $\pi^*:H^*_{\G/\HG}(\Bb/\HG)\to H^*_\G(\Bb)$
is the homomorphism induced by the projection
$\pi:\Bb\to\Bb/\HG$.
\item[(Rationality)]
If $\alpha\in H^*_\G(\Bb;\Q)$ then 
$\chi^{\Bb,\Ee,\Ss}(\alpha)\in\Q$. 
\end{description}
The Euler class is uniquely determined by the 
{\it (Functoriality)} and {\it (Thom class)} axioms.
\end{theorem}


\section{Rational invariants}\label{sec:rational}

Our next goal is to drop the hypothesis that $\G$ acts freely 
on $\mu^{-1}(\tau)$ and construct invariants 
for every central regular values $\tau\in Z(\g)$ 
of the moment map $\mu:M\to\g$.  
In this case we must deal with the 
presence of finite isotropy subgroups.  
We assume as before that the hypotheses~$(H1-3)$ 
are satisfied. 


\subsection{The setup}\label{sec:setup}

Let us fix the following data:
\begin{description}
\item[$\bullet$]
an equivariant homology class $B\in H_2(M_\G;\Z)$,
\item[$\bullet$]
a compact connected Riemann surface $(\Sigma,J_\Sigma,\dvol_\Sigma)$
and a principal $\G$-bundle $\pi:P\to\Sigma$ whose characteristic
class is the image of $B$ under the homomorphism
$H_2(M_\G;\Z)\to H_2(\BG;\Z)$, 
\item[$\bullet$]
a point $p_0\in P$ and an integer $k\ge 2$.
\end{description}
We emphasize that the purpose of fixing the point
$p_0$ is not in the definition of the evaluation 
map, but to obtain an action of the gauge group
on the classifying space $\EG$ of $\G$. 
Throughout we shall denote by 
$\Gg:=\Gg^{k+1,2}(P)$ the group of 
gauge transformations of $P$ of class $W^{k+1,2}$ 
and by 
$$
     \Gg_0 := \Gg^{k+1,2}_0(P)
     := \left\{g\in\Gg^{k+1,2}(P)\,|\,g(p_0)=\one\right\}
$$
the (normal) subgroup of based gauge transformation. 
When the need arises we shall think of the gauge group 
$\Gg$ as acting on $\EG$ by 
$
     g^*e:=g(p_0)^{-1}e.
$
So the subgroup $\Gg_0$ acts trivially on $\EG$.

The above data give rise to 
a $\G$-moduli problem $(\Bb,\Ee,\Ss)$ as follows. 
The Hilbert manifold $\Bb$ is the quotient
\begin{equation}\label{eq:B}
     \Bb := \frac{W^{k,2}_\G(P,M;B)\times\Aa^{k,2}(P)}
             {\Gg_0^{k+1,2}(P)},
\end{equation}
where 
$
     W^{k,2}_\G(P,M;B)
     :=\{(u\in W^{k,2}_\G(P,M)\,|\,[u]=B\}.
$
The Hilbert manifold $\Bb$ carries a $\G$-action
since the quotient group $\Gg/\Gg_0$ is isomorphic to $\G$. 
Consider the bundle $\Ee\to\Bb$ with fibres
\begin{equation}\label{eq:E}
     \Ee_{(u,A)}
     := W^{k-1,2}(\Sigma,\Lambda^{0,1}_{J_\Sigma}T^*\Sigma
         \otimes_J u^*TM/\G)
        \oplus W^{k-1,2}(\Sigma,\g_P).
\end{equation}
The action of the gauge group identifies $\Ee_{(u,A)}$
with $\Ee_{(g^{-1}u,g^*A)}$ for every $g\in\Gg$.
Thus $\Ee$ carries a $\G$-action.  
More precisely, the fibre of $\Ee$ over a point 
$[u,A]\in\Bb$ is the union of the spaces $\Ee_{(g^{-1}u,g^*A)}$
over all $g\in\Gg_0$ and any two such spaces are
identified by the action of the based gauge group. 
Then the quotient group $\G\cong\Gg/\Gg_0$ acts on both 
$\Ee$ and $\Bb$ and the projection is $\G$-equivariant.
For every Hamiltonian perturbation $H$
the left hand side of equation~(\ref{eq:vortex-ham}) 
defines a $\G$-equivariant section $\Ss:\Bb\to\Ee$ given by 
\begin{equation}\label{eq:S}
     \Ss([u,A]) := [\bar\p_{J,H,A}(u),*F_A+\mu(u)-\tau].
\end{equation}

\begin{lemma}\label{le:regular}
Assume~$(H1-3)$ and let $\tau\in Z(\g)$ be a regular
value of $\mu$. 

\smallskip
\noindent{\bf (i)}
The triple $(\Bb,\Ee,\Ss)$ defined by~(\ref{eq:B}),
(\ref{eq:E}), and (\ref{eq:S}) is a $\G$-moduli problem
(see Section~\ref{sec:moduli}) of index
$$
     \INDEX(\Ss) = (n-\dim\G)\chi(\Sigma)
     + 2\inner{c_1^\G(TM)}{B} = :2m.
$$

\smallskip
\noindent{\bf (ii)}
There exists a constant $\delta>0$
such that the $\G$-moduli problem $(\Bb,\Ee,\Ss)$
is regular whenever
\begin{equation}\label{eq:Hsmall}
      \inner{[\om-\mu+\tau]}{B} + \left\|\Om_H\right\|
      \le \delta\Vol(\Sigma).
\end{equation}
\end{lemma}

\begin{proof}
That $\Ee$ is a Hilbert space bundle over a Hilbert manifold
is a consequence of the local slice theorem
(for $W^{k,2}$ connections). That $\Ss$ is a Fredholm section 
follows from Proposition~\ref{prop:fredholm} and so does the 
index formula. That the zero set of $\Ss$ is compact
follows from Corollary~\ref{cor:compact}.
This proves~(i).  Assertion~(ii) follows from 
Lemma~\ref{le:irreducible}. 
\end{proof}

We can now evaluate the Euler class $\chi^{\Bb,\Ee,\Ss}$,
defined in Section~\ref{sec:moduli}, on equivariant cohomology
classes of $\Bb$. As in Section~\ref{sec:invariant},
such cohomology classes can be obtained by pulling 
back equivariant cohomology classes on $M$ with the 
evaluation map and equivariant cohomology classes
on $\Aa/\Gg_0$ with the projection onto the 
space of connections. More precisely, 
abbreviate $\Aa:=\Aa^{k,2}(P)$ and $\Gg=\Gg^{k+1,2}(P)$. 
Consider the action of the group $\Gg\times\G$ on the 
space 
$$
      \Xx := \Aa\times P\times \EG
$$
by 
$$
      (g,h)^*(A,p,e) :=
      (g^*A,pg(p)^{-1}h,g(p_0)^{-1}e).
$$

\begin{lemma}\label{le:theta}
There exists a continuous function $\theta:\Xx\to\EG$
which is $\Gg$-invariant and $\G$-equivariant.
Thus
\begin{equation}\label{eq:theta}
     \theta(g^*A,pg(p)^{-1}h,g(p_0)^{-1}e)
     = h^{-1}\theta(A,p,e)
\end{equation}
for $(A,p,e)\in\Xx$, $g\in\Gg$, and $h\in\G$.
Any two such maps $\theta_0,\theta_1:\Xx\to\EG$ 
are homotopic through maps satisfying~(\ref{eq:theta}). 
Moreover, $\theta$ can be chosen such that
$$
     \theta(A,p_0,e) = e.
$$
\end{lemma}

\begin{proof}
The group $\Gg\times\G$ acts freely on $\Xx$.  
Hence the quotient $\Xx/\Gg$ is a principal 
$\G$-bundle over $\Xx/(\Gg\times\G)$.
The classifying map of this bundle lifts to a $\G$-equivariant map
from $\Xx/\Gg$ to $\EG$.  The composition of this map 
with the projection $\Xx\to\Xx/\Gg$ is the required
map $\theta$. The last two assertions follow from the fact
that any two classifying maps are equivariantly homotopic.
\end{proof}

Note that we can identify the space $\Bb\times_\G\EG$ with 
the quotient of the space 
$
     W^{k,2}_\G(P,M;B)\times\Aa^{k,2}(P)\times\EG
$ 
by $\Gg=\Gg^{k+1,2}(P)$, where the action of the gauge group
is given by 
$
     g^*(u,A,e) := (g^{-1}u,g^*A,g(p_0)^{-1}e).
$
Hence there is an evaluation map
$$
     \ev_\G:(\Bb\times_\G\EG)\times\Sigma\to M\times_\G\EG,
$$
defined by 
$$
     \ev_\G([u,A,e],\pi(p)) := [u(p),\theta(A,p,e)],
$$
and a projection 
$$
     \pi_\Aa:\Bb\times_\G\EG
     \to\Aa_\Gg:=\Aa\times_{\Gg}\EG
$$
defined by 
$$
     \pi_\Aa([u,A,e]):= [A,e].
$$


\subsection{Definition of the invariants}\label{sec:inv-rat}

Let 
$$
     \alpha\in H^*(\Aa_\Gg),\qquad
     \beta_1,\dots,\beta_k\in H^*(M_\G),\qquad
     \gamma_1,\dots,\gamma_k\in H_*(\Sigma)
$$ 
such that 
$$
      \deg(\alpha) + \sum_{i=1}^k\deg(\beta_i) 
      - \sum_{i=1}^k\deg(\gamma_i)
      = 2m,
$$
and define 
\begin{eqnarray*}
&&
     \Phi_{B,\Sigma}^{M,\mu-\tau}
     (\alpha;\beta_1,\dots,\beta_k;\gamma_1,\dots,\gamma_k) \\ 
&&
     := \chi^{\Bb,\Ee,\Ss}\left(
        \pi_\Aa^*\alpha\smile\ev_\G^*\beta_1/\gamma_1\smile\cdots
        \smile\ev_\G^*\beta_k/\gamma_k\right).
\end{eqnarray*}
Here the map
$
     H^q((\Bb\times_\G\EG)\times\Sigma)\times H_i(\Sigma)
     \to H^{q-i}(\Bb\times_\G\EG):
     (\beta,\gamma)\mapsto\beta/\gamma
$
denotes the slant product, the $\G$-moduli problem
$(\Bb,\Ee,\Ss)$ is defined by~(\ref{eq:B}),
(\ref{eq:E}), and~(\ref{eq:S}), where 
the Hamiltonian perturbation $H$ satisfies~(\ref{eq:Hsmall}),
and the Euler class $\chi^{\Bb,\Ee,\Ss}:H^*_\G(\Bb)\to\R$ 
is defined in Section~\ref{sec:moduli}. 

\begin{theorem}\label{thm:invariant-rat}
Assume~$(H1-3)$ and let $\tau\in Z(\g)$ be a regular
value of $\mu$. The invariant 
$
     \Phi_{B,\Sigma}^{M,\mu-\tau}
     (\alpha;\beta_i;\gamma_i) 
$
is independent of the almost complex structure
$J$, the Hamiltonian perturbation $H$, 
the point $p_0\in P$, and the integer $k$
used to define it. It depends only on 
$(M,\om,\mu)$, on the genus of $\Sigma$,
on the component of $\tau$ in the (open) set of central
regular values of $\mu$, and on the (co)homology classes 
$B,\alpha,\beta_i,\gamma_i$.
\end{theorem}

\begin{proof}
The independence of $k$ follows from the fact that 
a finite dimensional reduction for $k=2$ is also a finite 
dimensional reduction for every $k>2$, and that the classifying 
map $\theta^{k,2}:\Xx^{k,2}\to\EG$ can be defined as the restriction of
the classifying map $\theta^{2,2}$ to the subspace 
$\Xx^{k,2}\subset\Xx^{2,2}$. 
The independence of $J$, $H$, $J_\Sigma$, 
$\dvol_\Sigma$, and $\tau$ follows from the 
{\it (Homotopy)} axiom for the Euler class. 

We prove the independence of the basepoint $p_0$.
Let $p_1\in P$ and suppose that $H=0$
and that $J$ is independent of the point $z\in\Sigma$.
Choose a diffeomorphism $\phi:\Sigma\to\Sigma$
that is isotopic to the identity and a $\G$-equivariant lift 
$\psi:P\to P$ such that $\psi(p_1)=p_0$. 
Then the $\G$-moduli problem 
with $p_0$, $J_\Sigma$, and $\dvol_\Sigma$
replaced by $p_1=\psi^{-1}(p_0)$, $\phi^*J_\Sigma$,
and $\phi^*\dvol_\Sigma$, respectively, is diffeomorphic
to the original one.  The diffeomorphism
is given by $[u,A]\mapsto[u\circ\psi,\psi^*A]$. 
Hence the invariants are the same. 
\end{proof}

\begin{remark}\label{rmk:Hsmall-rat}\rm
We emphasize again that the condition~(\ref{eq:Hsmall}) 
on the Hamiltonian perturbation is quite restrictive 
and that much more general regularity criteria 
are available in the abelian case.
For example, if $\G$ is abelian and acts linearly 
on $M=\C^n$ and the element $\tau_0\in\g$, defined 
by~(\ref{eq:tau0}), is a regular value of $\mu$ 
then, for every Hamiltonian perturbation and every
almost complex structure, the gauge group acts on the 
space of solutions of~(\ref{eq:vortex-ham})
with finite isotropy (see Remarks~\ref{rmk:regular}).  
So in this case the smallness condition~(\ref{eq:Hsmall}) 
on the Hamiltonian can be dropped. Such more general 
criteria can also be obtained in the nonabelian case. 
\end{remark}

\begin{theorem}\label{thm:equal}
Assume~$(H1-3)$ and let $\tau\in Z(\g)$ such that 
$\G$ acts freely on $\mu^{-1}(\tau)$.  Then the invariant 
$\Phi_{B,\Sigma}^{M,\mu-\tau}(\alpha;\beta_i;\gamma_i)$
defined in this section agrees with the one defined
in Section~\ref{sec:invariant}.
\end{theorem}

\begin{proof}
By Theorem~\ref{thm:moduli}, there exists a Hamiltonian 
perturbation $H$ that satisfies~(\ref{eq:Hsmall})
and is regular in the sense of Definition~\ref{def:regular}.
For such a perturbation the section $\Ss:\Bb\to\Ee$, 
defined by~(\ref{eq:S}), is transverse
to the zero section. Hence the result follows from the 
{\it (Transversality)} axiom for the Euler class.
\end{proof}


\subsection{Relations}\label{sec:relations}

There are two kinds of relations between the invariants 
defined in Section~\ref{sec:inv-rat}, namely those arising
from relations between the slant product and the cup product
and others arising from relations between certain universal
bundles in gauge theory. 

\begin{proposition}\label{prop:slant}
Let $\i:\Sigma\to\Sigma\times\Sigma$ denote
the inclusion of the diagonal and suppose that
$\beta,\beta',\beta''\in H^*(M_\G;\Z)$ and 
$\gamma,\gamma_i',\gamma_i''\in H_*(\Sigma;\Z)$
satisfy
$$
     \beta = \beta'\smile\beta'',\qquad
     \i_*\gamma = \sum_{i=1}^m\gamma_i'\otimes\gamma_i''.
$$
Then 
\begin{eqnarray*}
&&
     \Phi^{M,\mu-\tau}_{B,\Sigma}(\alpha;\beta,\beta_1,\dots,\beta_k;
     \gamma,\gamma_1,\dots,\gamma_k) \\
&&=
     \sum_{i=1}^m
     (-1)^{\deg(\gamma_i')\deg(\beta''/\gamma_i'')}
     \Phi^{M,\mu-\tau}_{B,\Sigma}(\alpha;\beta',\beta'',\beta_1,\dots,\beta_k;
     \gamma_i',\gamma_i'',\gamma_1,\dots,\gamma_k).
\end{eqnarray*}
In particular, if $\gamma_i=[\point]$ for every $i$, then 
$$
     \Phi_{B,\Sigma}^{M,\mu-\tau}
     (\alpha;\beta_1,\dots,\beta_k;\point,\dots,\point) 
     = \Phi_{B,\Sigma}^{M,\mu-\tau}
     (\alpha;\beta_1\smile\dots\smile\beta_k;\point).
$$
\end{proposition}

\begin{proof}
This follows from the formula
$$
     \int_\gamma\alpha'\smile\alpha''
     = \sum_{i=1}^m\int_{\gamma_i'}\alpha'\int_{\gamma_i''}\alpha''
$$
for $\alpha',\alpha''\in H^*(\Sigma;\Z)$.
\end{proof}

Let us denote by 
$$
     \theta_\Aa:\Aa_\Gg\times\Sigma\to\BG
$$
the map 
$$
     \theta_\Aa([A,e],\pi(p)) := [\theta(A,p,e)],
$$
where $\theta$ is as in Lemma~\ref{le:theta}.
This is a classifying map for the bundle 
$$
     \Pp_\Aa := (\Aa\times P\times\EG)/\Gg
     \to \Aa_\Gg\times\Sigma,
$$
where the gauge group $\Gg$ acts by
$
     g^*(A,p,e) := (g^*A,pg(p)^{-1},g(p_0)^{-1}e).
$
Recall that $\rho_M:M_\G\to\BG$ denotes the projection.

\begin{proposition}\label{prop:A}
For every $c\in H^*(\BG;\Z)$ and every $\gamma\in H_*(\Sigma;\Z)$, 
\begin{eqnarray*}
&&
     \Phi^{M,\mu-\tau}_{B,\Sigma}(\alpha;
     \rho_M^*c,\beta_1,\dots,\beta_k;
     \gamma,\gamma_1,\dots,\gamma_k) \\
&&=
     \Phi^{M,\mu-\tau}_{B,\Sigma}(\alpha\smile(\theta_\Aa^*c/\gamma);
     \beta_1,\dots,\beta_k;
     \gamma_1,\dots,\gamma_k).
\end{eqnarray*}
\end{proposition}

\begin{proof}
By definition of the maps, there is a commuting diagram
$$
\xymatrix{
    (\Bb\times_\G\EG)\times\Sigma\ar[d]_{\pi_\Aa\times\id}
    \ar[r]^{\;\;\;\;\;\ev_\G} &  
    M\times_\G\EG \ar[d]^{\rho_M} \\
    \Aa_\Gg\times\Sigma\ar[r]_{\theta_\Aa}& \BG}.
$$
Hence, for every class $c\in H^*(\BG;\Z)$ and every
$\gamma\in H_*(\Sigma;\Z)$, we have 
$$
     \ev_G^*\rho_M^*c/\gamma 
     = ((\pi_\Aa\times\id)^*\theta_\Aa^*c)/\gamma
     = \pi_\Aa^*(\theta_\Aa^*c/\gamma).
$$
This proves the proposition.
\end{proof}

Let us now consider the abelian case $\G=T$.
Then the constant gauge transformations act trivially
on $\Aa$. Hence there is a principal bundle 
$$
     \Pp_{\Aa/\Gg} := (\Aa\times P)/\Gg_0
     \to \Aa/\Gg\times\Sigma,
$$
where $\Gg_0$ is the based gauge group. Let us denote by 
$$
     \theta_{\Aa/\Gg}:\Aa/\Gg\times\Sigma\to\BG.
$$
It lifts to a map $\theta_0:\Aa\times P\to\EG$
that satisfies
$$
     \theta_0(g^*A,pg(p)^{-1}) = g(p_0)^{-1}\theta_0(A,p).
$$ 
Consider the homomorphism
$
     \mu_{\Aa/\Gg}:H^q(\BG;\Z)\times H_i(\Sigma;\Z)\to H^{q-i}(\Aa/\Gg;\Z)
$
defined by 
$$
     \mu_{\Aa/\Gg}(c,\gamma) := \theta_{\Aa/\Gg}^*c/\gamma
     \in H^*(\Aa/\Gg;\Z).
$$
Let $\pi_{\Aa/\Gg}:\Aa_\Gg\to\Aa/\Gg$ denote
the obvious projection.

\begin{proposition}\label{prop:A0}
Assume the abelian case $\G=T$. 
Then, for every $c\in H^*(\BG;\Z)$ and every $\gamma\in H_*(\Sigma;\Z)$
such that 
$$
     \deg(\gamma)>0,
$$
we have
\begin{eqnarray*}
&&
     \Phi^{M,\mu-\tau}_{B,\Sigma}(\alpha;
     \rho_M^*c,\beta_1,\dots,\beta_k;
     \gamma,\gamma_1,\dots,\gamma_k) \\
&&=
     \Phi^{M,\mu-\tau}_{B,\Sigma}
     (\alpha\smile\pi_{\Aa/\Gg}^*\mu_{\Aa/\Gg}(c,\gamma);
     \beta_1,\dots,\beta_k;
     \gamma_1,\dots,\gamma_k).
\end{eqnarray*}
\end{proposition}

\begin{proof}
In the abelian case the projection
$$
     \Aa_\Gg \to \Aa/\Gg\times\BG
$$
is a homeomorphism. 
We define the {\em tensor product} of two 
principal $T$-bundles $\pi_1:P_1\to X$ and $\pi_2:P_2\to X$
as the quotient
$$
     P_1\otimes P_2 
     := \left\{(p_1,p_2)\in P_1\times P_2\,|\,\pi_1(p_1)=\pi_2(p_2)
     \right\}/T
$$
by the diagonal action. With this notation 
$$
     \Pp_\Aa \cong 
     (\Pp_{A/\Gg}\times\BG)\otimes(\Aa/\Gg\times\Sigma\times\EG)
     \TO \Aa/\Gg\times\Sigma\times\BG.
$$
An explicit bundle isomorphism is 
$
   [A,p,e]\mapsto[([A,p]_0,[e]),([A],[p],e)],
$
where $[A,p]_0\in\Pp_{\Aa/\Gg}=(\Aa\times P)/\Gg_0$ 
denotes the equivalence class under the action of $\Gg_0$,
and $[A]$, respectively $[p]$ and $[e]$, 
denote the equivalence classes
under the action of $\Gg$, respectively $\G$. 
Hence
$$
     \deg(\gamma)>0\qquad\IMP\qquad
     \theta_\Aa^*c/\gamma
     = \pi_{\Aa/\Gg}^*(\theta_{\Aa/\Gg}^*c/\gamma)
     = \pi_{\Aa/\Gg}^*\mu_{\Aa/\Gg}(c,\gamma)
$$
and hence the result follows from 
Proposition~\ref{prop:A}.
\end{proof}

The classes $\mu_{A/\Gg}(c,\gamma)$ are easy to compute and 
they generate the cohomology of $\Aa/\Gg=\Aa/\Gg_0$.  
Let $\Lambda:=\exp^{-1}(\one)\subset\tt$ and denote by 
$\W\subset\tt$ the dual lattice (of elements
whose periods on $\Lambda$ are integer multiples
of $2\pi$). Then every $\w\in\W$ determines a homomorphism 
$$
      \rho_\w:T\to S^1, 
$$
given by $\rho_\w(\exp(\tau)):=\exp(i\inner{\w}{\tau})$ 
and hence complex line bundles 
$$
     \Ll^\w := \EG\times_{\rho_\w}\C \to \BG,\qquad
     L^\w:=P\times_{\rho_\w}\C\to\Sigma.
$$
The first Chern class of $\Ll^\w$ will be denoted by
$$
     c_\w := c_1(\Ll^\w) \in H^2(\BG;\Z).
$$
We describe the map $\gamma\mapsto\mu_{\Aa/\Gg}(c_\w,\gamma)$ 
explicitly. Every $\w\in\W$ and every loop 
$\gamma:S^1\to\Sigma$ determine a real valued 1-form 
on $\Aa$ given by
\begin{equation}\label{eq:1form}
     T_A\Aa=\Omega^1(\Sigma,\tt)\to\R:
     \alpha\mapsto-\frac{1}{2\pi}\int_\gamma\inner{\w}{\alpha}.
\end{equation}
This 1-form is closed and $\Gg$-invariant,
so it descends to a closed 1-form 
$
     \tilde\mu_\w(\gamma)\in\Omega^1(\Aa/\Gg).
$
Similarly, every 2-chain $\sigma$ on $\Sigma$ induces a
function $\tilde\mu(\sigma)$ on $\Aa/\Gg$ defined by
$$
     \tilde\mu_\w(\sigma)([A]) 
     := -\frac{1}{2\pi}\int_\sigma \inner{\w}{F_A}.
$$
By Stokes' Theorem, $d\tilde\mu_\w(\sigma)=\tilde\mu_\w(\p\sigma)$
and hence there are induced homomorphisms 
$\tilde\mu_\w:H_1(\Sigma;\Z)\to H^1(\Aa/\Gg;\Z)$ and
$\tilde\mu_\w:H_2(\Sigma;\Z)\to H^0(\Aa/\Gg;\Z)$. 
The following lemma asserts that 
$\tilde\mu_\w(\gamma)=\mu_{\Aa/\Gg}(c_\w,\gamma)$. 

\begin{lemma}\label{L:configA}
For every $\w\in\W$ the following holds. 
\begin{description}
\item[(i)] $\mu_{\Aa/\Gg}(c_\w,[pt])=0$ and
$\mu_{\Aa/\Gg}(c_\w,[\Sigma])=\inner{c_1(L^\w)}{[\Sigma]}
\in H^0(\Aa/\Gg;\Z)$. 
\item[(ii)] For $\gamma\in H_1(\Sigma;\Z)$ the class
$\mu_{\Aa/\Gg}(c_\w,\gamma)\in H^1(\Aa/\Gg;\Z)$ 
is represented by the closed 1-form on $\Aa/\Gg$ 
induced by~(\ref{eq:1form}).
The Poincar\'e dual of $\mu_{\Aa/\Gg}(c_\w,\gamma)$ 
is represented by the cycle 
$$
     C_{\Aa/\Gg}(\w,\gamma)
     := \left\{[A]\in\Aa/\Gg\;\Big|\;
        \int_\gamma\inner{\w}{A-A_0}=0\right\}
$$
with the orientation determined by $-\w$.  
\item[(iii)] 
The map 
$$
     \W\otimes H_1(\Sigma;\Z)\to H^1(\Aa/\Gg;\Z):
     (\w,\gamma)\mapsto\mu_{\Aa/\Gg}(c_\w,\gamma)
$$
induces an isomorphism from the exterior algebra on 
the free $\Z$-module $\W\otimes H_1(\Sigma;\Z)$
to $H^*(\Aa/\Gg;\Z)$.
\end{description}
\end{lemma}

\begin{proof}
The bundle 
$
    \Pp_{\Aa/\Gg}=(\Aa\times P)/\Gg_0\to\Aa/\Gg\times\Sigma
$
carries a universal connection
induced by the $\Gg_0$-invariant 1-form 
$\A\in\Omega^1(\Aa\times P,\tt)$, 
$$
    \A_{(A,p)}(\alpha,v) = A_p(v) + (d^*d)^{-1}d^*\alpha(p).
$$
Here $d^*d$ denotes the isomorphism
$$
    d^*d:\Om^0(\Sigma,\tt)\to
    \im(d^*:\Omega^1(\Sigma,\tt)\to\Om^0(\Sigma,\tt)).
$$
The curvature $F_\A$ of $\A$ is given by
$$
    (F_\A)_{([A],z)}\Bigl((\alpha,v),(\beta,w)\Bigr) 
    = F_A(v,w)+\alpha_z(w)-\beta_z(v)
$$ 
for $v,w\in T_z\Sigma$ and $\alpha,\beta\in\Omega^1(\Sigma,i\R)$. 
Let $[\gamma]\in H_1(\Sigma;\Z)$ and $[A]\in H_1(\Aa/\Gg;\Z)$ be
represented by loops $\gamma:\R/\Z\to\Sigma$ and $A:\R\to\Aa$ such that
$A(t+1)=g^*A(t)$ for some $g\in\Gg$. Since
the closed 2-form $-\inner{\w}{F_\A}/2\pi$ represents the cohomology
class $c_1(\theta_{\Aa/\Gg}^*\Ll^\w)=\theta_{\Aa/\Gg}^*c_\w$, we have
\begin{align*}
  \inner{\mu_{\Aa/\Gg}(c_\w,[\gamma])}{[A]} 
  &= - \frac{1}{2\pi}\int_0^1\int_0^1\inner{\w}
       {F_\A(\dot A(s),\dot\gamma(t))}\,ds\,dt \\ 
  &= - \int_0^1\Bigl(\frac{1}{2\pi}\int_\gamma\inner{\w}{\dot A(s)}
     \Bigr)\,ds \\
  &= \langle\tilde\mu_\w([\gamma]),[A]\rangle.
\end{align*}
This proves~(ii). 
Assertion~(i) is proved similarly. 
To prove~(iii) let $m:=\dim T$ and pick integer bases
$\w_1,\dots,\w_m$ of $\W$ and
$\gamma_1,\dots,\gamma_{2g}$ of $H_1(\Sigma;\Z)$.
Now let $\alpha_{ij}\in\Om^1(\Sigma,\tt)$ be a harmonic
$1$-form such that
$$
     -\frac{1}{2\pi}\int_{\gamma_i}\inner{\w_j}{\alpha_{i'j'}} 
     = \delta_{ii'}\delta_{jj'}
$$
for $i,i'=1,\dots,2g$ and $j,j'=1,\dots,m$. 
Then the map
$$
     \T^{2gm}=\R^{2gm}/\Z^{2gm}\to\Aa/\Gg:
     [t]\mapsto A_0+\sum_{i,j}t_{ij}\alpha_{ij}
$$
is a homotopy equivalence and identifies the 
cohomology class of the $1$-form $dt_{ij}$ 
with $\mu_{\Aa/\Gg}(c_{\w_j},\gamma_i)$.
\end{proof}


\section{Relative periodic orbits}\label{sec:periodic}

Let $\G$ be a compact Lie group and $(M,\om,\mu)$ be a 
symplectic manifold with a Hamiltonian $\G$-action. 
Let $\R\times M\to\R:(t,x)\mapsto H_t(x)=H_{t+1}(x)$
be a $\G$-invariant Hamiltonian.  A {\bf relative periodic orbit}
is a pair $(x_0,g_0)$, where $g_0\in\G$ and $x_0:\R\to M$ 
is a smooth function such that
\begin{equation}\label{eq:per}
      \dot x_0(t) + X_{H_t}(x_0(t)) = 0,\qquad
      x_0(t+1) = g_0x_0(t).
\end{equation}
It follows from the $\G$-invariance of $H_t$ that
the function $t\mapsto\mu(x_0(t))$ is constant for every relative 
periodic orbit $(x_0,g_0)$. The group $\G$ acts on the space of 
relative periodic orbits by 
$$
      g^*(x_0,g_0) := (g^{-1}x_0,g^{-1}g_0g).
$$
If $g_0$ belongs to the identity component of $\G$ then
there exists a smooth function $g:\R\to\G$ such that 
\begin{equation}\label{eq:gg0}
      g(t+1)=g_0g(t).
\end{equation}
Define $x:\R\to M$ and $\xi:\R\to\g$ by 
$$
      x(t) := g(t)^{-1}x_0(t),\qquad \xi(t) := g(t)^{-1}\dot g(t).
$$
Then 
\begin{equation}\label{eq:PER}
      \dot x(t) + X_{\xi(t)}(x(t)) + X_{H_t}(x(t)) = 0,\quad
      x(t+1)=x(t),\quad \xi(t+1)=\xi(t).
\end{equation}
A solution $(x,\xi)$ of~(\ref{eq:PER}) is called {\bf contractible}
if the loop $x:\R/\Z\to M$ is contractible.  A solution $(x_0,g_0)$
of~(\ref{eq:per}) is called {\bf contractible} if there exists a smooth
path $g:\R\to\G$ satisfying~(\ref{eq:gg0}) such that the loop
$g^{-1}x_0:\R/\Z\to M$ is contractible. 

\begin{remark}\label{rmk:per1}\rm
The loop group $L\G:=\Cinf(\R/\Z,\G)$ acts on the space of 
solutions of~(\ref{eq:PER}) by 
$$
     g^*(x,\xi) := (g^{-1}x,g^{-1}\dot g+g^{-1}\xi g).
$$
If $M$ is compact then this action preserves the 
space of contractible loops. The proof uses Floer 
homology (see for example~\cite[Chapter~10]{MS2}). 
\end{remark}

\begin{remark}\label{rmk:per2}\rm
If $(x,\xi)$ is a solution of~(\ref{eq:PER}) then 
$$
      \frac{d}{dt}\mu(x(t)) + [\xi(t),\mu(x(t))] = 0.
$$
In particular, if $\mu(x(t))\in Z(\g)$ for some $t\in\R$
then the function $t\mapsto\mu(x(t))$ is constant. 
\end{remark}

\begin{theorem}\label{thm:per}
Assume~$(H1-3)$ and suppose that $\tau\in Z(\g)$ is a 
regular value of $\mu$ such that $\mu^{-1}(\tau)\ne\emptyset$. 
Then, for every time dependent $\G$-invariant 
Hamil\-tonian $H_t=H_{t+1}:M\to\R$, 
there exists a contractible relative periodic 
orbit in $\mu^{-1}(\tau)$.
\end{theorem}

\begin{corollary}[Gromov]\label{cor:per1}
Let $(M,\om)$ be a compact symplectic manifold 
such that $\inner{[\om]}{\pi_2(M)}=0$.  Then 
every time-dependent $1$-periodic Hamiltonian system
on $M$ has a contractible periodic orbit.  
\end{corollary}

\begin{proof}
Theorem~\ref{thm:per} with $\G=\{\one\}$. 
\end{proof}

\begin{corollary}\label{cor:per2}\rm
Assume~$(H1-3)$ and suppose that $\G$ is abelian.
Then, for every time dependent $\G$-in\-va\-ri\-ant 
Hamiltonian $H_t:=H_{t+1}:M\to\R$ 
and every $\tau\in\g$ such that $\mu^{-1}(\tau)\ne\emptyset$,
there exists a contractible relative periodic orbit 
in $\mu^{-1}(\tau)$.
\end{corollary}

\begin{proof}
We may assume without loss of generality that $\G=T$
is a torus and $M$ is connected.  Then there exists a subgroup 
$\HG\subset\G$, called the principal orbit type, 
such that $\HG=\G_x$ for every $x$ in an open 
dense subset of $M$  (see Audin~\cite{AUDIN}).
Since $\HG\subset\G_x$ for every $x$ it follows that
$\inner{d\mu(x)v}{\eta}=0$ for all $v\in T_xM$ 
and all $\eta\in\h:=\Lie(\HG)$. Since $M$ is connected 
this shows that the image of $\mu$ is contained in an
affine subspace $\g_0\subset\g$ parallel to $\h^\perp$. 
The assertion about the principal orbit type now shows
that $\mu(M)$ is equal to the closure of its interior relative
to $\g_0$. Hence, for every $\tau\in\mu(M)$, there exists a
sequence $\tau_\nu\in\mu(M)$ converging to $\tau$ such that
$\tau_\nu$ is a regular value of the composition
$\mu_0:M\to\g_0$ of the moment map with the projection 
onto $\g_0$.  Now apply Theorem~\ref{thm:per} to the 
action of $\G/\HG$ on $M$ to obtain, for every $\nu$,
a contractible relative periodic orbit $(x_\nu,g_\nu)$ 
in $\mu^{-1}(\tau_\nu)$.  Every such sequence has a 
convergent subsequence. 
\end{proof}

\begin{example}\label{ex:zreg}\rm
Consider a Hamiltonian action of $\U(2)$ on $(M,\om)$
which factors through the determinant $\U(2)\to S^1$.
Then the moment map has no central regular values. 
\end{example}

\begin{conjecture}\label{con:per}
Assume~$(H1)$.  Then, for every time dependent 
$\G$-in\-va\-ri\-ant Hamiltonian $H_t:=H_{t+1}:M\to\R$ 
and every $\tau\in\g$ such that $\mu^{-1}(\tau)\ne\emptyset$,
there exists a contractible relative periodic orbit 
in $\mu^{-1}(\tau)$.
\end{conjecture}

\begin{example}\label{ex:per}\rm
Hypothesis~$(H1)$ cannot be removed in 
Conjecture~\ref{con:per}.  For example, consider
the case $\G=\{\one\}$ and $M=\T^2\times\R^2$ with the 
Hamiltonian function $H(x,y)=a_1y_1+a_2y_2$, where 
$a_1$ and $a_2$ are rationally independent. 
Then there are no (relative) periodic orbits
and the moment map is not proper.
\end{example}

\begin{remark}\label{rmk:per}\rm
It should be possible to extend the techniques developed 
in this paper to the case where~$(H3)$ is not satisfied,
however, the moduli spaces will then no longer be compact.
Such an extension should give rise to a proof of 
Conjecture~\ref{con:per} under hypotheses~$(H1)$ and~$(H2)$. 
\end{remark}

\begin{remark}\label{rmk:small}\rm
The proof of Theorem~\ref{thm:per} shows that 
the result continues to hold if hypothesis~$(H3)$ is 
replaced by the condition
$$
     \int_0^1(\max H_t-\min H_t)\,dt
     \le \frac12\min_{{\rm const}\ne v:S^2\to M\atop\bar\p_J(v)=0}
     \int_{S^2}v^*\om.
$$
\end{remark}

\begin{proof}[Proof of Theorem~\ref{thm:per}]
The proof is the analogue of Gromov's argument, with 
pseudoholomorphic curves replaced by the solutions of the 
perturbed symplectic vortex equations. 

Let $P:=S^2\times\G$ be the trivial bundle and 
$B:=0\in H_2(M\times_\G\EG;\Z)$.  We prove that 
\begin{equation}\label{eq:inv-tau}
      \Phi^{M,\mu-\tau}_{0,S^2}(\alpha)
      = \int_{\mu^{-1}(\tau)/\G}\alpha
\end{equation}
for every $\alpha\in\Om^{\dim M-2\dim \G}(M\times_\G\EG)$.
To see this note that, by Lemma~\ref{le:horizontal},
every solution $(u,A)\in\Cinf(S^2,M)\times\Om^1(S^2,\g)$ 
of the unperturbed equation~(\ref{eq:vortex}) over $S^2$ 
is horizontal. Since every flat $\G$-connection on the trivial 
bundle over $S^2$ is gauge equivalent to the zero connection,
it follows that every solution of~(\ref{eq:vortex})
is gauge equivalent to a solution of the form
$$
     u(z)\equiv x,\qquad A=0.
$$
For any such solution and any almost complex structure
$J\in\Jj_\G(M,\om)$ it follows from Remark~\ref{rmk:Tabla}
(with $H=0$) that the Cauchy--Riemann operator 
$D_{u,A}:\Cinf(S^2,u^*TM)\to\Om^{0,1}(S^2,u^*TM)$ is complex
linear.  Moreover, the bundle $u^*TM\to S^2$ is a direct sum of 
complex line bundles of degree zero.  Hence it follows 
from the Riemann--Roch theorem that $D_{u,A}$ is surjective. 
Combining these observations with Proposition~\ref{prop:DD*}
we find that the operator $\Dd_{u,A}$, defined by~(\ref{eq:DuA}),
is surjective.   Now consider the setup of Section~\ref{sec:rational}
and let $(\Bb,\Ee,\Ss)$ be the $\G$-moduli problem defined 
by~(\ref{eq:B}), (\ref{eq:E}), and~(\ref{eq:S}).  
Let $(B,E,S)$ be the finite dimensional $\G$-moduli 
problem defined by 
$$
     B := E := \mu^{-1}(\tau),\qquad S \equiv 0.
$$
Since $\Dd_{u,A}$ is surjective for every $[u,A]\in\Ss^{-1}(0)$,
the obvious inclusions $B\INTO\Bb$ and $E\INTO\Ee$ define 
a morphism from $(B,E,S)$ to $(\Bb,\Ee,\Ss)$ in the sense of 
Definition~\ref{def:morphism}. Hence~(\ref{eq:inv-tau}) 
follows from the {\it (Functoriality)} and 
{\it (Transversality)} axioms of the Euler class.

By Kirwan's theorem~\cite{KIRWAN}, the homomorphism 
$$
     H_\G^*(M)\to H^*_\G(\mu^{-1}(\tau))
$$
is surjective for every $\tau\in Z(\G)$. Since $\mu^{-1}(\tau)$
is nonempty and $\G$ acts with finite isotropy, there exists a 
$\G$-invariant horizontal volume form on $\mu^{-1}(\tau)$. 
This implies that there exists a $\G$-closed equivariant 
differential form $\alpha\in\Om^{\dim M-2\dim\G}_\G(M)$ 
such that 
$$
     \int_{\mu^{-1}(\tau)/\G}\alpha \ne 0.
$$
Hence, by Lemma~\ref{le:regular} and Theorem~\ref{thm:invariant-rat},
there exists a constant $\delta>0$ such that, 
for every metric on $S^2$, every $J\in\Jj_\G(M,\om)$,
and every compactly supported Hamiltonian perturbation 
$\widehat{H}\in\Om^1(S^2,\Cinf_\G(M))$,
\begin{equation}\label{eq:nonempty}
     \|\Om_{\widehat{H}}\|\le \delta\Vol(S^2)\qquad\IMP\qquad
     \Mm_{0,S^2}(\tau;J,\widehat{H})\ne\emptyset.
\end{equation}
For $T>0$ choose a metric on $S^2=\C\cup\{\infty\}$ such that 
the map $[-T,T]\times\R/\Z\to\C:(s,t)\mapsto e^{2\pi(s+it)}$
is an isometric embedding.  Let $\rho_T:[-T,T]\to[0,1]$ be a smooth
cutoff function sucht that $\pm\dot\rho_T(s)\ge 0$ for $\pm s\ge 0$ 
and $\rho_T(s)=1$ for $|s|\le T-1$.  
Fix a compactly supported $1$-periodic $\G$-invariant
Hamiltonian function $\R/\Z\times M\to\R:(t,x)\mapsto H_t(x)$. 
On the cylinder $[-T,T]\times\R/\Z$ consider the Hamiltonian 
perturbation $\widehat{H}_T:=\rho_T(s)H_t(x)\,dt$ and extend 
it by zero to all of $S^2$. The Hofer norm of the curvature of 
$\widehat{H}_T$ is given by
$$
     \|\Om_{\widehat{H}_T}\| = 2\|H\|,\qquad
     \|H\| := \int_0^1(\max H_t-\min H_t)\,dt.
$$
Hence it follows from~(\ref{eq:nonempty}), that 
$$
     \Mm_{0,S^2}(\tau;J,\widehat{H}_T)\ne\emptyset
$$ 
for $T$ sufficiently large.  This implies that for $T\ge T_0$ there 
exist functions $u=u_T:[-T,T]\times\R/\Z\to M$ and 
$\Phi=\Phi_T,\Psi=\Psi_T:[-T,T]\times\R/\Z\to\g$ such that 
\begin{equation}\label{eq:loc-per}
\begin{array}{rcl}
    \p_su+X_\Phi(u) 
    + J\left(\p_tu+X_\Psi(u)+\rho_T(s)X_{H_t}(u)\right) &= &0,  \\
    \p_s\Psi-\p_t\Phi+[\Phi,\Psi] + \mu(u)-\tau &= & 0,
\end{array}
\end{equation}
and 
\begin{equation}\label{eq:en-per}
    \int_0^1\int_{-T+1}^{T-1}\left(
    \left|\p_tu+L_u\Psi+X_{H_t}(u)\right|^2
    + \left|\mu(u)-\tau\right|^2\right)\,dsdt
    \le 2\|H\|.
\end{equation}
The inequality~(\ref{eq:en-per}) follows from the energy 
identity in Proposition~\ref{prop:energy}. 
Choose $s_T\in[-T+1,T-1]$ such that 
$$
     \int_0^1\left(
     \left|\p_tu+L_u\Psi+X_{H_t}(u)\right|^2
     + \left|\mu(u)-\tau\right|^2\right)(s_T,t)\,dt
     \le \frac{\|H\|}{T-1}.
$$
Gauge transforming the solution at $s=s_T$ we may assume, without
loss of generality, that 
$$
     \Psi_T(s_T,t)=:\xi_T,\qquad |\xi_T|\le c,
$$ 
where $c$ is the diameter of $\G$ with respect to 
our biinvariant metric. Namely, choose $g:\R/\Z\to\G$ 
such that 
$$
     \p_tg(t) + \Psi(s_T,t)g(t) = 0,\qquad
     g(0)=\one.
$$
Then write $g(1)^{-1}=\exp(\xi_T)$, where $|\xi_T|\le c$, 
and gauge transform $u_T$ and $\Psi_T$ by the product 
$g(t)\exp(t\xi_T)$.  Now (i.e. after gauge transforming)
define $x_T:\R/\Z\to M$ by 
$$
     x_T(t) := u_T(s_T,t).
$$
Then 
$$
     \lim_{T\to\infty}
     \int_0^1\left(
     \left|\dot x_T(t)+L_{x_T(t)}\xi_T+X_{H_t}(x_T(t))\right|^2
     + \left|\mu(x_T(t))-\tau\right|^2\right)\,dt
     = 0.
$$
This shows that the $L^2$ norm of $\dot x_T$ is bounded 
and so $x_T$ is bounded and equicontinuous. Hence, by the 
Arzela-Ascoli theorem, there exists a sequence $T_i\to\infty$ 
such that $x_{T_i}$ converges uniformly, $\dot x_{T_i}$
converges weakly in $L^2$, and $\xi_{T_i}$ converges in $\g$.
The limit $(x,\xi)$ is the required solution of~(\ref{eq:PER}).
Since $x_T$ is contractible for every $T$, so is $x$.
This proves the theorem for compactly supported Hamiltonian 
functions. The general case follows by cutting off the Hamiltonian
function outside of $\mu^{-1}(\tau)$.
\end{proof}


\section{Weighted projective space}\label{sec:proj}

Consider the symplectic manifold $M=\C^n$ with 
the standard symplectic form and the $S^1$-action
$$
     \lambda x = (\lambda^{\ell_1}x_1,\dots,\lambda^{\ell_n}x_n),
$$
where $\ell_1,\dots,\ell_n$ are positive integers. 
Then a moment map is given by 
\begin{equation}\label{eq:mu-n}
     \mu_\ell(x) = -\frac{i}{2}\sum_{\nu=1}^n \ell_\nu|x_\nu|^2.
\end{equation}
Suppose that $\Sigma$ has genus $g$ and let $P\to\Sigma$
be an $S^1$-bundle of degree $d$.  Consider the complex 
line bundle
$$
     L := P\times_{S^1}\C \to \Sigma,
$$
where $S^1$ acts on $P\times\C$ by 
$$
     \lambda^*(p,\zeta) = (p\lambda,\lambda^{-1}\zeta).
$$
Then the symplectic vortex equations~(\ref{eq:vortex}) 
with $\mu=\mu_\ell$ given by~(\ref{eq:mu-n}) 
can be written in the form
\begin{equation}\label{eq:vortex-n}
     \bar\p_Au_\nu=0,\qquad
     *iF_A + \sum_{\nu=1}^n\frac{\ell_\nu|u_\nu|^2}{2} 
     = \tau,
\end{equation}
where $u_\nu$ is a section of $L^{\otimes\ell_\nu}$ for $\nu=1,\dots,n$
and $A\in\Aa(L)$ is a Hermitian connection on $L$.  
Let us denote by $\Mm_{d,g}$ the space of gauge equivalence 
classes of solutions $(u_1,\dots,u_n,A)$ of~(\ref{eq:vortex-n}). 
The moduli space is nonempty only if
$$
     \tau > \frac{2\pi d}{\Vol(\Sigma)}.
$$ 
Moreover, $\Mm_{d,g}$ has virtual dimension
$$
     \dim\Mm_{d,g} 
     = 2\left(d\sum_{\nu=1}^n\ell_\nu - (n-1)(g-1)\right) =: 2m.
$$
For $d$ sufficiently large the dimension is positive.
We write 
$
     \Phi_{d,g}^{\C^n,\mu_\ell}
     := \Phi^{\C^n,\mu_\ell+i\tau}_{B,\Sigma}
$
for the invariant in the nontrivial chamber.
Let $c\in H^2_{S^1}(\C^n;\Z)\cong H^2({\rm B}S^1;\Z)\cong\Z$ 
denote the positive generator.

\begin{theorem}\label{thm:weights}
Assume
$$
     m:=\sum_{\nu=1}^n(d\ell_\nu-g+1) + g-1\ge 0.
$$  
Then, for $\alpha=0$, $\beta=c^m$, and $\gamma=\point$,
the invariant is 
$$
     \Phi_{d,g}^{\C^n,\mu_\ell}(0,c^m;\point) 
     = \left(\sum_{\nu=1}^n\ell_\nu\right)^g
       \prod_{\nu=1}^n\ell_\nu^{-d\ell_\nu+g-1}.
$$
\end{theorem}

In the case $\ell_\nu=1$ and $d>2g-2$ Theorem~\ref{thm:weights} 
was proved by Bertram--Daskalopoulos--Wentworth~\cite{BDW}. 
In this case the invariant $\Phi_{d,g}^{\C^n,\mu_\ell}(0,c^m;\point)$
corresponds to the Gromov--Witten invariant given by counting 
holomorphic maps $u:\Sigma\to\C P^{n-1}$ of degree $d$
(with a fixed complex structure on $\Sigma$ and a generic 
Hamiltonian perturbation) that pass at $m$ distinct 
specified points on $\Sigma$ through $m$ specified 
hyperplanes in $\C P^{n-1}$. For a proof of this 
correspondence in the case $\Sigma=S^2$ 
see~\cite{GaSa}. We emphasize that in the higher genus 
case the Hamiltonian perturbation is needed in order to destroy 
the constant holomorphic maps $\Sigma\to\C P^{n-1}$
which are not regular.  If one wants to work with the 
unperturbed Cauchy--Riemann equations one has to 
work with stable maps.

\begin{proof}[Proof of Theorem~\ref{thm:weights}]
We simplify~(\ref{eq:vortex-n}) and consider 
instead the equations
\begin{equation}\label{eq:flat}
     \bar\p_Au_\nu=0,\qquad
     *iF_A = \frac{2\pi d}{\Vol(\Sigma)},\qquad
     \sum_{\nu=1}^n\ell_\nu\left\|u_\nu\right\|_{L^2}^2=1
\end{equation}
for $u_\nu\in\Om^0(\Sigma,L^{\otimes\ell_\nu})$ and $A\in\Aa(L)$.
There are two ways to establish the correspondence between
equations~(\ref{eq:vortex-n}) and~(\ref{eq:flat}).  One can use the 
action of the complexified gauge group and the Kazdan--Warner equation,
or one can show that the corresponding $S^1$-moduli problems
are homotopic. We use the latter approach.
Consider the $1$-parameter family of equations
\begin{equation}\label{eq:eps}
     \bar\p_Au_\nu=0,\qquad
     *iF_A - \frac{2\pi d}{\Vol(\Sigma)} 
         = \frac{\eps}{2}\left(\frac{1}{\Vol(\Sigma)}
       -\sum_{\nu=1}^n\ell_\nu|u_\nu|^2\right)
\end{equation}
for $0\le\eps\le1$. For $\eps=0$ this equation is equivalent 
to~(\ref{eq:flat}) and for $\eps=1$ it is equivalent 
to~(\ref{eq:vortex-n}) with $\tau=(2\pi d+1/2)/\Vol(\Sigma)$.  
Note, in particular, that $\sum_\nu\ell_\nu\|u_\nu\|_{L^2}^2=1$ 
for every $\eps>0$ and every solution of~(\ref{eq:eps}). 
Thus we may formulate the $S^1$-moduli problems as follows.
We shall not bother with Sobolev completions and formulate
the problems in terms of smooth sections.  

Fix a point $z_0\in\Sigma$ and consider
the based gauge group
$$
     \Gg_0 := \left\{g\in\Cinf(\Sigma,S^1)\,|\,g(z_0)=1\right\}.
$$
Define $\Bb$ by 
$$
     \Bb := \left\{(A,u_1,\dots,u_n)\,|\,
     A\in\Aa(L),\,u_\nu\in\Om^0(\Sigma,L^{\otimes\ell_\nu}),\,
     \sum_{\nu=1}^n\ell_\nu\left\|u_\nu\right\|_{L^2}^2=1
     \right\}\Bigg/\Gg_0.
$$
The bundle $\Hh\to\Bb$ has fibres 
$$
     \Hh_{A,u} := \Om^{0,1}(\Sigma,L^{\otimes\ell_1})
     \oplus\cdots\oplus\Om^{0,1}(\Sigma,L^{\otimes\ell_n})
     \oplus\Om^0_0(\Sigma)
$$
over $[A,u]=[A,u_1,\dots,u_n]\in\Bb$, where $\Om^0_0(\Sigma)$ 
denotes the space of smooth real valued functions of mean value 
zero.  The section $\Ss_\eps:\Bb\to\Hh$ is given by 
$$
     \Ss_\eps(A,u)
     := \left(\bar\p_Au_1,\dots,\bar\p_Au_n,
         *iF_A - \frac{2\pi d}{\Vol(\Sigma)}
         - \frac{\eps}{2}\left(\frac{1}{\Vol(\Sigma)}
         - \sum_{\nu=1}^n \ell_\nu|u_\nu|^2\right)
        \right).
$$
With appropriate Sobolev completions this gives rise to a 
homotopy of regular $S^1$-moduli problems $(\Bb,\Hh,\Ss_\eps)$.  
In fact, if $d\ell_\nu>2g-2$ for every $\nu$ then, 
by Serre duality, $\Ss_\eps$ is transverse to
the zero section for every $\eps$ and so the zero sets 
of $\Ss_\eps$ give rise to a (trivial) cobordism from the 
moduli space of solutions of~(\ref{eq:flat}) to the 
moduli space of solutions of~(\ref{eq:vortex-n}).  
In any case (even without $d\ell_\nu>2g-2$) it follows from the {\it
(Homotopy)} axiom  
for the Euler class that our invariant is given by 
\begin{equation}\label{eq:invariant1}
     \Phi^{\C^n,\mu_\ell}_{d,g}(c^m)
     = \chi^{\Bb,\Hh,\Ss_1}(\pi^*c^m)
     = \chi^{\Bb,\Hh,\Ss_0}(\pi^*c^m),
\end{equation}
where $\pi:\Bb\times_{S^1}{\rm E}S^1\to{\rm B}S^1$ denotes the
obvious projection. We shall now compute the last 
term in~(\ref{eq:invariant1}) using a localization formula for circle
actions and the index theorem for families.

Fix a reference connection $A_0\in\Aa(L)$ and consider the space 
$$
     \Aa^{\rm coul} := \left\{A\in\Aa(L)\,\Big|\,
     *iF_A=\frac{2\pi d}{\Vol(\Sigma)},\,
     d^*(A-A_0)=0\right\}
$$
of projectively flat connections in Coulomb gauge relative to $A_0$.
The group
$$
     \Gg^{\rm coul} 
     := \left\{g\in\Cinf(\Sigma,S^1)\,|\,d^*(g^{-1}dg)=0\right\}
$$
of harmonic gauge transformations acts on $\Aa^{\rm coul}$ and the 
quotient 
$
     \Aa^{\rm coul}/\Gg^{\rm coul}
$ 
is diffeomorphic to the torus $\T^{2g}$ 
(the Jacobian of degree $d$ line bundles over $\Sigma$). 
An explicit diffeomorphism can be constructed as follows. 
Choose $2g$ embedded loops $\gamma_1,\dots,\gamma_{2g}$ in $\Sigma$ 
such that 
$$
     \gamma_j\cdot\gamma_{g+j}=1
$$
and $\gamma_j\cdot\gamma_{j'}=0$ for $j'\ne j\pm g$. 
Choose a dual basis $\alpha_j:=\PD(\gamma_j)\in H^1(\Sigma)$ 
of the space of harmonic $1$-forms so that
$$
    \int_\gamma\alpha_j = \gamma_j\cdot\gamma,\qquad
    \int_\Sigma\alpha_j\wedge\alpha_{g+j} = 1.
$$
Let $\Z^{2g}\to\Gg^{\rm coul}:k\mapsto g_k$ be a group homomorphism 
such that 
$$
     \frac{1}{2\pi i}\int_{\gamma_j}g_k^{-1}dg_k=k_j
$$
for every $k\in\Z^{2g}$ and every $j\in\{1,\dots,2g\}$. 
Then the map 
$
     \R^{2g}\to\Aa^{\rm coul}:t\mapsto A_t,
$
defined by 
$$
     A_t := A_0 + \sum_{j=1}^{2g}2\pi i t_j\alpha_j,
$$
descends to a diffeomorphism $\T^{2g}\to\Aa^{\rm coul}/\Gg^{\rm coul}$. 
Note that
\begin{equation}\label{eq:Ag}
     A_{t+k} = g_k^*A_t
\end{equation}
for $t\in\R^{2g}$ and $k\in\Z^{2g}$.
Now consider the action of $\Z^{2g}$ on $\R^{2g}\times L$ via 
$$
     k\cdot(t,[p,\zeta]) := (t+k,[p,g_k(p)^{-1}\zeta]).
$$
This action gives rise to a universal
line bundle 
$$
     \LL := \frac{\R^{2g}\times L}{\Z^{2g}}\to \T^{2g}\times\Sigma.
$$
For any integer $k\in\Z$ we denote by 
$
     \LL^k = \LL\otimes\cdots\otimes\LL
$
the $k$th tensor power of $\LL$. From now on we denote by $t$ 
an equivalence class in $\T^{2g}=\R^{2g}/\Z^{2g}$.  
For $t\in\T^{2g}$ denote by $\LL_t^k\to\Sigma$ the restriction 
of $\LL^k$ to $t\times\Sigma$. By~(\ref{eq:Ag}), the bundle 
$\LL_t^k$ is equipped with a connection $A_t^k$ and hence 
with a Cauchy--Riemann operator 
$$
     \bar\p_t^k:\Om^0(\Sigma,\LL_t^k)\to\Om^{0,1}(\Sigma,\LL_t^k).
$$
Denote the topological index (as a $K$-theory class) of this family of
Cauchy-Riemann operators by 
$$
     \IND^k := \bigcup_t\{t\}\times 
     \ker\,\bar\p_t^k\ominus\coker\bar\p_t^k
     \in K(\T^{2g}).
$$
Now consider the vector bundle
$$
     \E := \LL^{\ell_1}\oplus\cdots\oplus\LL^{\ell_n}
     \to \T^{2g}\times\Sigma.
$$ 
For $t\in\T^{2g}$ denote by $\E_t$ the restriction of $\E$ 
to $\{t\}\times\Sigma$ and by 
$$
     \bar\p_t:\Om^0(\Sigma,\E_t)\to\Om^{0,1}(\Sigma,\E_t)
$$
the corresponding Cauchy--Riemann operator. 
The $L^2$-norm of a section $u\in\Om^0(\Sigma,\E_t)$ is given by
$$
     \left\|u\right\|_{L^2} := 
     \sqrt{\sum_{\nu=1}^n\ell_\nu\left\|u_\nu\right\|_{L^2}^2},
$$
where the sections $u_\nu\in\Om^0(\Sigma,\LL_t^{\ell_\nu})$
denote the components of $u$.  This gives rise to an
$S^1$-moduli problem as follows.  Define 
$$
     \B:=\left\{(t,u)\,|\,t\in\T^{2g},\,
     u_\nu\in\Om^0(\Sigma,\E_t),\,
     \left\|u\right\|_{L^2} = 1\right\}.
$$
The circle $S^1$ acts on $\B$ by
$$
     \lambda^*(t,u_1,\dots,u_n)
     := (t,\lambda^{-\ell_1}u_1,\dots,\lambda^{-\ell_n}u_n)
$$
for 
$
     \lambda\in S^1.
$
The bundle $\HH\to\B$ has fibres 
$$
     \HH_{t,u} :=  \Om^{0,1}(\Sigma,\LL_t^{\ell_1})
     \oplus\cdots\oplus\Om^{0,1}(\Sigma,\LL_t^{\ell_n})
     = \Om^{0,1}(\Sigma,\E_t)
$$
and the section $\SS:\B\to\HH$ is given by 
$$
     \SS(t,u) := \bar\p_tu.
$$
The obvious embeddings define a morphism from 
$(\B,\HH,\SS)$ to $(\Bb,\Hh,\Ss_0)$ and so, by the 
{\it (Functoriality)} axiom for $S^1$-moduli
problems and~(\ref{eq:invariant1}), we have
$$
     \Phi^{\C^n,\mu_\ell}_{d,g}(c^m)
     = \chi^{\B,\HH,\SS}(\pi_\B^*c^m),
$$
where $\pi_\B:\B\times_{S^1}{\rm E}S^1\to{\rm B}S^1$
is the projection.  Here the action 
of $S^1$ on $\B\times{\rm E}S^1$ is given by 
$
     \lambda^*(t,u,e)
     = (t,\lambda^{-\ell_1}u_1,\dots,\lambda^{-\ell_n}u_n,\lambda^{-1}e).
$

Now the $S^1$-moduli problem 
$(\B,\HH,\SS)$ satisfies the hypotheses of the localization formula
for circle actions in~\cite[Theorem~11.1]{CMS} and we get  
\begin{equation}\label{eq:invariant2}
     \Phi^{\C^n,\mu_\ell}_{d,g}(c^m)
     = \int_{\T^{2g}}\frac{1}
       {\prod_{\nu=1}^n c(\IND^{\ell_\nu},\ell_\nu)}.
\end{equation}
Here $c(\IND,\cdot)$ denotes the Chern series of 
the $K$-theory class $\IND\in K(\T^{2g})$.
It is defined by 
$$
     c(\IND(\Dd),\eta) 
     := \sum_{j\ge0}\eta^{\INDEX(\Dd)-j}c_j(\IND(\Dd)),
$$
where $\INDEX(\Dd):=\dim\,\ker\,\Dd-\dim\,\coker\,\Dd$ is the Fredholm
index. 

The right hand side of~(\ref{eq:invariant2}) can be computed
by means of the Atiyah--Singer index theorem for families
(see~\cite{SAL2}). It asserts that 
$$
     \ch(\IND^k)
     = \int_\Sigma\td(T\Sigma)\ch(\LL^k)
     \in H^*(\T^{2g}).
$$
Here the Todd class and the Chern character of a line bundle $L$ with
first Chern class $c_1(L)=x$ are defined as
$$
     \td(L) := \frac{x}{1-e^{-x}}, \qquad \ch(L) := e^x.
$$
The Todd class of $T\Sigma$ is given by 
$$
     \td(T\Sigma) = 1+(1-g)\sigma,
$$
where $\sigma\in H^2(\Sigma;\Z)$ denotes the positive generator,
represented by a volume form with respect to which $\Sigma$
has volume one.  By Lemma~\ref{L:configA}, 
the first Chern class of $\LL^k$ is given by 
$$
     c_1(\LL^k) = k\left(
     \sum_{j=1}^{2g}\alpha_j\wedge\tau_j + d\sigma
     \right)
$$
where $\alpha_j=\PD(\gamma_j)\in H^1(\Sigma)$ and $\tau_j:=[dt_j]\in
H^1(\T^{2g};\Z)$. Let us denote by 
$$
     \Om := \sum_{j=1}^g\tau_j\wedge\tau_{g+j}
$$
the cohomology class of the standard symplectic form 
on $\T^{2g}$.   Then the Chern character of $\LL^k$
is given by 
\begin{eqnarray*}
     \ch(\LL^k)
&= &
     1 + c_1(\LL^k) + \frac{c_1(\LL^k)^2}{2} \\
&= &
     1 + dk\sigma - k^2\sigma\wedge\Om + k\sum_{j=1}^{2g}\alpha_j\wedge\tau_j.
\end{eqnarray*}
Hence 
$$
     \td(T\Sigma)\ch(\LL^k)
     = 1 + (dk+1-g)\sigma - k^2\sigma\wedge\Om +
     k\sum_{j=1}^{2g}\alpha_j\wedge\tau_j, 
$$
and integration over the fibre gives 
$$
     \ch(\IND^k)
     = \int_\Sigma\td(T\Sigma)\ch(\LL^k)
     = dk+1-g - k^2\Om.
$$
This implies (cf.~\cite{SAL2}) 
$$
     c_1\left(\IND^k\right) 
     = -k^2\Om,\qquad
     c_j\left(\IND^k\right)  
     = \frac{1}{j!}c_1\left(\IND^k\right)^j,
$$
and so the Chern series of $\IND^k$ is
given by 
$$
     c\left(\IND^k,\eta\right) 
     = \eta^{dk+1-g}\exp(-\eta^{-1}k^2\Om).
$$
Since the integral of $\Om^g/g!$ over $\T^{2g}$ is one, 
we obtain from~(\ref{eq:invariant2}) that 
\begin{eqnarray*}
     \Phi^{\C^n,\mu_\ell}_{d,g}(c^m)
&= &
     \int_{\T^{2g}}\frac{1}
     {\prod_{\nu=1}^nc(\IND^{\ell_\nu},\ell_\nu)}  \\
&= &
     \prod_{\nu=1}^n\ell_\nu^{-d\ell_\nu+g-1}
     \int_{\T^{2g}}\frac{1}
     {\prod_{\nu=1}^n\exp(-\ell_\nu\Om)}  \\
&= &
     \prod_{\nu=1}^n\ell_\nu^{-d\ell_\nu+g-1}
     \int_{\T^{2g}}
     \prod_{\nu=1}^n\exp(\ell_\nu\Om)  \\
&= &
     \prod_{\nu=1}^n\ell_\nu^{-d\ell_\nu+g-1}
     \int_{\T^{2g}}
     \exp\left(\sum_{\nu=1}^n\ell_\nu\Om\right)  \\
&= &
     \left(\sum_{\nu=1}^n\ell_\nu\right)^g
     \prod_{\nu=1}^n\ell_\nu^{-d\ell_\nu+g-1}.
\end{eqnarray*}
This proves the theorem.
\end{proof}


\section{Seiberg--Witten invariants}\label{sec:sw}

In this section we explain how the Seiberg--Witten invariants
of a product 
$$
     X=\Sigma\times S
$$ 
are related to our invariants is the case where either 
$S$ or $\Sigma$ is a sphere.  The relation will be established 
by considering the symplectic vortex equations over 
$\Sigma$ with a suitable target manifold $M_S$.  
The space $M_S$ is a symplectic manifold with a circle action 
and the quotient $M_S\dslash S^1$ is the $d$-fold symmetric product 
of~$S$.  In fact, the space $M_S$ itself consists of (gauge 
equivalence classes of) solutions to the vortex equations over $S$. 
It is a special case of the socalled {\it master space for the 
vortex equations} constructed in~\cite{BDW1}.
Here is how this works.

Let $(S,J_S,\dvol_S)$ be a compact Riemann surface 
of genus $g_S$ and $L\to S$ be a complex Hermitian line bundle of degree 
$$
     \deg(L)=d>2g_S-2.
$$
For a Hermitian connection $A\in\Aa(L)$ 
and a section $\Theta\in\Om^0(S,L)$ 
consider the vortex equations
\begin{equation}\label{eq:v}
     \bar\p_A\Theta=0,\qquad
     *iF_A + \frac{|\Theta|^2}{2} 
     - \frac{1}{2\Vol(S)}\int_S|\Theta|^2\,\dvol_S
     = \frac{2\pi d}{\Vol(S)}.
\end{equation}
The gauge group $\Gg_S:=\Cinf(S,S^1)$ acts on the space 
of solutions of~(\ref{eq:v}) and the action is free
whenever $\Theta\neq 0$. Fix a point $x_0\in S$
and consider the homomorphism $\rho_0:\Gg_S\to S^1$ defined by
$$
     \rho_0(g):=\exp(-\xi(x_0))g(x_0),\qquad
     d^*d\xi = d^*(g^{-1}dg),\qquad \int_S\xi\dvol_S=0.
$$
Its kernel is the subgroup $\Gg_{S0}\subset\Gg_S$ 
of all smooth maps $g:S\to S^1$ of the form $g=g_0\exp(\xi)$,
where $g_0:S\to S^1$ is a harmonic map
that vanishes at $x_0$ and $\xi:S\to i\R$ has mean value
zero. Thus the Lie algebra of $\Gg_{S0}$ is the space
of imaginary valued functions of mean value zero:
$$
     \Lie(\Gg_{S0}) = \Om^0_0(S,i\R)
     := \left\{\xi\in\Om^0(S,i\R)\,\Big|\,
       \int_S\xi\dvol_S=0\right\}.
$$
Let us denote the space of solutions of~(\ref{eq:v}) by 
$$
     \tilde M_S:=\left\{(A,\Theta)\in\Aa(L)\times\Om^0(S,L)\,\Big|\,
     A\mbox{ and }\Theta\mbox{ satisfy }(\ref{eq:v})\right\}
$$
and the quotient by the action of $\Gg_{S0}$ by
$$
     M_S := \tilde M_S/\Gg_{S0}.
$$
The tangent space of $M_S$ at a pair
$(A,\Theta)$ can be identified with the space 
of all pairs $(\alpha,\theta)\in\Om^1(S,i\R)\times\Om^0(S,L)$
that satisfy the linearized equation
\begin{eqnarray}\label{eq:v-lin}
     \bar\p_A\theta + \alpha^{0,1}\Theta = 0,\nonumber\\
     *id\alpha + \inner{\Theta}{\theta} 
     - \frac{1}{\Vol(S)}\int_S\inner{\Theta}{\theta}\dvol_S = 0,\\
     -d^*\alpha + i\inner{i\Theta}{\theta} 
     - \frac{1}{\Vol(S)}\int_Si\inner{i\Theta}{\theta}\dvol_S = 0.
     \nonumber
\end{eqnarray}
Here the last equation asserts that the pair $(\alpha,\theta)$
belongs to the local slice of the $\Gg_{S0}$-action, i.e. 
it is $L^2$ orthogonal to the $\Gg_{S0}$-orbit of the pair
$(A,\Theta)$.  The left hand side of~(\ref{eq:v-lin}) defines
a surjective Fredholm operator from 
$
     \Om^0(S,L)\oplus\Om^1(S,i\R)
$
to
$
     \Om^{0,1}(S,L)\oplus\Om^0_0(S,\C)
$
whenever $\Theta\ne0$.  The condition $d>2g_S-2$
guarantees surjectivity also in the case $\Theta=0$. So in this case 
$M_S$ is a manifold of dimension $2d+2$. 
Unfortunately, the case 
$d>2g_S-2$ is only interesting when $S$ has genus zero (see
Remark~\ref{rmk:genus} below). 
If $d\le 2g_S-2$ the space $M_S$ has singularities 
at the points where $\Theta=0$.  
In the case $d\le 2g_S-2$ and $g_S>0$ the space $M_S$
can be desingularized by a blowup construction,
however this leads to holomorphic spheres in the ambient
space $M_S$ and so our theory does not apply in its present form.
  
A symplectic form $\om_S$ on $M_S$ is given by 
\begin{equation}\label{eq:omS}
     \om_S((\alpha,\theta),(\alpha',\theta'))
     := - \int_S\alpha\wedge\alpha' 
        + \int_S\inner{i\theta}{\theta'}\dvol_S
\end{equation}
for two solutions $(\alpha,\theta)$ and $(\alpha',\theta')$
of~(\ref{eq:v-lin}). One can think of $M_S$ as the symplectic 
quotient of the space of all pairs $(A,\Theta)$ that satisfy
$\bar\p_A\Theta=0$ by the (Hamiltonian) $\Gg_{S0}$-action. 
The linear map
$$
     (\alpha,\theta)\mapsto(*\alpha,i\theta)
$$
on the space of solutions of~(\ref{eq:v-lin}) defines
a complex structure $J_S$ on $M_S$ that is compatible with $\om_S$.
Thus $(M_S,\om_S,J_S)$ is a K\"ahler manifold.

Now the circle $S^1$ acts on $M_S$ through the constant 
gauge transformations. This action is Hamiltonian with moment map
$$
    \mu_S(A,\Theta) = -\frac{i}{2\Vol(S)}\int_S|\Theta|^2\dvol_S.
$$
The factor $1/\Vol(S)$ arises from the fact that we
identify the circle with the subgroup of $\Gg_S$ 
of constant gauge transformations and use the standard
$L^2$ metric on $\Lie(\Gg_S)=\Om^0(S,i\R)$ to define the
moment map as a function with values in the Lie algebra, 
and not its dual.

Let us recall some standard facts about the space $M_S$ (see~\cite{BRADLOW1,GP}). 
There is a one-to-one correspondence between Hermitian connections 
$A\in\Aa(L)$ and holomorphic structures on $L$ via $A\mapsto\bar\p_A$.
Moreover, the Kazdan--Warner equation shows that every pair 
$(A,\Theta)\in\Aa(L)\times\Cinf(S,L)$ such that 
$\bar\p_A\Theta=0$ is complex gauge equivalent to a solution 
of~(\ref{eq:v}) by a gauge transformation of the form $g=e^f$
where $f:S\to\R$ has mean value zero (see Proposition~\ref{prop:KWS}).  
Hence the space $M_S$ can be identified with the space of $\Gg_{S0}^c$-gauge 
equivalence classes of the space of pairs $(A,\Theta)$ that satisfy 
$\bar\p_A\Theta=0$. For $d>2g_S-2$ it follows from Serre duality that
this is a vector bundle over the Picard variety of holomorphic bundles
of degree $d$ over $S$, ${\rm Pic}^d(S)\cong T^{2g}$, with fibre
$\C^{d+1-g_S}$. The circle acts trivially on the base and by the
standard action on the fibres.

This shows that the triple $(M_S,\om_S,\mu_S)$ satisfies hypotheses
$(H1-3)$, namely, $\mu_S$ is proper, the moment map is 
convex at infinity, and $\pi_2(M_S)=0$.  Moreover
every nonzero imaginary number is a regular value
of $\mu_S$ and $S^1$ acts freely on the preimage under
$\mu_S$.  The quotient is nonempty if and only if the imaginary
part is negative. 

\begin{remark}\rm
The symplectic quotient
$$
      \overline{M}_S := M_S\dslash S^1(-i/2\Vol(S))
$$
is a bundle over $T^{2g}$ with fibre $\C P^{d-g_S}$. On the other
hand, this quotient is the space of effective divisors of degree $d$
on $S$, so $\overline{M}_S \cong \Sym^d(S)$ is the $d$-fold symmetric
product of $S$. 
\end{remark}

The next theorem states that the invariants 
of the triple $(M_S,\om_S,\mu_S)$ for a Riemann
surface $\Sigma$ agree with the Seiberg--Witten
invariants of the product $\Sigma\times S$. 
We denote by $c\in H^2({\rm B}S^1;\Z)$ 
the positive generator and by 
$\pi_S:M_S\times_{S^1}{\rm E}S^1\to{\rm B}S^1$
the obvious projection.  For a nonnegative integer
$
     k \in Z
$
denote by
$$
     \Phi^{M_S,\mu_S}_{k,\Sigma}
     :=\Phi^{M_S,\mu_S+i/2}_{k,\Sigma}
$$
the invariant in the nontrivial chamber. Let 
$$
     E_{k,d}\to\Sigma\times S
$$ 
be the complex line bundle which has degree $k$ over $\Sigma$ and
degree $d$ over $S$, and denote by
$\gamma_{k,d}\in\Spin^c(\Sigma\times S)$ the spin$^c$ 
structure obtained by twisting the standard spin$^c$ structure
$\gamma_0$ (associated to the complex structure) by $E_{k,d}$.
If both $\Sigma$ and $S$ have positive genus the four-manifold
$\Sigma\times S$ has $b^+>1$ and carries a well-defined
Seiberg--Witten invariant 
$$
      SW_{\Sigma\times S}:\Spin^c(\Sigma\times S)\to\Z.
$$
If $\Sigma$ or $S$ is the sphere then $b^+=1$, so 
there are two chambers for the Seiberg--Witten invariants. 
In this case we denote by $\SW_{\Sigma\times S}$ 
the Seiberg--Witten invariant in the positive chamber, 
where ``positive'' is defined in the proof of
Theorem~\ref{thm:SW} below.  A result similar to the next
theorem was proved in~\cite{OT}.

\begin{theorem}\label{thm:SW}
Let $S$ and $\Sigma$ be a compact Riemann surfaces
of genera $g_S$ and $g_\Sigma$, respectively, and 
$k$, $d$ be nonnegative integers such that  
$$
     m := d(1-g_\Sigma) + k(1-g_S) + dk \geq 0,\qquad d > 2g_S-2.
$$
Then 
$$
    \Phi^{M_S,\mu_S}_{k,\Sigma}(0,\pi_S^*c^m;\point)
    = \SW_{\Sigma\times S}(\gamma_{k,d}).
$$
If $m<0$ then both invariants are zero. 
\end{theorem}

\begin{corollary}[\cite{LL,OO}]\label{cor:ruled}
Let $S$ be the Riemann sphere, $\Sigma$ a compact Riemann surface of
genus $g_\Sigma$ and $k$, $d$ be nonnegative integers such that 
$$
     m := d(1-g_\Sigma) + (d+1)k \geq 0.
$$
Then 
$$
     \SW_{\Sigma\times S}(\gamma_{k,d}) = (d+1)^{g_\Sigma}.
$$
\end{corollary}

\begin{proof}
Since $S$ is the Riemann sphere the manifold $M_S$ is 
diffeomorphic to $\C^{d+1}$ as a K\"ahler manifold with an $S^1$
action. Hence the result follows
from Theorem~\ref{thm:SW} and Theorem~\ref{thm:weights}
with $n=d+1$ and $\ell_1=\cdots=\ell_n=1$.
\end{proof}

\begin{corollary}\label{cor:sym}\rm
Let $S$ and $\Sigma$ be a compact Riemann surfaces
of genera $g_S$ and $g_\Sigma$, respectively, and 
$k$, $d$ be nonnegative integers such that  
$$
     m := d(1-g_\Sigma) + k(1-g_S) + dk \geq 0,\qquad 
     d > 2g_S-2,\qquad 
     k > 2g_\Sigma -2.
$$
Then 
$$
      \Phi^{M_S,\mu_S}_{k,\Sigma}(0,\pi_S^*c^m;\point)
      = \Phi^{M_\Sigma,\mu_\Sigma}_{d,S}(0,\pi_\Sigma^*c^m;\point),
$$
where $M_S$ is associated to a bundle of degree $d$ over $S$
via~(\ref{eq:v}) and $M_\Sigma$ is defined analogously, 
with $S$ and $d$ replaced by $\Sigma$ and $k$.
\end{corollary}

\begin{proof}
Interchange the roles of $\Sigma$ and $S$ in the 
proof of Theorem~\ref{thm:SW}. 
\end{proof}

\begin{remark}\label{rmk:genus}\rm
The statements of Theorem~\ref{thm:SW} and Corollary~\ref{cor:sym}
are only interesting when one of the two surfaces has genus zero.
Otherwise both invariants are zero. To see this, suppose that both
genera are positive. Then $\Sigma\times S$ is a minimal K\"ahler 
surface with $b^+>1$. It follows (see for example~\cite{SAL2}) 
that the Seiberg-Witten invariant is nonzero only for the 
canonical spin$^c$ structure and its dual, 
i.e.~for $d=k=0$ or $d=2g_S-2$ and $k=2g_\Sigma-2$. 
These cases are excluded by our hypotheses.   
\end{remark}

\begin{proof}[Proof of Theorem~\ref{thm:SW}]
Our proof follows the argument outlined in~\cite{CGS}.
The Seiberg--Witten equations for the spin$^c$
structure $\gamma_{k,d}$ on $X := \Sigma\times S$ have the form
\begin{equation}\label{eq:SW1}
      \bar\p_B\Theta_0 + {\bar\p_B}^*\Theta_2 = 0,\qquad
      F_B^{0,2} - \inner{\Theta_0}{\Theta_2} = 0,
\end{equation}
\begin{equation}\label{eq:SW2}
      i(F_B)_\Om + \frac{|\Theta_0|^2 - |\Theta_2|^2}{2} = \tau,
\end{equation}
where $\tau$ is a real number, 
$B\in\Aa(E)$ is a connection on $E := E_{k,d}$, 
$\Theta_0\in\Om^{0,0}(X,E)$, 
and $\Theta_2\in\Om^{0,2}(X,E)$.  
Here we denote by $\inner{\cdot}{\cdot}$ 
a Hermitian inner product on $E$, i.e. the inner product 
takes values in $\C$, it is complex anti-linear in the first variable
and complex linear in the second variable. 
In the last term the function $\Om^2(X,i\R)\to\Om^0(X,i\R)$,
$\eta\mapsto\eta_\Om$ is defined by 
$$
     \eta_\Om := *(\eta\wedge\Om),
$$
where $p_\Sigma:X\to\Sigma$ and 
$p_S:X\to S$ denote the projections,
$$
     \Om := p_\Sigma^*\dvol_\Sigma + p_S^*\dvol_S\in\Om^2(X)
$$
denotes the symplectic form, and $*$ denotes 
the Hodge $*$-operator on $X$. 

In the K\"ahler case it follows from~(\ref{eq:SW1})
that either $\Theta_0$ or $\Theta_2$ vanishes. 
The positive chamber for the Seiberg--Witten 
invariants corresponds to the condition
\begin{equation}\label{eq:taudk}
     \tau > \frac{2\pi k}{\Vol(\Sigma)} + \frac{2\pi d}{\Vol(S)}.
\end{equation}
In the K\"ahler case this condition implies $\Theta_2=0$. 

The $S^1$-moduli problem associated to 
equations (\ref{eq:SW1}) and (\ref{eq:SW2}) is defined as follows.  
As in the proof of Theorem~\ref{thm:weights} we shall
not explain the (obvious) Sobolev completions and describe 
the problem in terms of smooth data. Fix a point 
$(z_0,x_0)\in\Sigma\times S$ and denote by $\Gg_{X0}$ the 
based gauge group of all smooth functions $g:X\to S^1$
such that the restriction of $g$ to $\{z_0\}\times S$
belongs to the subgroup $\Gg_{S0}\subset\Gg_S$ determined by 
the point $x_0$:
$$
     \Gg_{X0} := \left\{g\in\Cinf(\Sigma\times S,S^1)\,|\,
     g|_{\{z_0\}\times S}\in\Gg_{S0}\right\}.
$$
Then the base $\Bb^{\SW}$ is the quotient 
$$
     \Bb^\SW
     := \frac{\left\{(B,\Theta_0,\Theta_2)\in
        \Aa(E)\times\Om^0(X,E)\times\Om^{0,2}(X,E)\,|\,
        (\ref{eq:SW2})\right\}}
        {\Gg_{X0}},
$$
the bundle $\Ee^\SW\to\Bb^\SW$ is given by
$$
     \Ee^\SW
     := \frac{\Bb^\SW\times\left(\Om^{0,1}(X,E)\oplus\Om^{0,2}(X)\right)}
        {\Gg_{X0}}
$$
and the section $\Ss^\SW:\Bb^\SW\to\Ee^\SW$ is given by 
$$
     \Ss^\SW(B,\Theta_0,\Theta_2) 
     := \left(\bar\p_B\Theta_0+\bar\p_B^*\Theta_2,
        F_B^{0,2} - \inner{\Theta_0}{\Theta_2}\right).
$$
With appropriate Sobolev completions this is a regular $S^1$-moduli problem
in the sense of Definition~\ref{def:moduli} and the 
Seiberg--Witten invariant can be expressed in the form
$$
     \SW_{\Sigma\times S}(\gamma_{k,d})
     = \chi^{\Bb^\SW,\Ee^\SW,\Ss^\SW}(\pi_{SW}^*c^m),
$$
where $m = d(1-g_\Sigma) + k(1-g_S) + dk$, 
$
     \pi_{SW}:\Bb^\SW\times_{S^1}{\rm E}S^1\to{\rm B}S^1
$
denotes the projection, and $c\in H^2({\rm B}S^1;\Z)$ 
is the positive generator. 

Now let us examine the symplectic vortex equations with values in $M_S$.
Let $P\to\Sigma$ be a principal $S^1$-bundle of degree $k$.
Consider the associated bundle 
$$
     L_P:=P\times_{S^1}\C\to\Sigma
$$
and denote by 
$$
     E_P := p_\Sigma^*L_P\otimes p_S^*L
$$
the corresponding bundle over $X=\Sigma\times S$, where 
$p_\Sigma:X\to\Sigma$ and $p_S:X\to S$
denote the projections.  In explicit
terms it can be represented as the quotient 
$$
     E_P := \frac{P\times L}{S^1},\qquad
     \lambda^*((z,p),(x,v)) := ((z,p\lambda),(x,\lambda^{-1}v)).
$$
This vector bundle has degree $k$ over $\Sigma$ and 
degree $d$ over $S$ and hence is isomorphic to $E$.
Henceforth we shall drop the subscript $P$ and write 
$E:=E_P$. 

The space $\Aa(P)\times\Cinf_{S^1}(P,M_S)$ embeds
into the space $\Aa(E)\times\Om^0(X,E)$ as follows.  
A connection $A_\Sigma$ on $P$ determines a connection 
$p_\Sigma^*A_\Sigma$ on $p_\Sigma^*L_P$.  An $S^1$-equivariant function
$u:P\to M_S$ consists of an $S^1$-invariant function 
$A:P\to\Aa(L)$ and an $S^1$-equivariant function 
$\Theta:P\to\Om^0(S,L)$.  The latter can be interpreted as a 
section of $E$ and the two connections together determine
a connection
$$
      p_\Sigma^*A_\Sigma\otimes 1 + 1\otimes p_S^*A \in\Aa(E).
$$
To understand this correspondence better let us choose 
holomorphic local coordinates $s+it\in U\subset\C$ on 
$\Sigma$ and a trivialization of $P$ along this coordinate chart. 
In such a trivialization the connection $A_\Sigma$ 
has the form $\Phi_1\,ds+\Psi_1\,dt$ where 
$
     \Phi_1,\Psi_1:U\to i\R.
$ 
The function $u$ is a map 
$U\to\Aa(L)\times\Om^0(S,L)$ denoted by 
$$
     U\to\Aa(L):(s,t)\mapsto A(s,t),\qquad
     U\to\Om^0(S,L):(s,t)\mapsto\Theta(s,t).
$$
The corresponding connection on $E$ is given in this local
frame by $A(s,t)+\Phi_1(s,t)\,ds+\Psi_1(s,t)\,dt$.
The pair $(A_\Sigma,u)$ satisfies equation~(\ref{eq:vortex})
if and only if there exist functions 
$
     \Phi_0,\Psi_0:U\to\Om^0_0(S,i\R)
$
such that, for all $s$ and $t$,
\begin{eqnarray}\label{eq:vortex-MS0}
     \bar\p_A\Theta &= &0, 
     \nonumber \\
     \p_s\Theta + (\Phi_0+\Phi_1)\Theta
     + i\left(\p_t\Theta + (\Psi_0+\Psi_1)\Theta\right) &= & 0,
     \nonumber \\
     \p_sA-d\Phi_0 + *_S(\p_tA-d\Psi_0) &= & 0, \\
     *_SiF_A + \frac{1}{2}|\Theta|^2 
     - \frac{1}{2\Vol(S)}\int_S|\Theta|^2\,\dvol_S 
     &= & \frac{2\pi d}{\Vol(S)},\nonumber \\
     \lambda^{-2}(\p_s\Psi_1-\p_t\Phi_1)
     - \frac{i}{2\Vol(S)}\int_S|\Theta|^2\,\dvol_S &= & 
     -i\left(\tau-\frac{2\pi d}{\Vol(S)}\right).
     \nonumber
\end{eqnarray}
Here $*_S$ denotes the Hodge $*$-operator on $S$ and 
$\lambda:U\to(0,\infty)$ represents the volume form 
$\lambda^2\,ds\wedge dt$ on $\Sigma$. The factor 
$\lambda^{-2}$ in the last term arises from the Hodge
$*$-operator on $\Sigma$. The functions $\Phi_0$ and $\Psi_0$ are
needed to project the terms in the second and third equation onto the
quotient by the gauge group $\Gg_{S0}$. Let us abbreviate 
$$
     \Phi:=\Phi_0+\Phi_1,\qquad\Psi:=\Psi_0+\Psi_1.
$$
Then~(\ref{eq:vortex-MS0}) can be written in the form 
\begin{eqnarray}\label{eq:vortex-MS1}
     \bar\p_A\Theta &= &0, 
     \nonumber \\
     \p_s\Theta + \Phi\Theta
     + i\left(\p_t\Theta + \Psi\Theta\right) &= & 0,
     \nonumber \\
     \p_sA-d\Phi + *_S(\p_tA-d\Psi) &= & 0, \\
     *_SiF_A + \frac{1}{2}|\Theta|^2 
     - \frac{1}{2\Vol(S)}\int_S|\Theta|^2\,\dvol_S 
     &= & \frac{2\pi d}{\Vol(S)},\nonumber \\
     \frac{i\lambda^{-2}}{\Vol(S)}
     \int_S\left(\p_s\Psi-\p_t\Phi\right)\,\dvol_S
     + \frac{1}{2\Vol(S)}\int_S|\Theta|^2\dvol_S &= & 
     \tau-\frac{2\pi d}{\Vol(S)}
     \nonumber
\end{eqnarray}
(cf.~\cite{CGS}). Now consider the connection 
$$
     B := A(s,t) + \Phi(s,t)\,ds + \Psi(s,t)\,dt
$$
on $E$ and think of $\Theta$ as a section of $E$.
Then the first two equations in~(\ref{eq:vortex-MS1}) are equivalent 
to $\bar\p_B\Theta=0$. The curvature of $B$ is the 
$2$-form
$$
     F_B = F_A + ds\wedge(\p_sA-d\Phi)  
           + dt\wedge(\p_tA-d\Psi)
           + (\p_s\Psi-\p_t\Phi)\,ds\wedge dt.
$$
The third equation in~(\ref{eq:vortex-MS1}) asserts that 
$F_B^{0,2}=0$. The last two equations can be written in the form
$$
     \frac{i\lambda^{-2}}{\Vol(S)}
     \int_S\left(\p_s\Psi-\p_t\Phi\right)\,\dvol_S
     + *_SiF_A + \frac{1}{2}|\Theta|^2 
     = \tau,
$$
or in terms of the connection $B$,
\begin{equation}\label{eq:MS1}
     \frac{1}{\Vol(S)}\int_S *iF_B
     + *i(F_B\wedge\dvol_\Sigma) + \frac{|\Theta|^2}{2} 
     = \tau.
\end{equation}
Here the integral denotes integration over the fibre. 
Hence a pair $(A_\Sigma,u)\in\Aa(P)\times\Cinf_{S^1}(P,M_S)$
satisfies the symplectic vortex equation~(\ref{eq:vortex}) if and only
if the corresponding pair $(B,\Theta)\in\Aa(E)\times\Om^0(X,E)$
satisfies equation (\ref{eq:MS1}) and
\begin{equation}\label{eq:MS2}
     \bar\p_B\Theta=0,\qquad
     F_B^{0,2}=0.
\end{equation}
Integration of equation (\ref{eq:MS1}) over $\Sigma$ yields
$$
     \frac{2\pi k}{\Vol(\Sigma)} + \frac{2\pi d}{\Vol(S)} +
     \frac{1}{2\Vol(\Sigma)}\int_\Sigma|\Theta|^2\dvol_\Sigma = \tau.
$$
If~(\ref{eq:taudk}) holds then $\Theta\ne0$
for every solution of~(\ref{eq:MS1}).  
Hence equations (\ref{eq:MS1}) and (\ref{eq:MS2}) give rise to an 
$S^1$-moduli problem as follows. The space $\Bb^0$ is the quotient
$$
     \Bb^0
     := \frac{\left\{(B,\Theta)\in\Aa(E)\times\Om^{0}(X,E)\,|\,
        (\ref{eq:MS1})\right\}}
        {\Gg_{X0}},
$$
the bundle $\Ee^0\to\Bb^0$ is given by 
$$
     \Ee^0
     := \frac{\Bb^0\times\Om^{0,1}(X,E)}{\Gg_{X0}},
$$
and the section $\Ss^0:\Bb^0\to\Ee^0$ is
$$
     \Ss^0(B,\Theta) := \bar\p_B\Theta.
$$
Since $\Theta\ne0$
the first equation in~(\ref{eq:MS1}) implies the second equation.
Hence the zero set of $\Ss^0$ is the space of gauge equivalence 
classes of solutions $(B,\Theta)$ of~(\ref{eq:MS1})
and~(\ref{eq:MS2}).

At first glance $\Ss^0$ doesn't look like a Fredholm section.
Note, however, that $\Ss^0$ is a two-dimensional Cauchy--Riemann 
operator in disguise.  The condition $\bar\p_B\Theta=0$
assserts, at the same time, that the restriction of 
$(B,\Theta)\in\Bb^0$ to every slice $\{z\}\times S$
belongs to the finite dimensional manifold $M_S$
and that, as a function $P\to M_S$, this map is a solution
of the (two dimensional) symplectic vortex equations. 
Hence, with appropriate Sobolev completions, 
the triple $(\Bb^0,\Ee^0,\Ss^0)$ is a regular $S^1$-moduli problem and 
$$
     \Phi^{M_S,\mu_S}_{k,\Sigma}(0,\pi_S^*c^m;\point)
     = \chi^{\Bb^0,\Ee^0,\Ss^0}(\pi_0^*c^m),
$$
where $\pi_0:\Bb^0\times_{S^1}{\rm E}S^1\to{\rm B}S^1$ is the
projection. 

A morphism from $(\Bb^0,\Ee^0,\Ss^0)$ to $(\Bb^\SW,\Ee^\SW,\Ss^\SW)$
can be defined as follows. The group of complex gauge 
transformations $g:X\to\C^*$ acts on the space of 
solutions of~(\ref{eq:MS1}) via
$$
     g^*(B,\Theta) 
     := (B + g^{-1}\bar\p g - \bar g^{-1}\p \bar g,g^{-1}\Theta).
$$
Given a pair $(B,\Theta)\in\Bb^0$, we look for a complex gauge
transformation of the form $g=e^f$, where $f:X\to\R$, 
such that the triple $(B_f,\Theta_f,0):=(e^f)^*(B,\Theta,0)$, 
given by
$$
     B_f = B + \bar\p f-\p f,\qquad
     \Theta_f = e^{-f}\Theta_0,
$$
satisfies equation (\ref{eq:SW2}): 
$$
     i(F_{B_f})_\Om + \frac{|\Theta_f|^2}{2} = \tau.
$$
A short computation yields
$$
     2i\p\bar\p f = -d_\Sigma^*d_\Sigma f\,\dvol_\Sigma -
     d_S^*d_Sf\,\dvol_S, 
$$
where $d_S:\Om^0(S)\to\Om^1(S)$ and 
$d_\Sigma:\Om^0(\Sigma)\to\Om^1(\Sigma)$ are the respective
differentials and $d_S^*$ and $d_\Sigma^*$ their
$L^2$-adjoints. Therefore $(2i\p\bar\p f)_\Om = -d^*df$, and equation
(\ref{eq:SW2}) for $(B_f,\Theta_f,0)$ is 
equivalent to the {\it Kazdan-Warner equation}
$$
     -d^*df + \frac{|\Theta|^2}{2}e^{-2f} = \tau-i(F_B)_\Om.
$$
It follows from the theorem of Kazdan and Warner
(\cite{KW}, see also Appendix~\ref{app:KW})
that this equation has a unique solution $f$ whenever
$$
     \tau 
     > \frac{1}{\Vol(\Sigma)\Vol(S)} 
       \int_{\Sigma\times S}iF_{B_0}\wedge\Om
     = \frac{2\pi k}{\Vol(\Sigma)}
       + \frac{2\pi d}{\Vol(S)}
$$
(see~(\ref{eq:taudk})). So we have constructed a map
$$
     \Bb^0\to\Bb^\SW:(B,\Theta) \mapsto (B_f,\Theta_f,0).
$$
We claim that the image of this map is the submanifold
of all triples of the form $(B,\Theta,0)\in\Bb^\SW$.

A left inverse $\Bb^\SW\to\Bb^0$ can be constructed 
as follows.  Given a triple $(B,\Theta,0)\in\Bb^\SW$ 
we must find a complex gauge
transformation of the form $g=e^f$, where $f:X\to\R$, 
such that the pair $(B_f,\Theta_f)$, given by
$$
     B_f := B + \bar\p f-\p f,\qquad
     \Theta_f := e^{-f}\Theta,
$$
satisfies~(\ref{eq:MS1}):
$$
     \frac{1}{\Vol(S)}\int_S *iF_{B_f}
     + *i(F_{B_f}\wedge\dvol_\Sigma) + \frac{|\Theta_f|^2}{2} 
     = \tau.
$$
This translates into the equation
$$
     - d_\Sigma^*d_\Sigma f_\Sigma
     - d_S^*d_Sf + e^{-2f}\frac{|\Theta|^2}{2}
     = \tau - \frac{1}{\Vol(S)}\int_S *iF_B
       - *i(F_B\wedge\dvol_\Sigma),
$$
where 
$$
     f_\Sigma:=\frac{1}{\Vol(S)}\int_Sf\,\dvol_S:\Sigma\to\R.
$$
By Theorem~\ref{thm:KW} of the appendix, this equation has a unique solution 
$f\in C^0(\Sigma,W^{2,p}(S))$. If $B$ and $\Theta$ are smooth
one checks easily that $f$ is smooth.  This shows that for every pair
$(B,\Theta,0)\in\Bb^\SW$ there exists a unique complex gauge 
transformation of the form $g=e^f$ such that 
$g^*(B,\Theta)\in\Bb^0$. That this map is a left inverse of
the map $\Bb^0\to\Bb^\SW$ follows from the uniqueness statement in
Theorem~\ref{thm:KW}: Let $(B,\Theta)\in\Bb^0$ and $g=e^f$ be a
complex gauge transformation such that $g^*(B,\Theta)\in\Bb^0$. Then
$f$ satisfies the equation
$$
     - d_\Sigma^*d_\Sigma f_\Sigma
     - d_S^*d_Sf + e^{-2f}\frac{|\Theta|^2}{2}
     = \frac{|\Theta|^2}{2},
$$
and $f=0$ by uniqueness. 

It follows that the map $\Bb^0\to\Bb^\SW$ 
defines an embedding of Fr\'echet manifolds 
and lifts naturally to an embedding 
of $\Ee^0$ into $\Ee^\SW$ which intertwines 
the two sections and idenitifies the kernels and 
cokernels of the linearized operators along the zero 
set of $\Ss^0$. The proof of~\cite[Theorem~7.4]{CMS} 
shows that there exists a finite dimensional reduction 
$(B^0,E^0,S^0)$ of $(\Bb^0,\Ee^0,\Ss^0)$ in the smooth
category (not involving Sobolev completions). 
The composition of the inclusion $B^0\to\Bb^0$
with the inclusion $\Bb^0\to\Bb^\SW$ (and of their 
lifts to the vector bundles) now defines a morphism of 
$S^1$-moduli problems as in Definition~\ref{def:morphism}.  
With this established, the result follows from 
the {\it (Functoriality)} axiom for the Euler class.  
\end{proof}

\begin{remark}\label{rmk:lefschetz}\rm
There should be an analogue of Theorem~\ref{thm:SW}
in the case where the product $\Sigma\times S$ 
is replaced by a topological Lefschetz fibration
$X\to S^2$ on a symplectic manifold~\cite{DON1,DON2}.
Here $\Sigma$ should be replaced by $S^2$ and $S$ by the 
generic fibre of $X$.  To carry this out one has to overcome 
several major technical difficulties. The interesting case is 
where the degree $d$ of the bundle over the fibre satisfies
$d\le 2g_S-2$, and so the space $M_S$ has singularities.
Moroever, one has to deal with the singularities
of the fibration as in~\cite{DONS}.  In addition, 
the complex techniques with the Kazdan--Warner equation
only work in the K\"ahler case.  In the nonintegrable 
case the correspondence between the Seiberg--Witten 
equations~(\ref{eq:SW1}), (\ref{eq:SW2}) and the 
symplectic vortex equations~(\ref{eq:MS1}), (\ref{eq:MS2})
is much more subtle and requires a hard adiabatic limit analysis
as in the proof of the Atiyah--Floer conjecture~\cite{DS}
or as in~\cite{GaSa} (see~\cite{SAL1} for an outline 
of the Seiberg--Witten analogue). If this program can be 
carried out then, combined with the work of Donaldson--Smith 
in~\cite{DONS}, it might lead to an alternative proof of 
Taubes' theorem~\cite{TAUBES1,TAUBES2,TAUBES3} about the relation 
between the Seiberg--Witten and the Gromov invariants. 
\end{remark}


\appendix

\section{The coupled Kazdan--Warner equation}\label{app:KW}

Let $(\Sigma,J_\Sigma,\dvol_\Sigma)$ and $(S,J_S,\dvol_S)$
be compact connected Riemann surfaces. Fix a constant $p>1$. 
Given a function $u\in L^p(\Sigma\times S)$ we define 
$u_\Sigma\in L^p(\Sigma)$ by 
$$
    u_\Sigma(z) := \frac{1}{\Vol(S)}\int_Su(z,\cdot)\,\dvol_S
$$
for $z\in\Sigma$. In the following we shall denote 
by $d_S:\Om^0(S)\to\Om^1(S)$ and 
$d_\Sigma:\Om^0(\Sigma)\to\Om^1(\Sigma)$ the respective
differentials and by $d_S^*$ and $d_\Sigma^*$ their
$L^2$-adjoint operators. 

\begin{theorem}\label{thm:KW}
Let $p>1$ and $f,h\in C^0(\Sigma\times S)$
such that 
$$
     h\ge 0,\qquad
     \int_{\Sigma\times S}h > 0,\qquad
     \int_{\Sigma\times S}f > 0.
$$
Then there exists a unique function 
$u\in C^0(\Sigma,W^{2,p}(S))$ such that 
$u_\Sigma\in W^{2,p}(\Sigma)$ and
\begin{equation}\label{eq:KW}
    d_\Sigma^*d_\Sigma u_\Sigma + d_S^*d_Su
    + e^uh = f.
\end{equation}
Moreover, if $h$ and $f$ are smooth then so is the 
unique solution $u$ of~(\ref{eq:KW}). 
\end{theorem}

The proof of the theorem is based on a lemma and two propositions.

\begin{lemma}\label{le:KW0}
Let $S$ be a compact Riemann surface.
Then there exists a constant $c_S>0$ such that the 
following holds. Let $p>1$, $C\ge 0$, and $0<a\le A$.
If $h\in C^0(S)$ and $u\in W^{2,p}(S)$ satisfy 
$$
     h\ge 0,\qquad
     h_0:=\frac{1}{\Vol(S)}\int_Sh\,\dvol_S>0,
$$
and
\begin{equation}\label{eq:KW0}
    a-e^uh
    \le d_S^*d_Su \le 
    A-f^{-1}(e^u)h
\end{equation}
almost everywhere, where $f(r):=re^{Cr}$, then 
$$
    \log\left(\frac{a}{\left\|h\right\|_{L^\infty}}\right)
    \le u
    \le \log\left(\frac{A}{h_0}\right)
        + \frac{A}{h_0}
        \left(C+c_S\left\|h\right\|_{L^\infty}\right).
$$
\end{lemma}

\begin{proof}
Assume first that $u$ and $h$ are smooth. 
Choose $c_S>0$ such that
$$
     \int_Sv\,\dvol_S=0\qquad\IMP\qquad
     2\left\|v\right\|_{L^\infty}\le c_S\left\|d_S^*d_Sv\right\|_{L^\infty}
$$
for every $v\in\Cinf(S)$. Let $v\in\Cinf(S)$ 
be the unique solution of the equation
$$
     d_S^*d_Sv = h-h_0,\qquad \int_Sv\,\dvol_S=0.
$$
Since $\left\|h-h_0\right\|_{L^\infty}\le\left\|h\right\|_{L^\infty}$
it follows that 
$$
     \max_S v-\min_S v
     \le c_S\left\|h\right\|_{L^\infty}.
$$
Now fix a constant $\eps>0$ and denote 
$$
     w_\eps := \frac{A+\eps}{h_0}v + u.
$$
Choose $x_\eps\in S$ such that 
$
     w_\eps(x_\eps) = \sup_S w_\eps.
$
Then
\begin{eqnarray*}
     0
&\le &
     d_S^*d_Sw_\eps(x_\eps)   \\
&= &
     \frac{A+\eps}{h_0}d_S^*d_Sv(x_\eps)
     + d_S^*d_Su(x_\eps)   \\
&\le &
     \frac{A+\eps}{h_0}\left(h(x_\eps)-h_0\right)
     + A - h(x_\eps)f^{-1}\left(e^{u(x_\eps)}\right)  \\
&= &
     - \eps 
     + h(x_\eps)\left(\frac{A+\eps}{h_0} 
       - f^{-1}\left(e^{u(x_\eps)}\right)\right).
\end{eqnarray*}
It follows that $h(x_\eps)>0$ and
$
     f^{-1}\left(e^{u(x_\eps)}\right)
     < (A+\eps)/h_0. 
$
Since $f$ is strictly monotone, this implies
$$
     u(x_\eps)
     < \log\left(f\left(\frac{A+\eps}{h_0}\right)\right)
     = \log\left(\frac{A+\eps}{h_0}\right)
       + C\frac{A+\eps}{h_0}.
$$
Since $w_\eps(x)\le w_\eps(x_\eps)$ for all $x\in S$
it follows that
\begin{eqnarray*}
    u(x) 
&\le &
    u(x_\eps) 
    + \frac{A+\eps}{h_0}\left(v(x_\eps)-v(x)\right)  \\
&\le &
    \log\left(\frac{A+\eps}{h_0}\right)
    + \frac{A+\eps}{h_0}\left(
      C + \max_Sv-\min_Sv\right)  \\
&\le &
    \log\left(\frac{A+\eps}{h_0}\right)
    + \frac{A+\eps}{h_0}\left(C+c_S\left\|h\right\|_{L^\infty}\right).
\end{eqnarray*}
This holds for every $\eps>0$ and every $x\in S$.
Hence
$$
    \sup_Su\le \log\left(\frac{A}{h_0}\right)
    + \frac{A}{h_0}\left(C+c_S\left\|h\right\|_{L^\infty}\right).
$$
To prove the first inequality we choose $x_0\in S$ 
such that 
$
     u(x_0) = \inf_S u.
$
Then
$$
     0
     \ge d_S^*d_Su(x_0)   
     \ge a - e^{u(x_0)}h(x_0)
     \ge a - e^{u(x_0)}\left\|h\right\|_{L^\infty}
$$
and hence
$$
    \inf_Su = u(x_0) 
    \ge \log\left(\frac{a}{\left\|h\right\|_{L^\infty}}\right).
$$
This proves the lemma in the smooth case. 

Now suppose that $h\in C^0(S)$ and $u\in W^{2,p}(S)$ satisfy
the hypotheses of the lemma.   Then $u$ is continuous
and~(\ref{eq:KW0}) shows that $d_S^*d_Su\in L^\infty(S)$.  
Choose sequences $a_\nu\to a$ and $A_\nu\to A$ such that 
$$
      0 < a_\nu < a \le A < A_\nu.
$$
Then there exist sequences of smooth functions 
$u_\nu,h_\nu\in\Cinf(S)$ such that $h_\nu$ converges
uniformly to $h$, $u_\nu$ converges to $u$ in the 
$W^{2,p}$-norm, $h_\nu\ge 0$, and 
$$
     a_\nu-e^{u_\nu}h_\nu
     \le d_S^*d_Su_\nu 
     \le A_\nu - f^{-1}(e^{u_\nu}h_\nu).
$$
To see this, we may first choose a sequence $w_\nu\in\Cinf(S)$
converging to $d_S^*d_Su$ in the $L^p$ norm and satisfying
$a_\nu-e^uh < w_\nu < A_\nu - f^{-1}(e^uh)$.  Then define
$u_\nu$ as the solution of the equation $d_S^*d_Su_\nu=w_\nu$ 
with $\int_S(u_\nu-u)\,\dvol_S=0$ and choose any sequence
$h_\nu\in\Cinf(S)$ converging uniformly to $h$ to obtain
the required estimate for $u_\nu$ and $h_\nu$.  
It then follows that $u_\nu$ and $h_\nu$ satisfy the hypotheses 
of the lemma with $a$ and $A$ replaced by $a_\nu$ and $A_\nu$, 
respectively. Hence they satisfy the conclusion and so the 
required estimate for $u$ and $h$ follows by taking 
the limit $\nu\to\infty$.
\end{proof}

\begin{proposition}\label{prop:KWS}
Let $p>1$. For every $t\in\R$ and every $h\in C^0(S)$ 
such that $h\ge 0$ there exists a unique solution $u\in W^{2,p}(S)$
of the equation
\begin{equation}\label{eq:KWS}
     d_S^*d_Su + e^uh = \frac{1}{\Vol(S)}\int_Se^uh\,\dvol_S,\qquad
     \frac{1}{\Vol(S)}\int_Su\,\dvol_S = t.
\end{equation}
Moreover, if $h\in W^{k,p}(S)$ for some integer $k\ge 1$ then 
$u\in W^{k+2,p}(S)$. If $kp>2$ then the map $(h,t)\mapsto u$ 
which assigns to each pair $(h,t)\in W^{k,p}(S)\times\R$ 
that satisfies $h\ge 0$ the unique solution 
$u\in W^{k+2,p}(S)$ of~(\ref{eq:KWS}) extends 
to a smooth map between open subsets of Banach spaces.
\end{proposition}

\begin{proof}
The proof has three steps.

\medskip
\noindent{\bf Step~1.}
{\it For every $p>1$ and every $c>0$ there exists a constant 
$c_p>0$ such that, if $h\in C^0(S)$ and $t\in\R$ 
satisfy
\begin{equation}\label{eq:htc}
      h\ge 0,\qquad 
      h_0 := \frac{1}{\Vol(S)}\int_Sh\,\dvol_S > \frac{1}{c},\qquad
      \left\|h\right\|_{L^\infty}\le c,\qquad
     |t|\le c,
\end{equation}
then
$$
     \left\|u\right\|_{W^{2,p}} \le c_p
$$
for every solution $u\in W^{2,p}(S)$ of~(\ref{eq:KWS}).}

\medskip
\noindent
Let $u$ be a solution of~(\ref{eq:KWS}) and denote
$$
     a := \frac{1}{\Vol(S)}\int_Se^uh\,\dvol_S.
$$
Then, by~(\ref{eq:htc}) and Lemma~\ref{le:KW0} with $C=0$ and $A=a$, 
we have
$$
     \log\left(\frac{a}{c}\right)
     \le u \le 
     \log\left(\frac{a}{h_0}\right) 
     + \frac{cc_Sa}{h_0}.
$$
Integrating the first inequality over $S$ gives 
$$
     \log\left(\frac{a}{c}\right) 
     \le \frac{1}{\Vol(S)}\int_Su\,\dvol_S
     =   t
     \le c
$$
and hence 
$
     a\le ce^c.
$
Moreover,
$$
     e^u \le \frac{a}{h_0}e^{cc_Sa/h_0}
     \le ac e^{c^2c_Sa}
     \le c^2e^ce^{c^3c_Se^c}.
$$
Hence $e^u$ satisfies a uniform upper bound, 
depending only on $S$ and $c$.  Hence there exists a 
constant $c'=c'(c)>0$ such that $\|d_S^*d_Su\|_{L^\infty}\le c'$
for every solution of~(\ref{eq:KWS}). 
Since $\Vol(S)^{-1}\int_Su\,\dvol_S=t\in[-c,c]$, 
Step~1 follows from elliptic regularity for the 
Laplace operator on $S$. 

\medskip
\noindent{\bf Step~2.}
{\it Consider the Banach spaces 
$$
     \Xx := W^{2,p}(S),\qquad \Yy := L^p_0(S)\times\R,
$$
where $L^p_0(S)$ denotes the space of $L^p$-functions 
on $S$ with mean value zero. For $h\in C^0(S)$ 
define $\Ff_h:\Xx\to\Yy$ by 
$$
     \Ff_h(u) 
     := \left(d_S^*d_Su + e^uh - \frac{1}{\Vol(S)}\int_Se^uh\,\dvol_S,
     \frac{1}{\Vol(S)}\int_Su\,\dvol_S\right).
$$
If $h\ge0$ then the differential $d\Ff_h(u):\Xx\to\Yy$ 
is a Banach space isomorphism for every $u\in\Xx$.}

\medskip
\noindent
The differential of $\Ff_h$ is given by 
$$
     d\Ff_h(u)\xi 
     = \left(d_S^*d_S\xi + e^uh\xi - \frac{1}{\Vol(S)}\int_Se^uh\xi\,\dvol_S,
       \frac{1}{\Vol(S)}\int_S\xi\,\dvol_S\right).
$$
Hence $d\Ff_h(u):\Xx\to\Yy$ is a Fredholm operator 
of index zero. Multiplying the first component of $d\Ff_h(u)\xi$ 
by $\xi$ and integrating over $S$ we find that the kernel of 
$d\Ff_h(u)$ consists of all functions $\xi\in W^{2,p}(S)$ that
satisfy 
$$
     \int_S\left|d_S\xi\right|^2\,\dvol_S 
     + \int_S e^uh\left|\xi\right|^2\,\dvol_S = 0,\qquad
     \int_S\xi\,\dvol_S = 0.
$$
Hence $d\Ff_h(u)$ is bijective whenever $h\ge 0$. 

\medskip
\noindent{\bf Step~3.}
{\it We prove the proposition.}

\medskip
\noindent
If $h=0$ then every solution of~(\ref{eq:KWS}) is constant
and hence $u\equiv t$ is the only solution. Now assume $h=1$
and let $u\in W^{2,p}(S)$ be a solution of~(\ref{eq:KWS}). 
Then 
$$
     d_S^*d_Su+e^u=\frac{1}{\Vol(S)}\int_Se^u\,\dvol_S
$$
and hence, by Lemma~\ref{le:KW0}, 
$$
     e^u\ge \frac{1}{\Vol(S)}\int_Se^u\,\dvol_S.
$$
This implies $e^u\equiv\mbox{ constant}$ and hence 
$u\equiv t$. Thus we have proved the existence and uniqueness 
statament in the cases $h=0$ and $h=1$.  
Now let $h\in C^0(S)$ be any nonnegative function 
such that $\int_Sh\,\dvol_S>0$, and define 
$$
     h_\eps := (1-\eps)h + \eps.
$$
We prove that the number of solutions of~(\ref{eq:KWS})
with $h$ replaced by $h_\eps$ is independent of $\eps$.
To see this consider the set 
$$
     \Mm := \left\{(\eps,u)\,|\,0\le\eps\le1,\,u\in W^{2,p}(S),\,
     \Ff_{h_\eps}(u)=(0,t)\right\}.
$$
By Step~2, this set is a smooth $1$-manifold with boundary
and the projection $\Mm\to[0,1]:(\eps,u)\mapsto\eps$ is a 
submersion.  That $\Mm$ is compact follows from the fact that, 
by Step~1, there exists a constant $c_p>0$ such that 
$$
     (\eps,u)\in\Mm\qquad\IMP\qquad \left\|u\right\|_{W^{2,p}}\le c_p. 
$$
Hence every sequence $(\eps_i,u_i)\in\Mm$ has a subsequence such 
that $u_i$ converges in $C^0(S)$ and $\eps_i$ converges. 
Hence, for this subsequence, $e^{u_i}h_{\eps_i}$ converges 
in $C^0(S)$ and so, by elliptic regularity for the Laplace
operator on $S$, $u_i$ converges in $W^{2,p}(S)$ (for any $p>1$). 
Thus $\Mm$ is compact and so the number $\#\Ff_{h_\eps}^{-1}(0,t)$
is independent of $\eps\in[0,1]$.  For $\eps=1$ this number is one
and this proves the existence and uniqueness statement.  
That $h\in W^{k,p}$ implies $u\in W^{k+2,p}$ follows from
elliptic regularity for the Laplace operator.
That the map $(h,t)\mapsto u$ is smooth follows 
from the implicit function theorem and Step~2. 
\end{proof}

\begin{proposition}\label{prop:KW}
Let $a>0$ and $\Hh\subset C^0(\Sigma\times S)$ be a compact 
set such that $h \ge 0$ for every $h\in\Hh$. 
Then there exists a constant 
$\delta=\delta(\Sigma,S,\Hh,a)>0$ such that 
the following holds. If $p>1$ and $u\in C^0(\Sigma,W^{2,p}(S))$ 
satisfies $u_\Sigma\in W^{2,p}(\Sigma)$ and
\begin{equation}\label{eq:KWa}
    d_\Sigma^*d_\Sigma u_\Sigma + d_S^*d_Su
    + e^uh = a
\end{equation}
for some $h\in\Hh$ then $d_\Sigma^*d_\Sigma u_\Sigma$ 
is continuous and
$$
    \delta h_\Sigma\le a-d_\Sigma^*d_\Sigma u_\Sigma
    \le \delta^{-1}\sup_Sh,\qquad
    \delta\le e^u\le\delta^{-1}.
$$
\end{proposition}

\begin{proof}
The proof has three steps.

\medskip
\noindent{\bf Step~1.}
{\it Let $c_S$ be the constant of Lemma~\ref{le:KW0}
and choose $c>0$ such that 
$$
     h\in\Hh\qquad\IMP\qquad \left\|h\right\|_{L^\infty}\le c.
$$
Define $f:[0,\infty)\to[0,\infty)$ by
$
     f(r) := r e^{cc_Sr}.
$
Then
$$
     f^{-1}(e^u)h_\Sigma
     \le a - d_\Sigma^*d_\Sigma u_\Sigma
     \le e^u\sup_Sh
$$
for every $h\in\Hh$ and every solution $u$ of~(\ref{eq:KWa}).}

\medskip
\noindent
Integrating~(\ref{eq:KWa}) over $S$ we obtain
$$
     a - d_\Sigma^*d_\Sigma u_\Sigma 
     = \frac{1}{\Vol(S)}\int_S e^uh\,\dvol_S
     \ge 0.
$$
In particular, $d_\Sigma^*d_\Sigma u_\Sigma$ is continuous.
If $h|_{\{z\}\times S}\equiv 0$ then
$a=d_\Sigma^*d_\Sigma u_\Sigma(z)$ and $h_\Sigma(z)=0$.
So the assertion of Step~1 holds trivially at the point $z$.
If $h|_{\{z\}\times S}\not\equiv 0$ then
$a-d_\Sigma^*d_\Sigma u_\Sigma(z)>0$ and so the restrictions 
of $h$ and $u$ to $\{z\}\times S$ satisfy the requirements
of Lemma~\ref{le:KW0} with $C=0$ and $a=A$ replaced by 
the constant $a-d_\Sigma^*d_\Sigma u_\Sigma(z)$. Hence
$$
    \log\left(\frac{a-d_\Sigma^*d_\Sigma u_\Sigma}
    {\sup_Sh}\right)
    \le u \le 
    \log\left(\frac{a-d_\Sigma^*d_\Sigma u_\Sigma}
    {h_\Sigma}\right) 
    + \frac{cc_S\left(a-d_\Sigma^*d_\Sigma u_\Sigma\right)}
    {h_\Sigma}.
$$
This implies the assertion of Step~1 in the case $h_\Sigma(z)\ne 0$.

\medskip
\noindent{\bf Step~2.}
{\it There exists a constant 
$\delta=\delta(\Sigma,S,\Hh,a)>0$ 
such that 
$$
    \delta\le e^{u_\Sigma}\le \delta^{-1}
$$
for every $h\in\Hh$ and every solution $u$ of~(\ref{eq:KWa}).}

\medskip
\noindent
By the proof of Step~1, we have
$$
    \log\left(\frac{a-d_\Sigma^*d_\Sigma u_\Sigma}{c}\right)
    \le u_\Sigma
    \le 
    \log\left(\frac{a-d_\Sigma^*d_\Sigma u_\Sigma}
    {h_\Sigma}\right)
    + \frac{cc_S\left(a-d_\Sigma^*d_\Sigma u_\Sigma\right)}
    {h_\Sigma}
$$
whenever $h_\Sigma(z)>0$ and hence
$$
    a - ce^{u_\Sigma}
    \le d_\Sigma^*d_\Sigma u_\Sigma
    \le a - f^{-1}(e^{u_\Sigma})h_\Sigma.
$$
Hence $u_\Sigma\in W^{2,p}(\Sigma)$ satisfies the second 
inequality in~(\ref{eq:KW0}) with $A=a$, $C=cc_S$, and $h$ replaced 
by $h_\Sigma$. It satisfies the first inequality with $h$ replaced 
by $c$.  Hence, by Lemma~\ref{le:KW0},
$$
    \log\left(\frac{a}{c}\right)
    \le u_\Sigma \le
    \log\left(\frac{a\Vol(\Sigma)}
    {\int_\Sigma h_\Sigma\,\dvol_\Sigma}\right) 
    + \frac{c\left(c_S+c_\Sigma\right)a\Vol(\Sigma)}
      {\int_\Sigma h_\Sigma\,\dvol_\Sigma}.
$$
This proves Step~2.

\medskip
\noindent{\bf Step~3.}
{\it We prove the proposition.}
     
\medskip
\noindent
For $h\in\Hh$ and $z\in\Sigma$ denote by $h_z:S\to\R$ the function 
$h_z(x):=h(z,x)$ and let $\Tt_{h_z}:\R\to W^{2,p}(S)$ be the 
map which assigns to every $t\in\R$ the unique solution 
$u=\Tt_{h_z}(t)\in W^{2,p}(S)$ of~(\ref{eq:KWS}) 
with $h$ replaced by $h_z$.  Then every solution 
$u:\Sigma\times S\to\R$ of~(\ref{eq:KWa}) satisfies
$$
     u(z,\cdot) = \Tt_{h_z}(u_\Sigma(z)).
$$
By Proposition~\ref{prop:KWS}, the map
$
     \Hh\times\Sigma\times\R\to W^{2,p}(S):
     (h,z,t)\mapsto\Tt_{h_z}(t)
$
is continuous. Since $\Hh$ is compact it follows that there exists
an $\eps>0$ such that 
$$
     h\in\Hh,\quad z\in\Sigma,\quad \delta\le e^t\le\delta^{-1}
     \qquad\IMP\qquad
     \left\|\Tt_{h_z}(t)\right\|_{L^\infty(S)}\le |\log(\eps)|.
$$
This implies $\eps\le e^u\le\eps^{-1}$ for every $h\in\Hh$
and every solution $u$ of~(\ref{eq:KWa}). The inequality
for $a-d_\Sigma^*d_\Sigma u_\Sigma$ now follows from Step~1. 
\end{proof}

\begin{proof}[Proof of Theorem~\ref{thm:KW}]
The proof has four steps.

\medskip
\noindent{\bf Step~1.}
{\it It suffices to prove the theorem if $f$ is constant.}

\medskip
\noindent
Let $v_\Sigma\in W^{2,p}(\Sigma)$ be a
solution of the equation
$$
     d_\Sigma^*d_\Sigma v_\Sigma = f_\Sigma-a,\qquad
     a := \frac{1}{\Vol(\Sigma)}\int_\Sigma f_\Sigma\,\dvol_\Sigma,
$$
and let $v\in C^0(\Sigma,W^{2,p}(S))$ be the unique solution
of the equation
$$
     d_S^*d_Sv = f-f_\Sigma,\qquad
     \frac{1}{\Vol(S)}\int_Sv\,\dvol_S = v_\Sigma.
$$
This equation is understood pointwise for $z\in\Sigma$. 
Then $v$ is continuous and
$$
     d_\Sigma^*d_\Sigma v_\Sigma + d_S^*d_Sv = f-a.
$$
Note that if $f$ is smooth then so is $v$.
Moreover, $u$ is a solution of~(\ref{eq:KWa}) with $h$ 
replaced by $e^vh$ if and only if $u+v$ is a solution 
of~(\ref{eq:KW}).  

\medskip
\noindent
{\bf Step~2.}
{\it Let $h\in C^0(S)$ such that $h\ge 0$ and define 
$f_h:\R\to\R$ by 
$$
     f_h(t) := \frac{1}{\Vol(S)}\int_Se^uh\,\dvol_S,
$$
where $u\in W^{2,p}(S)$ is the unique solution of~(\ref{eq:KWS}). 
Then 
$
     f_h'(t)\ge 0
$
for every $t\in\R$ with equality if and only if $h\equiv0$.}

\smallbreak

\medskip
\noindent
Let $u\in W^{2,p}(S)$ be the unique solution of~(\ref{eq:KWS})
and $\xi\in W^{2,p}(S)$ be the unique solution of the equation
$$
     d_S^*d_S\xi + e^uh\xi 
     = \frac{1}{\Vol(S)}\int_Se^uh\xi\,\dvol_S,\qquad
     \frac{1}{\Vol(S)}\int_S\xi\,\dvol_S = 1.
$$
Then
$$
     f_h'(t)
     = \frac{1}{\Vol(S)}\int_Se^uh\xi\,\dvol_S
     = \frac{1}{\Vol(S)}\int_S\left(
       \left|d_S\xi\right|^2 + e^uh\left|\xi\right|^2
       \right)\,\dvol_S
     \ge 0.
$$
Equality implies that $\xi\equiv1$ and $h\equiv0$.

\medskip
\noindent{\bf Step~3.}
{\it Given a nonzero continuous function 
$h:\Sigma\times S\to[0,\infty)$ define 
$\Ff_h:W^{2,p}(\Sigma)\to L^p(\Sigma)$ by 
$$
     \Ff_h(u_\Sigma)(z) 
     := d_\Sigma^*d_\Sigma u_\Sigma(z)
        + f_{h_z}(u_\Sigma(z))
$$
for $z\in\Sigma$, where $h_z:=h(z,\cdot)\in C^0(S)$.
Then $d\Ff_h(u_\Sigma):W^{2,p}(\Sigma)\to L^p(\Sigma)$
is a Banach space isomorphism for every 
$u_\Sigma\in W^{2,p}(\Sigma)$.}

\medskip
\noindent
This follows directly from Step~2. 

\medskip
\noindent{\bf Step~4.}
{\it We prove the theorem.}

\medskip
\noindent
By Step~1 we may assume $f\equiv a$. 
Assume first that $h\equiv1$.  We claim that in
this case $u\equiv\log(a)$ is the only solution 
of~(\ref{eq:KWa}). To see this, note that, by 
Proposition~\ref{prop:KWS}, the restriction of $u$ 
to each fibre $\{z\}\times S$ is constant, hence 
$u=u_\Sigma$ and $d_\Sigma^*d_\Sigma u_\Sigma + e^u=a$,
and hence, again by Proposition~\ref{prop:KWS}, 
$u=u_\Sigma=\log(a)$.  Now let $h\in C^0(\Sigma\times S)$
be any nonzero nonnegative function. Note that $\Ff_h(u_\Sigma)=a$ iff
$u$ is a solution of (\ref{eq:KWa}). Define 
$h_\eps\in C^0(\Sigma\times S)$ by 
$
     h_\eps := (1-\eps)h + \eps
$
and consider the set $\Mm\subset[0,1]\times W^{2,p}(\Sigma)$
given by 
$$
     \Mm := \left\{(\eps,u_\Sigma)\,|\,0\le\eps\le1,\,
     u_\Sigma\in W^{2,p}(\Sigma),\,\Ff_{h_\eps}(u_\Sigma)=a\right\}.
$$
By Step~3, this is a $1$-manifold with boundary and the 
projection $\Mm\to[0,1]:(\eps,u_\Sigma)\mapsto\eps$ is a 
submersion.  To prove that $\Mm$ is compact note that, 
by Proposition~\ref{prop:KW}, there exists a constant $c>0$ 
such that 
$$
     \left\|u_\Sigma\right\|_{W^{2,p}}\le c
$$
for every
$
     (\eps,u_\Sigma)\in\Mm.
$
Hence every sequence $(\eps_i,u_i)\in\Mm$ has a subsequence
such that $\eps_i$ converges and $u_i$ converges in $C^0(\Sigma)$.
The equation 
$$
     d_\Sigma^*d_\Sigma u_i(z)
        + f_{h_{\eps_i,z}}(u_i(z)) = a
$$
now shows that $u_i$ converges in $W^{2,p}(\Sigma)$.
Hence $\Mm$ is compact and so the number $\#\Ff_{h_\eps}^{-1}(a)$
is independent of $\eps$.  For $\eps=1$ we have seen that
this number is one.  This proves the existence and uniqueness 
statement.  

Now suppose that $h$ is smooth. Then the function 
$\Sigma\times\R\to\R:(z,t)\mapsto f_{h_z}(t)$
is smooth and hence, by a standard elliptic
bootstrapping argument, the unique solution 
$u_\Sigma:\Sigma\to\R$ of the equation
$$
     d_\Sigma^*d_\Sigma u_\Sigma(z) + f_{h_z}(u_\Sigma(z)) = a
$$
is smooth. Hence, by Lemma~\ref{le:KW0}, the unique solution
$u:\Sigma\times S\to\R$ of the equation
$$
     d_S^*d_Su + e^uh = \frac{1}{\Vol(S)}\int_Se^uh\,\dvol_S,\qquad
     \frac{1}{\Vol(S)}\int_Su\,\dvol_S = u_\Sigma,
$$
is smooth.  This proves the theorem.
\end{proof}


\section{The local slice theorem}\label{app:slice}

Let $G$ be a compact Lie group and 
$P\to X$ be a principal $G$-bundle over a
compact $n$-manifold $X$. For $p>n/2$ denote by 
$
     \Gg^{k+1,p} = \Gg^{k+1,p}(P)
$
the space of all $W^{k+1,p}$-sections of the bundle
$P\times_{\ad}\G\to X$. Fix a smooth reference 
connection $\hat A\in\Aa(P)$ and denote by 
$$
     \Aa^{1,p}(P)
     := \left\{\hat A+\alpha\,|\,\alpha\in
     W^{1,p}(X,T^*X\otimes\g_P)\right\}
$$
the space of $W^{1,p}$-connections. 
This space is independent of the connection~$\hat A$. 

\begin{theorem}\label{thm:slice}
Let $p,q$ be positive real number such that 
\begin{equation}\label{eq:pq}
     q\ge p>\frac{n}{2},\qquad q>n,\qquad
     \mbox{if }p\le n\mbox{ then }
     q < \frac{np}{n-p}.
\end{equation}
Then, for every $A_0\in\Aa^{1,p}(P)$
and every positive constant $c_0$, 
there exist positive constants $c$ and $\delta$ 
such that the following holds. 
If $A\in\Aa^{1,p}(P)$ satisfies
$$
     \left\|A-A_0\right\|_{W^{1,p}} \le c_0,\qquad
     \left\|A-A_0\right\|_{L^q} \le \delta
$$
then there exists a gauge transformation 
$g\in\Gg^{2,p}(P)$ such that 
$$
      d_{A_0}^*(g^*A-A_0) = 0
$$
and
$$
      \left\|g^*A-A\right\|_{L^q}
      \le c\left\|A-A_0\right\|_{L^q},
      \qquad
      \left\|g^*A-A\right\|_{W^{1,p}}
      \le c\left\|A-A_0\right\|_{W^{1,p}}.
$$
\end{theorem}

\begin{lemma}\label{le:product}
Let $p,q,r$ be positive real numbers such that
$$
     r\le p,\qquad r\le q,\qquad
     \frac{1}{p}+\frac{1}{q} < \frac{1}{n}+\frac{1}{r}.
$$
Then there exists a constant $c>0$ such that
$$
     \left\|fg\right\|_{W^{1,r}}
     \le c\left\|f\right\|_{W^{1,p}}\left\|g\right\|_{W^{1,q}}.
$$
for $f,g\in\Cinf_0(\R^n)$.  
In particular, this holds when $p$ and $q$ 
satisfy~(\ref{eq:pq}) and $r=p$. 
It also holds when $p=q=r>n$. 
\end{lemma}

\begin{proof}
By H\"older's inequality and the product rule, we have
$$
     \left\|fg\right\|_{W^{1,r}}
     \le \left\|f\right\|_{W^{1,p}}\left\|g\right\|_{L^{rp/(p-r)}}
     + \left\|f\right\|_{L^{rq/(q-r)}}\left\|g\right\|_{W^{1,q}}.
$$
If $q>n$ the Sobolev embedding theorem asserts
that the $L^{rp/(p-r)}$-norm of $g$ can be estimated
from above by the $W^{1,q}$-norm.  
The same holds for $q=n$ since then 
it follows from the hypotheses that $p>r$.
If $q<n$ we have $r<p$ and 
$$
     \frac{rp}{p-r} 
     = \left(\frac{1}{r}-\frac{1}{p}\right)^{-1}
     < \left(\frac{1}{q}-\frac{1}{n}\right)^{-1}
     = \frac{nq}{n-q}
$$
and hence the $L^{rp/(p-r)}$-norm of $g$ can 
again be estimated from above by the $W^{1,q}$-norm. 
Similarly, the $L^{rq/(q-r)}$-norm of $f$ can 
again be estimated from above by the $W^{1,p}$-norm. 
\end{proof}

\begin{lemma}\label{le:laplace1}
If $A\in\Aa^{1,p}(P)$ and $p>n/2$ then the following holds
for every $r>1$.
\begin{description}
\item[(i)]
If $p<n$ assume in addition $r<np/(n-p)$. Then the operator
$d_A:\W^{1,r}(X,\g_P)\to L^r(X,T^*X\otimes\g_P)$
is a compact perturbation of $d_{\hat A}$.
Similarly for $d_A^*$.
\item[(ii)]
For $r\le p$ the operator
$d_A:\W^{2,r}(X,\g_P)\to W^{1,r}(X,T^*X\otimes\g_P)$
is a compact perturbation of $d_{\hat A}$.
Similarly for $d_A^*$.
\item[(iii)]
For $r\le p$ the operator 
$d_A^*d_A:\W^{2,r}(X,\g_P)\to L^r(X,\g_P)$
is a compact perturbation of 
$d_{\hat A}^*d_{\hat A}$.
\end{description}
\end{lemma}

\begin{proof}
For $\xi\in\Om^0(X,\g_P)$ and $\alpha\in\Om^1(X,\g_P)$
we have
$$
     d_A\xi - d_{\hat A}\xi = [(A-\hat A),\xi],\qquad
     d_A^*\alpha - d_{\hat A}^*\alpha 
     =  *[*(A-\hat A)\wedge \alpha].
$$
Assume first that $p<n$. Then $r<np/(n-p)$
and hence there exists a real number $s>1$ such that
$
     1/s + (n-p)/np = 1/r.
$
Since $2p<n$ it follows that 
$
     s < nr/(n-r)
$
whenever $r<n$. Hence the Sobolev embedding theorem 
asserts that the inclusion 
$W^{1,r}(X,\g_P)\INTO L^s(X,\g_P)$ is compact.
It also asserts that 
$A-\hat A\in L^{np/(n-p)}(X,T^*X\otimes\g_P)$
and hence, by H\"older's inequality,
the operator 
$
     L^s(X,\g_P)\to L^r(X,T^*X\otimes\g_P):
     \xi\mapsto [(A-\hat A),\xi]
$
is bounded. Hence the composition with the inclusion
$W^{1,r}\INTO L^s$ is compact.
If $p\ge n$ choose any number $s>r$ such that
the inclusion $W^{1,r}(X,\g_P)\INTO L^s(X,\g_P)$
is compact and use the fact that
$A-\hat A\in L^{rs/(s-r)}(X,T^*X\otimes\g_P)$.
This proves~(i). 

We prove~(ii).
By Lemma~\ref{le:product} the operator
$$
     W^{1,s}(X,\g_P)\to W^{1,r}(X,T^*X\otimes\g_P):
     \xi\mapsto[(A-\hat A),\xi]
$$
is bounded whenever 
$
     r\le p,
$
$
     r\le s,
$
and
$
     1/p + 1/s < 1/n + 1/r.
$
If $r>n$ then $p>n$ and we may choose $s=r$. 
If $r\le n$ then, since $2p>n$, we have 
$$
     \frac{1}{p}+\frac{n-r}{nr} < \frac{1}{n}+\frac{1}{r}
$$
and hence may choose $s$ such that $r\le s<nr/(n-r)$.
In either case the Sobolev embedding theorem 
asserts that the inclusion 
$W^{2,r}(X,\g_P)\INTO W^{1,s}(X,\g_P)$
is compact.  This proves~(ii).
Assertion~(iii) follows directly from~(i) and~(ii). 
\end{proof}

\begin{lemma}\label{le:laplace2}
Suppose $p$ and $q$ satisfy~(\ref{eq:pq})
and let $A\in\Aa^{1,p}(P)$.  Then there 
exists a constant $c=c(A)=c(A,p,q)>0$ such that,
for every $\alpha\in W^{1,p}(X,T^*X\otimes\g_P)$,
there exists a $\xi\in W^{2,p}(X,\g_P)$
such that 
\begin{equation}\label{eq:laplace}
     d_A^*d_A\xi = d_A^*\alpha
\end{equation}
and 
\begin{equation}\label{eq:estimate}
     \left\|\xi\right\|_{W^{2,p}}
     \le c\left\|d_A^*\alpha\right\|_{L^p},\qquad
     \left\|\xi\right\|_{W^{1,q}}
     \le c\left\|\alpha\right\|_{L^q}.
\end{equation}
\end{lemma}

\begin{proof}
Let $r:=q/(q-1)$ so that
$$
     \frac{1}{q}+\frac{1}{r} = 1.
$$
By Lemma~\ref{le:laplace1}, the operator
$
     d_A:W^{1,s}(X,\g_P)\to L^s(X,T^*X\otimes\g_P)
$
is a compact perturbation of $d_{\hat A}$
for $s=q$ and $s=r$. Let 
$$
     W^{-1,q}(X,\g_P) := (W^{1,r}(X,\g_P))^*
$$
and denote by 
$d_A^*:L^q(X,T^*X\otimes\g_P)\to W^{-1,q}(X,\g_P)$
the dual operator of 
$d_A:W^{1,r}(X,\g_P)\to L^r(X,T^*X\otimes\g_P)$.
Then 
\begin{equation}\label{eq:laplace1}
     d_A^*d_A:W^{1,q}(X,\g_P)\to W^{-1,q}(X,\g_P)
\end{equation}
is a compact perturbation of $d_{\hat A}^*d_{\hat A}$
and hence is a Fredholm operator of index zero. 
Likewise, it follows from Lemma~\ref{le:laplace1} 
that the operator
\begin{equation}\label{eq:laplace2}
     d_A^*d_A:W^{2,p}(X,\g_P)\to L^p(X,\g_P)
\end{equation}
is a compact perturbation of $d_{\hat A}^*d_{\hat A}$
and hence is also a Fredholm operator of index zero. 
The operator~(\ref{eq:laplace1}) is a natural extension
of~(\ref{eq:laplace2}).  
Taking the $L^2$-inner product of $d_A^*d_A\xi$
with $\xi$ for $\xi\in W^{2,p}(X,\g_P)$ 
we see that the kernel of~(\ref{eq:laplace2}) 
is the finite dimensional subspace
$$
     H^0(X,A):=\ker\,d_A\subset W^{2,p}(X,\g_P).
$$
The operator~(\ref{eq:laplace1}) has the same 
kernel, because every $\xi\in W^{1,q}(X,\g_P)$
with $d_A\xi=0$ lies in $W^{2,p}(X,\g_P)$. 
Choose a complement $E'$ of $H^0(X,A)$
in the Sobolev space $W^{1,q}(X,\g_P)$.  Then 
$$
     E := E'\cap W^{2,p}(X,\g_P)
$$
is a complement of $H^0(X,A)$
in $W^{2,p}(X,\g_P)$.  
Let $F'$ denote the image of the operator 
$d_A^*:L^q(X,T^*X\otimes\g_P)\to W^{-1,q}(X,\g_P)$
and $F$ denote the image of the operator
$d_A^*:W^{1,p}(X,T^*X\otimes\g_P)\to L^p(X,\g_P)$.
Then $H^0(X,A)\subset W^{1,r}(X,\g_P)$ annihilates
$F'$ and is $L^2$-orthogonal to $F$.
Moreover $F'$ contains the image of~(\ref{eq:laplace1})
and $F$ contains the image of~(\ref{eq:laplace2}).
Hence, for dimensional reasons, 
$F'$ is equal to the image of~(\ref{eq:laplace1})
and $F$ is equal to the image of~(\ref{eq:laplace2}).
Thus $d_A^*d_A$ is a Banach space isomorphism 
from $E$ to $F$ and extends to a Banach space
isomorphism from $E'$ to $F'$

Now let $\alpha\in W^{1,p}(X,T^*X\otimes\g_P)$.
Then $d_A^*\alpha\in F$ and hence there exists 
a unique $\xi\in E$ that satisfies~(\ref{eq:laplace}).
By the open mapping theorem, 
this solution of~(\ref{eq:laplace}) satisfies
$$
     \left\|\xi\right\|_{W^{2,p}}
     \le \left\|(d_A^*d_A)^{-1}\right\|_{\Ll(F,E)}
     \left\|d_A^*\alpha\right\|_{L^p}.
$$
Since $\xi\in E'$ it also satisfies
$$
     \left\|\xi\right\|_{W^{1,q}}
     \le \left\|(d_A^*d_A)^{-1}\right\|_{\Ll(F',E')}
     \left\|d_A^*\alpha\right\|_{W^{-1,q}}.
$$
Now
\begin{eqnarray*}
     \left\|d_A^*\alpha\right\|_{W^{-1,q}}
&= &
     \sup_{\eta\ne 0}
     \frac{\inner{\eta}{d_A^*\alpha}_{W^{1,r},W^{-1,q}}}
          {\left\|\eta\right\|_{W^{1,r}}}  \\
&= &
     \sup_{\eta\ne 0}
     \frac{\inner{d_A\eta}{\alpha}_{L^r,L^q}}
          {\left\|\eta\right\|_{W^{1,r}}}  \\
&\le &
     \sup_{\eta\ne 0}
     \frac{\left\|d_A\eta\right\|_{L^r}
           \left\|\alpha\right\|_{L^q}}
          {\left\|\eta\right\|_{W^{1,r}}}  \\
&\le &
     c\left\|\alpha\right\|_{L^q},
\end{eqnarray*}
where the constant $c$ depends only on $A$ and $r$. 
Hence $\xi$ satisfies~(\ref{eq:estimate}). 
\end{proof}

\begin{lemma}\label{le:quadratic}
Suppose $p$ and $q$ satisfy~(\ref{eq:pq})
and fix a constant $c_0>0$. 
Then there exists a constant $c=c(c_0,p,q)>0$ such that
the following holds.  If $\xi\in W^{2,p}(X,\g_P)$ 
satisfies
$
    \left\|\xi\right\|_{W^{2,p}}
    \le c_0
$
then, for every $A\in\Aa^{1,p}(P)$, we have
$$
     \left\|\exp(\xi)^*A-A-d_A\xi\right\|_{W^{1,p}}
     \le c\left(1+\left\|A-\hat A\right\|_{W^{1,p}}\right)
         \left\|\xi\right\|_{W^{1,q}}
         \left\|\xi\right\|_{W^{2,p}},
$$
$$
     \left\|\exp(\xi)^*A-A\right\|_{W^{1,p}}
     \le c\left(1+\left\|A-\hat A\right\|_{W^{1,p}}\right)
     \left\|\xi\right\|_{W^{2,p}},
$$
$$
     \left\|\exp(\xi)^*A-A\right\|_{L^q}
     \le c\left(1+\left\|A-\hat A\right\|_{W^{1,p}}\right)
     \left\|\xi\right\|_{W^{1,q}}.
$$
\end{lemma}

\begin{proof}   
The function 
$
     \alpha(t):=\exp(t\xi)^*A-A
$
satisfies the differential equation
$\dot\alpha(t) = d_A\xi - [\xi,\alpha(t)]$
and $\alpha(0) = 0$.
Hence
$$
     \alpha(t)
     = \sum_{k=0}^\infty
     \frac{(-1)^kt^{k+1}}{(k+1)!}\ad(\xi)^kd_A\xi,
$$
and hence 
\begin{equation}\label{eq:exp}
     \exp(\xi)^*A-A-d_A\xi 
     = \sum_{k=1}^\infty
       \frac{(-1)^k}{(k+1)!}\ad(\xi)^kd_A\xi.
\end{equation}
Now
\begin{eqnarray*}
     \left\|d_A\xi\right\|_{L^q}
&\le &
     \left\|d_{\hat A}\xi\right\|_{L^q} 
     + \left\|[(A-\hat A),\xi]\right\|_{L^q} \\
&\le &
     \left\|d_{\hat A}\xi\right\|_{L^q}
     + \left\|A-\hat A\right\|_{L^q} 
       \left\|\xi\right\|_{L^\infty}  \\
&\le &
     c\left(1+\left\|A-\hat A\right\|_{W^{1,p}}\right)
     \left\|\xi\right\|_{W^{1,q}} 
\end{eqnarray*}
and, by Lemma~\ref{le:product},
\begin{eqnarray*}
     \left\|d_A\xi\right\|_{W^{1,p}}
&\le &
     \left\|d_{\hat A}\xi\right\|_{W^{1,p}} 
     + \left\|[(A-\hat A),\xi]\right\|_{W^{1,p}} \\
&\le &
     \left\|d_{\hat A}\xi\right\|_{W^{1,p}}
     + c'\left\|A-\hat A\right\|_{W^{1,p}} 
        \left\|\xi\right\|_{W^{1,q}}  \\
&\le &
     c''\left(1+\left\|A-\hat A\right\|_{W^{1,p}}\right)
     \left\|\xi\right\|_{W^{2,p}}.
\end{eqnarray*}
Hence the assertion follows from~(\ref{eq:exp})
and Lemma~\ref{le:product}.
\end{proof}

\begin{proof}[Proof of Theorem~\ref{thm:slice}]
The proof is by Newton's iteration.   

\medskip
\noindent{\bf Step~1}
{\it 
Fix a connection $A_0\in\Aa^{1,p}(P)$ and a constant $c_0>0$.
Let $c(A_0)$ be the constant of Lemma~\ref{le:laplace2} 
with $A$ replaced by $A_0$.  
Then there exists a constant $c_1=c_1(A_0,c_0)$
such that the following holds.
If $A\in\Aa^{1,p}(P)$ such that
$$
     \left\|A-A_0\right\|_{W^{1,p}}\le c_0
$$
and $\xi\in W^{2,p}(X,\g_P)$ 
is a solution of the equation
$$
      d_{A_0}^*d_{A_0}\xi = d_{A_0}^*(A_0-A)
$$
such that 
\begin{equation}\label{eq:step0}
     \left\|\xi\right\|_{W^{2,p}}
     \le c(A_0)\left\|d_{A_0}^*(A-A_0)\right\|_{L^p},\qquad
     \left\|\xi\right\|_{W^{1,q}}
     \le c(A_0)\left\|A-A_0\right\|_{L^q},
\end{equation}
then $g:=\exp(\xi)$ and $A_1:=g^*A$ satisfy}
\begin{equation}\label{eq:step1}
     \left\|A_1-A_0\right\|_{L^q}
     + \left\|d_{A_0}^*(A_1-A_0)\right\|_{L^p}
     \le c_1\left\|A-A_0\right\|_{L^q},
\end{equation}
\begin{equation}\label{eq:step1'}
     \left\|A_1-A_0\right\|_{W^{1,p}}
     \le c_1\left\|A-A_0\right\|_{W^{1,p}}.
\end{equation}

\medskip
\noindent
By Lemma~\ref{le:laplace1} and~(\ref{eq:step0}), 
we have 
$$
    \left\|\xi\right\|_{W^{2,p}}
    \le c_2\left\|A-A_0\right\|_{W^{1,p}}
    \le c_0c_2.
$$
for some constant $c_2=c_2(A_0)>0$.
Now let $c_3$ be the constant of
Lemma~\ref{le:quadratic}, with 
$A$ replaced by $A_0$ and 
$c_0$ replaced by $c_0c_2$. 
Then 
\begin{eqnarray*}
     \left\|A_1-A\right\|_{W^{1,p}}
&= &
     \left\|g^*A-A\right\|_{W^{1,p}}  \\
&\le &
     c_3\left(1+\left\|A-\hat A\right\|_{W^{1,p}}\right)
     \left\|\xi\right\|_{W^{2,p}} \\
&\le &
     c_2c_3\left(1+ c_0 + \left\|A_0-\hat A\right\|_{W^{1,p}}\right)
     \left\|A-A_0\right\|_{W^{1,p}}  \\
&\le &
     c_4\left\|A-A_0\right\|_{W^{1,p}}
\end{eqnarray*}
for some constant $c_4=c_4(A_0,c_0)$. 
This proves~(\ref{eq:step1'}).
Similarly,
$$
     \left\|A_1-A\right\|_{L^q}
     \le c_5\left\|A-A_0\right\|_{L^q}.
$$
for some constant $c_5=c_5(A_0,c_0)$.
 
Now consider the identity
\begin{eqnarray*}
     d_{A_0}^*(A_1-A_0)
&= &
     d_{A_0}^*(g^*A-A_0)  \\
&= &
     d_{A_0}^*(g^*A-A) 
     + d_{A_0}^*(A-A_0) \\
&= &
     d_{A_0}^*(g^*A-A-d_{A_0}\xi) \\
&= &
     d_{A_0}^*(g^*A-A-d_A\xi) 
     + d_{A_0}^*[(A-A_0)\wedge\xi].
\end{eqnarray*}
We have 
$$
     d_{A_0}^*[(A-A_0)\wedge\xi]
     = [d_{A_0}^*(A-A_0)\wedge\xi] 
       + *[d_{A_0}\xi\wedge *(A-A_0)].
$$
Let $r:=qp/(q-p)$ so that 
$$
     \frac{1}{q}+\frac{1}{r}=\frac{1}{p}.
$$ 
If $p<n$ then $r<np/(n-p)$
and hence there is a Sobolev embedding
$W^{2,p}\INTO W^{1,r}$. 
For $p\ge n$ such an embedding exists 
as well.  Hence, in either case,
\begin{eqnarray*}
&&
     \left\|d_{A_0}^*[(A-A_0)\wedge\xi]\right\|_{L^p}  \\
&&\qquad\le 
     c_6\biggl(\left\|d_{A_0}^*(A-A_0)\right\|_{L^p}
        \left\|\xi\right\|_{L^\infty}
     + \left\|A-A_0\right\|_{L^q}
        \left\|d_{A_0}\xi\right\|_{L^r}\biggr)  \\
&&\qquad\le 
    c_7\biggl(\left\|d_{A_0}^*(A-A_0)\right\|_{L^p}
          \left\|\xi\right\|_{W^{1,q}}
     + \left\|A-A_0\right\|_{L^q}
          \left\|\xi\right\|_{W^{2,p}}\biggr)  \\
&&\qquad\le 
    c_8\left\|A-A_0\right\|_{W^{1,p}}
       \left\|A-A_0\right\|_{L^q}.
\end{eqnarray*}
Moreover, it follows from Lemma~\ref{le:quadratic} that
\begin{eqnarray*}
     \left\|d_{A_0}^*(g^*A-A-d_A\xi)\right\|_{L^p} 
&\le & 
     c_9\left\|g^*A-A-d_A\xi\right\|_{W^{1,p}}   \\
&\le & 
     c_{10}\left(1+\left\|A-\hat A\right\|_{W^{1,p}}\right)
     \left\|\xi\right\|_{W^{2,p}}
     \left\|\xi\right\|_{W^{1,q}}  \\
&\le & 
     c_{11}\left\|A-A_0\right\|_{W^{1,p}}
           \left\|A-A_0\right\|_{L^q}.
\end{eqnarray*}
These two estimates imply
$$
     \left\|d_{A_0}^*(A_1-A_0)\right\|_{L^p}
     \le  c_0(c_8+c_{11})\left\|A-A_0\right\|_{L^q}
$$
and this proves~(\ref{eq:step1}).

\medskip
\noindent{\bf Step~2}
{\it 
Define the sequence $A_2,A_3,\dots$ inductively by
$$
      A_{\nu+1} = g_\nu^*A_\nu,\qquad 
      g_\nu = \exp(\xi_\nu),
$$
where $\xi_\nu\in W^{2,p}(X,\g_P)$ is chosen such that
$$
      {d_{A_0}}^*d_{A_0}\xi_\nu = {d_{A_0}}^*(A_0-A_\nu)
$$
and, with $c=c(A_0)$, 
\begin{equation}\label{eq:step2}
     \left\|\xi_\nu\right\|_{W^{2,p}}
     \le c
     \left\|d_{A_0}^*(A_\nu-A_0)\right\|_{L^p},\qquad
     \left\|\xi_\nu\right\|_{W^{1,q}}
     \le c
     \left\|A_\nu-A_0\right\|_{L^q}.
\end{equation}
There exist constants $\delta>0$ and $c_{12}>0$
such that the following holds. If 
$$
    \left\|A-A_0\right\|_{L^q}<\delta
$$
then, for every $\nu\ge 1$,}
\begin{equation}\label{eq:conv}
     \left\|d_{A_0}^*(A_\nu-A_0)\right\|_{L^p}
     \le 2^{1-\nu}c_1\left\|A-A_0\right\|_{L^q},
\end{equation}
\begin{equation}\label{eq:bound}
     \left\|A_\nu-A_0\right\|_{W^{1,p}}
     \le 2c_0c_1,
\end{equation}
\begin{equation}\label{eq:quad}
      \left\|d_{A_0}^*(A_{\nu+1}-A_0)\right\|_{L^p} 
      \le c_{12}
      \left\|A_\nu-A_0\right\|_{L^q}
      \left\|{d_{A_0}}^*(A_\nu-A_0)\right\|_{L^p}.
\end{equation}

\medskip
\noindent
For $\nu=1$ the inequalities~(\ref{eq:conv})
and~(\ref{eq:bound}) were established in Step~1.
Let $\nu\ge 1$ and assume, by induction, 
that~(\ref{eq:conv}) and~(\ref{eq:bound}) 
have been established with 
$\nu$ replaced by $j\in\{1,\dots,\nu\}$.
We prove first that~(\ref{eq:quad}) holds under
these assumptions. 
As in the proof of Step~1, we have
$$
     d_{A_0}^*(A_{\nu+1}-A_0)
     = d_{A_0}^*(g_\nu^*A_\nu-A_\nu-d_{A_\nu}\xi_\nu) 
     + d_{A_0}^*[(A_\nu-A_0)\wedge\xi_\nu]
$$
and
\begin{eqnarray*}
&&
     \left\|d_{A_0}^*[(A_\nu-A_0)\wedge\xi_\nu]\right\|_{L^p}  \\
&&\qquad\le 
     c_{13}\biggl(\left\|d_{A_0}^*(A_\nu-A_0)\right\|_{L^p}
          \left\|\xi_\nu\right\|_{W^{1,q}}
          + \left\|A_\nu-A_0\right\|_{L^q}
          \left\|\xi_\nu\right\|_{W^{2,p}}\biggr)  \\
&&\qquad\le 
     c_{14}\left\|A_\nu-A_0\right\|_{L^q}
     \left\|d_{A_0}^*(A_\nu-A_0)\right\|_{L^p}.
\end{eqnarray*}
It follows from~(\ref{eq:step2}),
the induction hypothesis, and
Lemma~\ref{le:laplace1}, that
\begin{eqnarray}\label{eq:xij}
     \left\|\xi_j\right\|_{W^{2,p}}
&\le &
     c(A_0)\left\|d_{A_0}^*(A_j-A_0)\right\|_{L^p}  
     \nonumber \\
&\le &
     c'(A_0)\left\|A_j-A_0\right\|_{W^{1,p}}  \\
&\le &
     2c_0c_1c'(A_0)
     \nonumber
\end{eqnarray}
for $j=1,\dots,\nu$. 
Now we can apply Lemma~\ref{le:quadratic},
with $c_0$ replaced by $2c_0c_1c'(A_0)$, 
to obtain, for some positive constant 
$c_{15}=c_{15}(A_0,c_0)$,
\begin{eqnarray*}
&&
     \left\|d_{A_0}^*({g_\nu}^*A_\nu-A_\nu-d_{A_\nu}\xi_\nu)\right\|_{L^p}
     \\
&&\qquad\le 
     c_{15}\biggl(1+\left\|A_\nu-\hat A\right\|_{W^{1,p}}\biggr)
     \left\|\xi\right\|_{W^{1,q}} \left\|\xi\right\|_{W^{2,p}} \\
&&\qquad\le
     c_{16}
     \left\|A_\nu-A_0\right\|_{L^q}
     \left\|d_{A_0}^*(A-A_0)\right\|_{L^p},
\end{eqnarray*}
where 
$
     c_{16}
     := c(A_0)^2c_{15}
       (1+2c_0c_1+\|A_0-\hat A\|_{W^{1,p}}).
$
Hence~(\ref{eq:quad}) holds with 
$c_{12}:=c_{14}+c_{16}$. 

Now we prove that~(\ref{eq:bound})
holds with $\nu$ replaced by $\nu+1$. 
By~(\ref{eq:xij}), the section $\xi=\xi_j$
satisfies the hypotheses of Lemma~\ref{le:quadratic}
for $j=1,\dots,\nu$, with $c_0$ replaced
by $2c_0c_1c'(A_0)$.  Hence 
\begin{eqnarray*}
     \left\|A_{j+1}-A_j\right\|_{W^{1,p}}
&\le &
     c_{15}
     \biggl(1+\left\|A_j-\hat A\right\|_{W^{1,p}}\biggr)
     \left\|\xi_j\right\|_{W^{2,p}}  \\
&\le &
     c(A_0)c_{15}
     \biggl(1+\left\|A_j-\hat A\right\|_{W^{1,p}}\biggr)
     \left\|d_{A_0}^*(A_j-A_0)\right\|_{L^p}
\end{eqnarray*}
for $j=1,\dots,\nu$. Hence, by~(\ref{eq:conv}), 
\begin{equation}\label{eq:Aj}
     \left\|A_{j+1}-A_j\right\|_{W^{1,p}}
     \le c_{17}2^{-j}\left\|A-A_0\right\|_{L^q}.
\end{equation}
for $j=1,\dots,\nu$, where
$
     c_{17}:=2c(A_0)c_{15}(1+2c_0c_1+\|A_0-\hat A\|_{W^{1,p}})
             c_1.
$
If 
$$
     c_{17}\left\|A-A_0\right\|_{L^q}\le c_0c_1
$$
then
\begin{eqnarray*}
     \left\|A_{\nu+1}-A_0\right\|_{W^{1,p}}
&\le &
     \sum_{j=1}^{\nu}\left\|A_{j+1}-A_j\right\|_{W^{1,p}}
        + \left\|A_1-A_0\right\|_{W^{1,p}} \\
&\le &
     c_{17}\left\|A-A_0\right\|_{L^q}
     + c_1\left\|A-A_0\right\|_{W^{1,p}} \\
&\le &
     2c_0c_1.
\end{eqnarray*}
This proves~(\ref{eq:bound}) with $\nu$
replaced by $\nu+1$.  

Now we shall use~(\ref{eq:quad}) and the 
induction hypothesis to prove that~(\ref{eq:conv})
holds with $\nu$ replaced by $\nu+1$. 
Since~(\ref{eq:Aj}) holds for $j=1,\dots,\nu-1$, 
we have
\begin{eqnarray}\label{eq:Aq}
     \left\|A_\nu-A_0\right\|_{L^q}
&\le &
     \sum_{j=1}^{\nu-1}\left\|A_{j+1}-A_j\right\|_{L^q}
        + \left\|A_1-A_0\right\|_{L^q} 
     \nonumber \\
&\le &
     c_{18}\sum_{j=1}^{\nu-1}\left\|A_{j+1}-A_j\right\|_{W^{1,p}}
        + \left\|A_1-A_0\right\|_{L^q} \\
&\le &
     \left(c_{17}c_{18}+c_1\right)
     \left\|A-A_0\right\|_{L^q}.
     \nonumber
\end{eqnarray}
Here $c_{18}$ is the constant in the Sobolev
embedding $W^{1,p}\INTO L^q$. If 
$$
     c_{12}(c_{17}c_{18}+c_1)
     \left\|A-A_0\right\|_{L^q}
     \le \frac{1}{2}
$$
then, by~(\ref{eq:quad}),
\begin{eqnarray*}
     \left\|d_{A_0}^*(A_{\nu+1}-A_0)\right\|_{L^p} 
&\le &
     c_{12}\left\|A_\nu-A_0\right\|_{L^q} 
     \left\|d_{A_0}^*(A_\nu-A_0)\right\|_{L^p}  \\
&\le &
     c_{12}(c_{17}c_{18}+c_1)
     \left\|A-A_0\right\|_{L^q} 
     \left\|d_{A_0}^*(A_\nu-A_0)\right\|_{L^p}  \\
&\le &
     \frac{1}{2}\left\|d_{A_0}^*(A_\nu-A_0)\right\|_{L^p}.
\end{eqnarray*}
This proves~(\ref{eq:conv}) with $\nu$ replaced by $\nu+1$.

\smallbreak

\medskip
\noindent{\bf Step~3}
{\it We prove the theorem.}

\medskip
\noindent
By~(\ref{eq:Aj}), the sequence
$A_\nu$ converges strongly in the $W^{1,p}$-norm
and the limit connection 
$$
      A_\infty := \lim_{\nu\to\infty} A_\nu \in\Aa^{1,p}(P)
$$
satisfies
$$
     \left\|A_\infty-A_0\right\|_{W^{1,p}}
     \le \sum_{j=0}^\infty
         \left\|A_{j+1}-A_j\right\|_{W^{1,p}}
     \le (c_1+c_{17}c_{18})\left\|A-A_0\right\|_{W^{1,p}}.
$$
Moreover, by~(\ref{eq:Aq}),
$$
     \left\|A_\infty-A_0\right\|_{L^q}
     \le (c_1+c_{17}c_{18})\left\|A-A_0\right\|_{L^q}
$$
and, by~(\ref{eq:conv}), 
$$
     d_{A_0}^*(A_\infty-A_0)
     = \lim_{\nu\to\infty}d_{A_0}^*(A_\nu-A_0) = 0.
$$
Write
$$
      A_\nu = h_\nu^*A,\qquad h_\nu = g g_1 g_2\cdots g_\nu.
$$
Consider the identity
\begin{equation}\label{eq:h}
      dh_\nu = h_\nu A_\nu - Ah_\nu
\end{equation}
in a local frame.  
The right hand side of~(\ref{eq:h}) is bounded
in $L^q$ and hence $h_\nu$ is bounded
in $W^{1,q}$. Now the product inequality
of Lemma~\ref{le:product} with $r=p$
shows that $h_\nu A_\nu-Ah_\nu$ is bounded
in $W^{1,p}$ and, by~(\ref{eq:h}), 
$h_\nu$ is bounded in $W^{2,p}$.  
Hence $h_\nu$ has a subsequence,
still denoted by $h_\nu$, which converges
in the $W^{1,q}$-norm. Since $A_\nu$
converges in the $W^{1,p}$-norm
it follows from Lemma~\ref{le:product}
with $r=p$ that $h_\nu A_\nu-Ah_\nu$
converges in the $W^{1,p}$-norm.
By~(\ref{eq:h}), $h_\nu$ converges in the 
$W^{2,p}$-norm.  The limit 
$$
     h_\infty:=\lim_{\nu\to\infty}h_\nu\in\Gg^{2,p}(P)
$$
satisfies 
$$
     A_\infty = \lim_{\nu\to\infty}h_\nu^*A = h_\infty^*A.
$$
This proves the theorem.
\end{proof}



\begin{thebibliography}{99}
\small

\bibitem{AUDIN} M.~Audin,
{\it The Topology of Torus Actions on Symplectic Manifolds},
Birk\-h\"auser, Basel, 1991.

\bibitem{BDW} A.~Bertram, G.~Daskalopoulos, R.~Wentworth, 
Gromov invariants for holomorphic maps 
from Riemann surfaces into Grassmannians, 
{\it J.~Amer.~Math.~Soc.} {\bf 9} (1996), 529--571.

\bibitem{BRADLOW1} S.~Bradlow, 
Vortices in holomorphic line bundles over closed K\"ahler manifolds, 
{\it Comm.~Math.~Phys.} {\bf 135} (1990), 1--17.

\bibitem{BRADLOW2} S.~Bradlow,
Special metrics and stability for holomorphic bundles with 
global sections, {\it J.~Diff.~Geom.} {\bf 33} (1991), 169--213.

\bibitem{BDW1} S.~Bradlow, G.~Daskalopoulos, R.~Wentworth,
Birational equivalence of vortex moduli, 
{\it Topology} {\bf 35} (1996), 731--748.

\bibitem{CGS} K.~Cieliebak, A.R.~Gaio, D.A.~Salamon,
J-holomorphic curves, moment maps, and invariants
of Hamiltonian group actions,
{\it IMRN} {\bf 10} (2000), 831--882.

\bibitem{CMS} K.~Cieliebak, I.~Mundet i Riera, D.A.~Salamon,
Equivariant moduli problems and the Euler class,
Preprint, ETH-Z\"urich, October 2001.
To appear in {\it Topology}.

\bibitem{DON1} S.K.~Donaldson,
Symplectic submanifolds and almost complex ge\-o\-me\-try,
{\it J.~Diff.~Geom.} {\bf 44} (1996), 666--705.

\bibitem{DON2} S.K.~Donaldson,
Lefschetz fibrations on symplectic manifolds,
Preprint, 1999, to appear in {\it J. Diff. Geo.}

\bibitem{DONS} S.K.~Donaldson, I.~Smith,
Lefschetz pencils and the canonical class for symplectic
4-manifolds, Preprint, December 2000,
{\it math.~SG/0012067}.

\bibitem{DS}  S.~Dostoglou, D.A.~Salamon,
Self-dual instantons and holomorphic curves,
{\it Annals of Math.} {\bf 139} (1994), 581--640.

\bibitem{EG} Y.~Eliashberg, M.~Gromov, 
Convex symplectic manifolds, 
in: {\it Several Complex Variables and Complex Geometry}
(Santa Cruz, CA, 1989), Part 2, 
{\it Proc.~Sympos.~Pure Math.} \textbf{52}, 
Amer.~Math.~Soc., Providence, RI (1991), 135--162.

\bibitem{FO} K.~Fukaya, K.~Ono,
Arnold conjecture and Gromov--Witten invariants for 
general symplectic manifolds, Preprint, February 1996,
Summary in: {\it Proceddings of the Taniguchi Conference 
on Mathematics}, edited by T.~Sunada and M.~Maruyama,
{\it Advanced Studies in Pure Math.} {\bf 31},
Math.~Society of Japan (2001), 75--91. 

\bibitem{GAIO} A.R.~Gaio,
J-holomorphic curves and moment maps,
PhD thesis, University of Warwick, 1999.

\bibitem{GaSa} A.R.~Gaio, D.A.~Salamon,
J-holomorphic curves, moment maps, and adiabatic limits,
Preprint, June 2001.

\bibitem{GP} O.~Garcia-Prada, 
A direct existence proof for the vortex equations over a 
compact Riemann surface, 
{\it Bulletin L.M.S.} {\bf 26} (1994), 88--96.

\bibitem{GRO} M.~Gromov,
Pseudo holomorphic curves in symplectic manifolds,
{\it Invent.~Math.} {\bf 82} (1985), 307--347.

\bibitem{KW} J.~Kazdan, F.~Warner, 
Curvature functions for compact 2-man\-i\-folds, 
{\it Ann.~Math.} {\bf 99} (1978), 14--47.

\bibitem{KIRWAN} F.~Kirwan,
{\it Cohomology of Quotients in Symplectic and 
Algebraic Geometry}, Princeton University Press, 1984.
 
\bibitem{LRR} A.M.~Li, J.~Robbin, Y.~Ruan,
{\it Virtual Muduli Cycles and Gromov--Witten Invariants},
in preparation. 

\bibitem{LL} T.J.~Li, A.~Liu,
General wall-crossing formula,
{\it Math.~Res.~Letters} {\bf 2} (1995), 797--810.

\bibitem{LT}  J.~Li, G.~Tian,
Virtual moduli cycles and Gromov--Witten invariants
of algebraic varieties,
{\it J.~Amer.~Math.~Soc.} {\bf 11} (1998), 119--174. 

\bibitem{MS1} D.~McDuff, D.A.~Salamon, 
{\it $J$-holomorphic Curves and Quantum Cohomology},
University Lecture Series {\bf 6}, 
Amer.~Math.~Soc., Providence, RI, 1994.  

\bibitem{MS2} D.~McDuff, D.A.~Salamon, 
{\it Introduction to Symplectic Topology},
Oxford University Press, 1995.  

\bibitem{MUNDET} I.~Mundet i Riera,
Yang-Mills-Higgs theory for symplectic fibrations,
PhD thesis, Madrid, April 1999. 

\bibitem{MUNDET2} I.~Mundet i Riera,
Hamiltonian Gromov--Witten invariants,
Pre\-print {\it math.~SG/0002121}.
To appear in {\it Topology}.

\bibitem{OT} C.~Okonek, A.~Teleman,
Gauge theoretical equivariant Gromov--Witten invariants
and the full Seiberg--Witten invariants of ruled surfaces.
Preprint {\it math.~SG/0102119}.

\bibitem{OO} K.~Ono, H.~Ohta,
Notes on symplectic $4$-manifolds with $b^+=1$, I and II, 
{\it Internat.~J.~of Math.} {\bf 7} (1996), 755--770.

\bibitem{RUAN} Y.~Ruan,
Virtual neighborhoods and pseudoholomorphic curves,  
in: {\it Topics in Symplectic $4$-Manifolds} (Irvine CA 1996), 
Internat.~Press, Cambridge, MA, 1998, 101--116.

\bibitem{SAL1} D.A.~Salamon,
Seiberg--Witten invariants of mapping tori, symplectic fixed points,
and Lefschetz numbers, 
{\it Turkish Journal of Mathematics} {\bf 23} (1999), 117--143.

\bibitem{SAL2} D.A.~Salamon,
{\it Spin Geometry and Seiberg--Witten Invariants},
in preparation. 

\bibitem{UHL} K.~Uhlenbeck,
Connections with $L^p$ bounds on the curvature,
{\it Commun.~Math.~Phys.} {\bf 83} (1982), 31--42.

\bibitem{TAUBES1} C.H.~Taubes,
The Seiberg-Witten and the Gromov invariants,
{\it Math.~Res.~Letters} {\bf 2} (1995), 221--238.

\bibitem{TAUBES2} C.H.~Taubes,
$\SW\IMP\Gr$, from the Seiberg-Witten equations to 
pseudoholomorphic curves, 
{\it J.~Amer.~Math.~Soc.} {\bf 9} (1996), 845--918.

\bibitem{TAUBES3} C.H.~Taubes,
{\it Seiberg-Witten and Gromov Invariants for Symplectic 
Four-Manifolds}, International Press, 2000.

\bibitem{T} M.~Thaddeus,
Stable pairs, linear systems, and the Verlinde algebra,
{\it Inv.~Math.} {\bf 117} (1994), 317--353. 

\bibitem{WEHRHEIM} K.~Wehrheim,
{\it Uhlenbeck Compactness: An Exposition}, 
in preparation.

\end{thebibliography}
\end{document}